\documentclass[12 pt]{amsart}
\usepackage{amsthm, amsfonts, amssymb, color}

\usepackage{mathrsfs}
\usepackage{amsmath}
\usepackage{amstext, amsxtra}
\usepackage{txfonts}
\usepackage[colorlinks, linkcolor=black, citecolor=blue, pagebackref, hypertexnames=false]{hyperref}
\usepackage{comment}
\usepackage{soul}

\allowdisplaybreaks
\newcommand{\bl}{\color{blue}}
\usepackage{geometry}
\newtheorem{theorem}{Theorem}[section]
\newtheorem{proposition}[theorem]{Proposition}

\newtheorem{lemma}[theorem]{Lemma}
\newtheorem{definition}[theorem]{Definition}
\newtheorem{example}[theorem]{Example}
\newtheorem{assumption}{Assumption}[section]
\newtheorem{remark}[theorem]{Remark}

\newcommand{\s}{ L}
\newcommand{\eq}{\begin{equation}}
\newcommand{\eeq}{\end{equation}}
\newcommand{\eqa}{\begin{eqnarray}}
\newcommand{\eeqa}{\end{eqnarray}}

\newcommand\N{\mathbb{N}}
\newcommand\R{\mathbb{R}}

\newcommand\F{\mathcal{F}}

\numberwithin{equation}{section}

\title   [] { On The large Time Asymptotics of  Schr\"odinger type equations with General Data}

\author{   Avy Soffer }

\address{Department of Mathematics\\
Rutgers University\\
110 Frelinghuysen Rd.\\
Piscataway, NJ, 08854, USA}

\email{soffer@math.rutgers.edu}

\author{Xiaoxu Wu}
\address{Mathematical Sciences Institute\\ Australia National University\\ 
Canberra 2601, Australia.}
 \email{xiaoxu.wu@anu.edu.au}

\thanks{
A.Soffer is supported in part by Simons Foundation Grant number 851844
}

\begin{document}

\begin{abstract}
For the Schr\"odinger equation with a general interaction term, which may be linear or nonlinear,
time dependent and including charge transfer potentials, we prove the global solutions are asymptotically given by a free wave and a  weakly localized part.
The proof is based on constructing in an adapted way the Free Channel Wave Operator, and further tools from the recent works \cite{Liu-Sof1,Liu-Sof2,SW2020}.
This work generalizes the results of the first part of  \cite{Liu-Sof1,Liu-Sof2} to arbitrary dimension, and non-radial data.
\end {abstract}
\maketitle
\section{Introduction}
The analysis of dispersive wave equations and systems is of critical importance in the study of evolution equations in Physics and Geometry.
It is well known that the asymptotic solutions of such equations, if they exist, show a dizzying zoo of possible solutions.
Besides the "free wave", which corresponds to a solution of the equation without interaction terms, a multitude of other solutions may appear.
Such solutions are localized around possibly moving center of mass. They include nonlinear bound states, solitons, breathers, hedgehogs, vortices etc...
The analysis of such equations is usually done on a case by case basis, due to this complexity. \cite{Sof}
A natural question then follows: is it true that in general, solutions of dispersive equations converge in appropriate norm ($L^2$ or $\mathcal{H}^1_x$) to a free wave and independently moving localized parts{ (localized in space)}?
In fact this is precisely the statement of Asymptotic Completeness in the case of $N$-body Scattering. In this case the possible outgoing clusters are clearly identified, as bound states of subsystems.
But when the interaction term includes time dependent potentials (even localized in space) and more general nonlinear terms, we do not have an {\emph{a priory}} knowledge of the possible asymptotic states.

In the case of time independent interaction terms, one can use spectral theory. The scattering states evolve from the continuous spectrum, and the localized part is formed by the point spectrum.
Once the interaction is time dependent/nonlinear that is not possible.
In fact, there are no general scattering results for localized time dependent potentials.
  The exceptions are charge transfer Hamiltonians \cite{Yaj1, Gr1,Wul,Per, RSS},
decaying in time potentials and small potentials \cite{ How, RS}, time periodic potentials \cite{Yaj2, How} and random (in time) potentials \cite{BeSof}.
See also \cite{BeSof1,BeSof2}.
 For potentials with asymptotic energy distribution more could be done \cite{SS1988}.
A recent progress for more general localized potentials without smallness assumptions is obtained in \cite{SW2020}.

Turning to the nonlinear case, Tao \cite{T1,T2,T3} has shown that the asymptotic decomposition holds for NLS with inter-critical nonlinearities, in 3 or higher space dimensions, in the case of radial initial data.
In particular, in a sufficiently high dimension, and with an interaction that is a sum of a smooth compactly supported potential and a repulsive nonlinearity, Tao was able to show that the localized part is smooth
and localized.
In other cases, Tao showed the localized part is only weakly localized and smooth.
Tao's work uses direct estimates of the incoming and outgoing parts of the solution to control the nonlinear part, via Duhamel representation.
In a {certain} sense, it is in the spirit of Enss' work. See also \cite{Rod-Tao}.

In contrast, the new approach of Liu-Soffer \cite{Liu-Sof1,Liu-Sof2} is based on proving {\emph{a priory}} estimates on the full dynamics, which hold in a suitably localized regions of the extended phase-space.
In this way it was possible to show the asymptotic decomposition for general localized interactions, including time and space dependent ones, which are localized. Radial initial data is assumed.
More detailed information is obtained on the localized part of the solution.
 Besides being smooth, its expanding part (if it exists) can grow at most like $|x|\leq C\sqrt t$ (for $t\geq1$ and some constant $C>0$), and furthermore, is concentrated
in a thin set of the extended phase-space.
The free part of the solution concentrates on the \emph{propagation set} where $x=vt, v=2P,$ and $P$ being the dual to the space variable, the momentum, is given by the operator $-i\nabla_x.$
The weakly localized part is found to be localized in the regions where $$|x|/t^{\alpha}\sim 1  \quad\textbf{and} \quad |P|\sim  t^{-\alpha}, \quad\forall ~0<\alpha \leq 1/2.$$
It therefore shows that the spreading part follows a self similar pattern.
The method of proof is based on three main parts: first, construct the Free Channel Wave Operator. Then prove localization of the remainder localized part, and use it to prove the smoothness of the localized part.
Finally, by using further propagation estimates which are adapted to localized solutions, Liu and Soffer proved the concentration on thin sets of the phase-space corresponding to self similar solutions.
It should be emphasized that the spreading localized solutions, if they exist, were shown to have a non-small nuclei part around the origin. This is true for both the results of Tao \cite{T1,T2,T3} and Liu-Soffer \cite{Liu-Sof1,Liu-Sof2}.
Therefore, these are not purely self-similar solutions, as appear in the special cases of critical nonlinearities. See e.g. \cite{S2021, DKM2}.

We will follow here this point of view.
The key tool from scattering theory that is used to study multichannel scattering is the notion of \emph{channel wave operator}, which we denote by
\begin{equation}\label{omegaa}
\Omega_a^* \psi(0)\equiv s\text{-}\lim\limits_{t \to \infty} e^{iH_a t}U(t,0)\psi(0),\quad \psi(0)\in \mathcal{H},
\end{equation}
where $U(t,0)$ denotes the evolution operator of the full dynamics on a Hilbert space $\mathcal{H}$ and $H_a$ stands for the generator of one possible asymptotic dynamics, for example, $H_a$ can be $-\Delta_x$. Here, \eqref{omegaa} is well-defined whenever the strong limit exists.

Here the limit is taken in the strong sense in $L^2.$
Note that since $U(t,0)$ is nonlinear in general, then so is the wave operator $\Omega_a^*.$
$U(t,0)\psi(0)$ is the solution of the dispersive equation with initial data $\psi(0)$ and dynamics (linear or nonlinear) $U(t)\equiv U(t,0)$ generated by a Hamiltonian $H(t).$
The asymptotic dynamics is generated by a Hamiltonian $H_a$ for a given channel denoted by $a.$
In this work we will only construct the{ \bf free channel wave operator}, where $H_a=-\Delta_x.$

A crucial observation is that one can modify the definition of the Channel wave operators to
\begin{equation}
\Omega_a^* \psi(0)\equiv s\text{-}\lim\limits_{t \to \infty} e^{iH_a t}J_a U(t,0)\psi(0),\quad \psi(0)\in \mathcal{H},
\end{equation}
provided 
\begin{equation}\label{wlim}
w\text{-}\lim\limits_{t \to \infty}e^{iH_a t}(1-J_a)  U(t,0)\psi(0)=0.
\end{equation}
Here $J_a$ denotes some operator with norm $1$ and {we take $J_a=\mathcal F_c(\frac{|x-2Pt|}{t^\alpha}\leq 1)$ for the free channel wave operator. See~\eqref{def: F}-\eqref{eq: Fj<} for the definition of $\mathcal F_c$. See also \cite {SS1987}}. This construction can be easily generalized to the case {where} the {\it asymptotic} dynamics is nonlinear.
In practice, we should choose $J_a$ to be a member of a partition of unity which is supported on the extended phase space where the solution is expected to converge; to be useful, it should also be decaying (in some vague sense) on the support of the interaction that couples the channel $a$ to the rest of the solution.

Now, to prove that the limit exists we use Cook's method. For this, we need to show the integrability of the derivative w.r.t. time of the vector $e^{iH_a t}J_a U(t,0)\psi(0)$ in $\mathcal{H}$. Taking the derivative (w.r.t. time) gives two types of terms: with $\partial_t[U(t,0)]=(-i)H(t)U(t,0)$,
\begin{align}
\partial_t[e^{iH_a t}J_a U(t,0)]\psi(0)=&e^{iH_a t}D_{H_a}(J_a) U(t,0)\psi(0)-e^{iH_a t}iJ_a(H(t)-H_a) U(t,0)\psi(0).
\end{align}
Here the operator $D_H(B)$ denotes 
\begin{align}
D_H (B)\equiv i[H,B]+\frac{\partial B}{\partial t}.
\end{align}
For example, with $\psi(t)\equiv U(t,0)\psi(0)$, when $H_a=H_0\equiv-\Delta_x$ and the interaction $\mathcal N(x,t,|\psi(t)|)=H(t)-(-\Delta_x)$, $\partial_t[ e^{iH_a t}J_a U(t,0)]\psi(0)$ reads
\begin{equation}
\partial_t[ e^{iH_0 t}J_a U(t,0)]\psi(0)=e^{iH_0 t}D_{H_0}(J_a) \psi(t)-ie^{iH_0 t}J_a \mathcal N(x,t,|\psi(t)|)\psi(t).
\end{equation}

By choosing
 $$
 J_a= F(\frac{|x|}{t^{\alpha}}\geq 1),
 $$
where $F$ denotes a smooth-cut off function or a smooth characteristic function, it is easy to see that such $J_a$ satisfies our requirement, as on its support the interaction term vanishes like $t^{-m\alpha}$
for a localized interaction vanishing like $|x|^{-m}$ at infinity. Furthermore, it is not hard to prove that the identity \eqref{wlim} holds true by using duality and the dispersive estimate of free flows. However, the Heisenberg Derivative part coming from $D_H$ is not necessarily integrable in time, under the full dynamics.
The solution can have a part that stays on the boundary of the support of $F$, or revisit it for infinitely many times.
To resolve this problem, as was done in the $N$-body case \cite{SS1987} and in the general nonlinear case \cite{Liu-Sof1,Liu-Sof2},
we further microlocalize the partition of unity, such that on the boundary, the solution can be shown to decay sufficently fast in time $t$ (by propagation estimates).
In \cite{SS1987} these boundaries are cones in the configuration space, and then one needs to microlocalize the momentum to point either out or into the cone.
In \cite{Liu-Sof1,Liu-Sof2}  one microlocalizes the partition $F$ by localizing on the incoming/outgoing parts of the solution.
This microlocalization needs to be done in a way that allows proving \emph{propagation estimates} there \cite{SS1987,D-Ger}.
It should be clear by now, that this method is tied to a distinguished point in space, and requires the interaction term to be localized around it.
The function $F$ can only annihilate a localized term, and the notion of incoming and outgoing is tied to the choice of origin.
Therefore, in order to go to the general initial data case, we need a more general type of constructions.
This is the content of this work.

The key new construction is a free channel wave operator, with a different type of localization in the phase space.
This localization is constructed by projecting in the phase-space on a neighborhood of the thin propagation set in the extended phase space.
As the free wave concentrates where $x=2Pt$, we use the projection, with $H_a=H_0\equiv-\Delta_x$,
$$
J_{\text{free}}\equiv J_a=\mathcal F_c(\frac{|x-2Pt|}{t^{\alpha}}\leq 1), \quad \quad \text{ for some }\alpha \in (0,1)\text{ and }t>1,
$$
where $\mathcal F_c$ denotes a smooth cut-off function. Here, the subscript of $\F_c$ stands for the "conformal multiplier" and we define $\F_c(\frac{|x-2Pt|}{t^\alpha}\leq 1)$ as an operator on $L^2$ or $\mathcal{H}$ by using the spectral theorem and the self-adjointness of $|x-2tP|$. It is a property of the free dynamics that the solution vanishes outside the support of $\F_c$
as time goes to infinity.
The fundamental {property} of this operator that we use {is} the following {equation}:
\begin{equation}\label{id: Fc}
e^{-iH_0 t}\F_c(\frac{|x|}{t^{\alpha}}\leq 1)e^{iH_0 t}=\F_c(\frac{|x-2Pt|}{t^{\alpha}}\leq 1).
\end{equation}
See Section \ref{app: phase-space} for detailed computations. 

Throughout this paper, $C$ will denote a constant and may vary from one line to another. We write $\lesssim$ or $\gtrsim$ whenever $A\leq CB$ or $CA\geq B$ for some constant $C>0$. We write $A\lesssim_a B$ or $A\gtrsim_a B$ if  $A\leq C_aB$ or $C_aA\geq B$ for some constant $C_a>0$ which depends on parameter $a$. 

A useful property of $\F_c$ is the following inequality 
\begin{align}
&\|\F_c(\frac{|x-2tP|}{t^\alpha}\leq 1) \phi\|_{L^p_x(\mathbb{R}^n)}\nonumber\\
=&\|e^{-i|x|^2/4t}\F_c(\frac{|x-2tP|}{t^\alpha}\leq 1) \phi\|_{L^p_x(\mathbb{R}^n)}\nonumber\\
\lesssim& \|P e^{-i|x|^2/4t}\F_c(\frac{|x-2tP|}{t^\alpha}\leq 1) \phi\|_{L^2_x(\mathbb{R}^n)}^a\|e^{-i|x|^2/4t}\F_c(\frac{|x-2tP|}{t^\alpha}\leq 1) \phi\|_{L^2_x(\mathbb{R}^n)}^{1-a}\nonumber\\
\lesssim &\| (1/t) |2Pt-x| \F_c(\frac{|x-2tP|}{t^\alpha}\leq 1) \phi\|_{L^2_x(\mathbb{R}^n)}^a \|\F_c(\frac{|x-2tP|}{t^\alpha}\leq 1) \phi\|_{L^2_x(\mathbb{R}^n)}^{{1-a}}\nonumber\\
\lesssim &t^{(-1+\alpha)a}\|\F_c(\frac{|x-2tP|}{t^\alpha}\leq 1) \phi\|_{L^2_x(\mathbb{R}^n)}\nonumber\\
\lesssim &t^{(-1+\alpha)a}\|\F(\frac{|x|}{t^{\alpha}}\leq 1) e^{itH_0}\phi\|_{L^2_x(\mathbb{R}^n)},
\end{align}
where $n$ is the space dimension, $|x|$ denotes the length of $x$ in $\mathbb{R}^n$ and we have used Nirenberg-Sobolev type inequality, {the unitarity of $e^{-itH_0}$ on $L^2_x(\mathbb{R}^n)$ and Eqs. 
\begin{equation}
    P e^{-i|x|^2/4t}= e^{-i|x|^2/4t}P-e^{-i|x|^2/4t}\frac{x}{2}\cdot P=\frac{e^{-i|x|^2/4t}}{2}(P-\frac{x}{2t})
\end{equation}
and~\eqref{id: Fc}}. Here, the constants $p>2$ and $a$ depend on {the dimension of the space}.
For example, in three space dimensions, $p=6$, $a=1.$
Furthermore, the Heisenberg Derivative of this operator is positive:
\begin{equation}
D_{H_0} \F_c(\frac{|x-2Pt|}{t^{\alpha}}\leq 1)= -\alpha \frac{|x-2Pt|}{t^{1+\alpha}} \F_c' \geq 0.
\end{equation}
This is due to the fact that $D_{H_0}(|x-2Pt|^2)=0$.\par
The operator $\F_c(\frac{|x-2Pt|}{t^{\alpha}}\leq 1)$ and its functions have a long history.
In fact, the operator $|x-2tP|^2$ is the multiplier that gives the conformal identity for Schr\"odinger equations.
Then $\F_c(\frac{|x-2Pt|}{t^{\alpha}}\leq 1)$ was used to prove sharp propagation estimates in \cite{Sig,SSjams,SSInvention,SSDuke, Gr, Dere}.
In a completely different way it was used in \cite {Lind-S1,Lind-S2}.
Using propagation estimates similar to \cite{SS1987}, the problem of showing the existence of the free channel wave operator, defined in terms of the above $\F_c,$ is reduced to proving the propagation estimate that follows from using $\F_c$ as a propagation observable.
Since the Heisenberg derivative is positive, it remains to verify for what interaction terms the following is true:
$$
\int_1^{\infty} \|\F_c \mathcal N(x,t,|\psi|) U(t,0)\psi(0)\|_{L^2_x(\mathbb{R}^n)} dt < \infty.
$$
See Section \ref{sec: prop} for detailed discussions.

\section{The Scattering Problem and Results}

Let $H_0:=-\Delta_x$. We consider the general class of Nonlinear Schr\"odinger type equations of the form:
\eq
\begin{cases}
i\partial_t\psi(x,t)=H_0\psi(x,t)+\mathcal{N}(x,t,\psi(x,t))\psi(x,t) \\
\psi(x,0)=\psi_0\in \mathcal{H}^a_x(\mathbb{R}^n)
\end{cases}, \quad (x,t)\in \mathbb{R}^n\times \mathbb{R}\label{SE}
\eeq
with space dimension $n\geq 1$, where $\mathcal{H}_x^a\equiv\mathcal{H}_x^a(\mathbb{R}^n), a\in [0,1],$ denotes the $L^2$ Sobolev space of order $a$. $\mathcal{N}$ is NOT assumed to be real.

\subsection{Assumptions and examples.} {The equations under consideration are well studied \cite{Str, Caz}.  
Let
\[
B \equiv C^{1}(\mathbb{R},L^{2}) + C^{1}(\mathbb{R},L^{\infty})
\]
denote the space of functions that are continuously differentiable in $t$ with values in $L^{2}+L^{\infty}$.  
In the linear case
\[
\mathcal{N}(|x|,t,|\psi|) \equiv V(x,t), 
\qquad \psi(0)\in H^{1},
\]
with $V \in B$, global existence in $H^{1}$ holds, and the solution is given by a unitary group
\[
U(t,0)\psi(0)=\psi(t),
\]
see, for example, \cite{Yos} and \cite[Thm.~X.70, Thm.~X.71]{R-SII}.

When $\mathcal{N}(\psi)$ is purely nonlinear, i.e.\ $\mathcal{N}(\psi)= f(\psi)$, the existence theory is well developed in both the energy-subcritical and critical regimes; see \cite{Caz, Bour, Tao}.  
In the supercritical setting, weak global solutions exist in the defocusing case \cite{Str}. For lower-power nonlinearities that vanish sufficiently rapidly near the origin, standard energy estimates yield global existence, even in the presence of additional potential terms \cite{Caz}.

}

We consider solutions $\psi(t)\equiv \psi(x,t)$ of system \eqref{SE} which exist globally in $t\in\mathbb{R}$ and are uniformly bounded in $\mathcal{H}^a_x$. {The solution $\psi(t)$ is understood as a weak solution.}
\begin{assumption}\label{asp: global}There exists a positive constant $C>0$ such that
\eq
E :=\sup\limits_{t\in \mathbb{R}}\|\psi(t)\|_{\mathcal{H}^a_x}\leq C<\infty\label{con: H}
\eeq
is valid for some $a\in [0,1]$. 
\end{assumption}
Let $\langle \cdot\rangle: \mathbb{R}^n\to \mathbb{R}, x\mapsto \sqrt{|x|^2+1}$, denote the Japanese brackets. We consider the interaction $\mathcal{N}(x,t,\psi(t))$, which falls into one of the following categories:
\begin{enumerate}
\item(\emph{Space localized potentials}):\,For $n\geq 1$ and some $\delta >1$,
\eq
\langle x\rangle^\delta\mathcal{N}(x,t,\psi(t))\in L^\infty_{x,t}(\mathbb{R}^{n+1}).\label{local}
\eeq
{\item (\emph{$\s^p$ potentials}):\,For $n\geq 3$, 
\eq
\mathcal N(x,t,\psi(t))\psi(t)\in \s^\infty_tL^1_x(\mathbb{R}^{n+1}).\label{Lpresult}
\eeq
Our $L^p$ potentials cover the following models, which are proved in Examples~\ref{example1} and~\ref{example2} {of} Section~\ref{sec: example}:
\begin{enumerate}
\item (\emph{Charge transfer Hamiltonians}):\label{charge} Let Assumption~\ref{asp: global} hold. When the space dimension $ n\geq 3$, the charge transfer interaction $\mathcal{N}(x,t,\psi)=\sum\limits_{j=1}^N V_j(x-tv_j,t)$, where $V_j(x,t)\in \s^\infty_t\s^2_x(\mathbb{R}^{n+1})$, $j=1,\cdots,N$, with $v_j\neq v_l$ if $j\neq l$, satisfies condition~\eqref{Lpresult}.
\item \label{Nonlinear}(\emph{Purely nonlinear {interaction terms}}): Let Assumption~\ref{asp: global} hold {with} $a=1$. When space dimension $n\geq3,$ $\mathcal{N}(x,t,\psi)=I(|\psi|)$, with $I(|\psi|)$ satisfying the estimate
\eq\label{L1: ineq}
\|I(|\psi|)\psi\|_{L^1_x(\mathbb{R}^n)}\lesssim_{\|\psi\|_{\mathcal{H}^1_x}} 1,
\eeq
satisfies condition~\eqref{Lpresult}.
\end{enumerate}}
\end{enumerate}

Here, typical examples for purely nonlinear interactions (interactions which depend only on $\psi(t)$) are polynomial nonlinearities (see Example~\ref{example3} {of} Section~\ref{sec: example}): 
\eq
I(|\psi|)=P(|\psi|), \qquad n\geq 3,
\eeq
where $P(z)$ denotes a polynomial of degree $N$ with $P(0)=0$. The main assumption to be verified in this case, is that the energy identity implies the $L^1$ condition (see~\eqref{L1: ineq}). See Example~\ref{example3} for a detailed discussion.

\begin{assumption}[Space localized potentials]\label{asp: local}Assume Eq.~\eqref{local} is valid for some $\delta>1$.
    
\end{assumption}
\begin{assumption}[$\s^p$ potentials]\label{asp: nonlocal}Assume Eq.~\eqref{Lpresult} is valid. 
    
\end{assumption}
Let $W^{k,p}_x(\mathbb{R}^n), 1\leq p\leq \infty,$ denote the $L^p$ Sobolev space of order $k$. We refer to $\mathcal{N}(x,t,\psi(t))$ as a charge-transfer interaction if $\mathcal{N}$ is linear and if there exists a positive integer $N\geq 2$ and $N$ vectors $v_j\in \mathbb{R}^n, j=1,\cdots,N,$ ($v_j\neq v_l$ if $j\neq l$) such that 
\eq\label{charge: expression}
\mathcal{N}(x,t,\psi(t))=\sum\limits_{j=1}^N V_j(x-tv_j, t),\quad 
\eeq
where $V_j(x,t), j=1,\cdots,N,$ are functions localized in $x$ variable (see Assumption~\ref{asp: charge}). 
{\begin{assumption}[Charge-transfer potentials]\label{asp: charge} There exists $\delta\geq n+1$ such that $\mathcal{N}(x,t,\psi(t))$ satisfies ~\eqref{charge: expression} with $V_j(x,t), j=1,\cdots,N,$ satisfying
\eq
\sup\limits_{t\in \mathbb{R}} \|\langle x\rangle^\delta V_j(x,t)\|_{W^{1,\infty}_x(\mathbb{R}^{n})}<\infty.\label{Vj}
\eeq
    
\end{assumption}}

{\subsection{Examples}\label{sec: example}
\begin{example}\label{example1} Let $N\in \mathbb{N}^+$ be a positive integer. When $\mathcal{N}(x,t,\psi)=\sum\limits_{j=1}^N V_j(x-tv_j,t)$ with $V_j(x,t)\in \s^\infty_t\s^2_x(\mathbb{R}^{n+1}), j=1,\cdots,N$, Assumption~\ref{asp: nonlocal} is satisfied if Assumption~\ref{asp: global} holds true. 
    
\end{example}
\begin{proof}To compute $\| \mathcal{N}(x,t,\psi(t))\psi(t)\|_{L^\infty_tL^1_x(\mathbb{R}^{n+1})}$, we find, by H\"older's inequality 
\begin{align}
    \| \mathcal{N}(x,t,\psi(t))\psi(t)\|_{L^\infty_tL^1_x(\mathbb{R}^{n+1})}\leq & \|\mathcal{N}(x,t,\psi(t)) \|_{L^\infty_tL^2_x(\mathbb{R}^{n+1})}\|\psi(t)\|_{L^\infty_tL^2_x(\mathbb{R}^{n+1})}\nonumber\\
    \lesssim_E & \sum\limits_{j=1}^N \| V_j(x,t)\|_{L^\infty_tL^2_x(\mathbb{R}^{n+1})}< \infty.
\end{align}

\end{proof}
\begin{example}\label{example2} Suppose that Assumption~\ref{asp: global} is valid with $a=1$. When space dimension $n\geq3,$ $\mathcal{N}(x,t,\psi(t))=I(|\psi(t)|)$, with $I(|\psi|)$ satisfying the estimate
\eq\label{L1: ineq: ex}
\|I(|\psi|)\psi\|_{L^1_x(\mathbb{R}^n)}\lesssim_{\|\psi\|_{\mathcal{H}^1_x}} 1,
\eeq
satisfies condition~\eqref{Lpresult}.
    
\end{example}
\begin{proof}It follows from that, by H\"older's inequality, Assumption~\eqref{asp: global} and Eq.~\eqref{L1: ineq: ex}, 
\begin{align}
    \| \mathcal{N}(x,t,\psi(t))\psi(t)\|_{L^\infty_tL^1_x(\mathbb{R}^{n+1})}
    =& \| I(\psi(t))\psi(t)\|_{L^\infty_tL^1_x(\mathbb{R}^{n+1})}\nonumber\\
    \lesssim_E & 1.
\end{align}
    
\end{proof}
\begin{example}\label{example3}Assume that Assumption~\ref{asp: global} is valid with $a=1$ and let $N\in \mathbb{N}^+$ be a positive integer. When space dimension $n\geq 3$, $\mathcal{N}(x,t,\psi(t))=P(|\psi(t)|),$ where $P(z)$ denotes a polynomial of degree $N$ with $P(0)=0$, satisfies condition~\eqref{Lpresult} provided that
\eq
\|\psi(0)\|_{L^{N+2}_x(\mathbb{R}^n)}<\infty.\label{E(0)}
\eeq
    
\end{example}
\begin{proof} The energy of this purely nonlinear system is given by
\eq
E(\psi(t)):=(\psi(t),(-\Delta_x+\tilde P(|\psi(t)|)\psi(t))_{L^2_x(\mathbb{R}^n)}, 
\eeq
where $\tilde P(z)$ is another polynomial of the same degree as $P(z)$, defined by, with $P'(k)\equiv \frac{d}{dk}[P(k)]$,
\eq
\tilde P(k)=P(k)-\frac{\int_0^k u^2P'(u)du}{k^2},\quad k\in \mathbb{R}.\label{energy: P}
\eeq
By {condition}~\eqref{E(0)} and interpolation inequality, 
\eq
E(\psi(0))<\infty.
\eeq
By Assumption~\ref{asp: global}, this implies that there is a energy conservation law. That is, $E(\psi(t))=E(\psi(0))$ for all $t\in \mathbb{R}$. Let 
\eq
p(z)=\sum\limits_{j=1}^N a_jz^j. \label{pz}
\eeq
By Eqs.~\eqref{energy: P} and~\eqref{pz}, H\"older's inequality and interpolation inequality, we obtain
\begin{align}
   & |a_N|^{1/(N+2)} \| \psi(t)\|_{L^{N+2}_x(\mathbb{R}^n)}\nonumber\\
   \lesssim & \left(\|\psi(t)\|_{H^1_x(\mathbb{R}^n)}^2+\sum\limits_{j=1}^{N-1} \| \psi(t)\|_{L^{j+2}_x(\mathbb{R}^n)}^{j+2}\right)^{1/(N+2)}\nonumber\\
    \lesssim & \left(\|\psi(t)\|_{H^1_x(\mathbb{R}^n)}^2+ \sum\limits_{j=1}^{N-1} \|\psi(t)\|_{L^2_x(\mathbb{R}^n)}^{b_j(N+2)}\|\psi(t)\|_{L^{N+2}_x(\mathbb{R}^n)}^{(1-b_j)(N+2)}\right)^{1/(N+2)}\nonumber\\
    \lesssim& \|\psi(t)\|_{H^1_x(\mathbb{R}^n)}^{2/(N+2)}+\|\psi(t)\|_{L^2_x(\mathbb{R}^n)}^{b_j}\|\psi(t)\|_{L^{N+2}_x(\mathbb{R}^n)}^{(1-b_j)}\vert_{j=N-1}+\|\psi(t)\|_{L^2_x(\mathbb{R}^n)}\label{aN}
\end{align}
where $b_j\in (0,1), j=1,\cdots,N-1,$ are numbers (determined by H\"older's inequality) satisfying 
\eq
\frac{1}{j+2}=\frac{b_j}{2}+\frac{1-b_j}{N+2}.
\eeq
With $a_N\neq 0$, since $\|\psi(t)\|_{L^{N+2}_x(\mathbb{R}^n)}$ in the right-hand side of Eq.~\eqref{aN} is sub-linear, Eq.~\eqref{aN} implies (The constant $E$ is defined in Assumption~\ref{asp: global})
\eq
\| \psi(t)\|_{L^{N+2}_x(\mathbb{R}^n)}\lesssim_{E}  1,
\eeq
which leads to estimate, by interpolation, 
\begin{align}
    \| \mathcal{N}(x,t,\psi(t))\psi(t)\|_{L^\infty_tL^1_x(\mathbb{R}^{n+1})}=&\| P(|\psi(t)|)\psi(t)\|_{L^\infty_tL^1_x(\mathbb{R}^{n+1})}\nonumber\\
    \lesssim & \| \psi(t)\|_{L^\infty_tL^{N+1}_x(\mathbb{R}^{n+1})}^{N+1}+\| \psi(t)\|_{L^\infty_tL^{2}_x(\mathbb{R}^{n+1})}^{2} \lesssim_E  1.
\end{align}
This completes the proof.

\end{proof}}

\subsection{Main results.} Let $\F_j(\lambda), j=c,1,2,$ denote smooth characteristic functions of the interval $[1,+\infty)$ satisfying
\eq\label{def: F}
\F_j(\lambda)=\begin{cases}
    1 & \text{ when }\lambda \geq 1\\
    0 & \text{ when }\lambda <\frac{1}{2}
\end{cases}, \quad j=c,1,2.
\eeq
{Here $F_1$ and $F_c$ denote the same cut-off function, but we use different symbols for convenience: $F_1$ is used for the frequency cut-off and $F_c$ for the spatial cut-off (or functions of $x - 2tP$).} For each $j\in \{c,1,2\}$, we define 
\eq
 \F_j(\lambda>a):=\F_j(\lambda/a)
\eeq
and
\eq\label{eq: Fj<}
\F_j(\lambda\leq a ):=1-\F_j(\lambda/a).
\eeq
In this paper, we restrict our discussion to the case where $t \geq 0$. The case where $t<0$ is treated similarly. Here are our main results. {Our first result states that when $\mathcal{N}(x,t,\psi(t))$ is a $L^p$ potential}, the adapted free channel wave operator, defined in~\eqref{wave-4}, exists on $L^2_x(\mathbb{R}^n), n\geq 3$. 
\begin{theorem}[$L^p$ potentials]\label{thm1}{ 
{Assume $\psi(t)$ is} the solution to system~\eqref{SE}. Let Assumptions~\ref{asp: global} and~\ref{asp: nonlocal} hold. }When the space dimension $n\geq3$, for all $\alpha\in(0,1-2/n)$, the adapted free channel wave operator acting on the initial data $\psi(0)$, defined as
\eq
\Omega_\alpha^*\psi(0):=s\text{-}\lim\limits_{t\to \infty} e^{itH_0} \F_c(\frac{|x-2tP|}{t^\alpha}\leq 1)\psi(t)\label{wave-4}
\eeq
exists in $L^2_x(\mathbb{R}^n)$. Furthermore, $\Omega_\alpha^*\psi(0)$ is independent on the choice of $\alpha$ in the following sense: for all $\alpha,\alpha'\in (0, 1-2/n)$, 
\eq\label{nonlocal weak}
\Omega_{\alpha}^*\psi(0)=\Omega_{\alpha'}^*\psi(0).
\eeq
\end{theorem}

\begin{remark}
The condition on the nonlinearity here is not sharp. One expects that the same proof can be applied to the abstract version, {where the free
dynamics is governed not by $H_0$ but by some other Fourier multiplier $\omega(P)$,} which involves a more general interaction, see Proposition~\ref{Abs}.
\end{remark}
{\begin{remark}{To control the non-free part, which has localization properties to be defined later, even in three or higher spatial dimensions, we require that the interaction term be localized or of the charge transfer type, as detailed in Theorems~\ref{thm} and~\ref{thm3}. It is important to note that the interaction can also be nonlinear.}
\end{remark}}
\begin{remark}\label{below-2}
In three or higher space dimensions, since we don't need the space localization to control the interaction term, the theorem  applies to the general nonlinear systems without using the spherical symmetry assumption.
\end{remark}
{\begin{remark}
   If Assumption~\ref{asp: global} is replaced by the weaker hypothesis that 
there exists merely a global weak solution in $L^{2}$, we expect that our proof of the 
existence of the adapted wave operator acting on initial data in $L^{2}$ continues to hold 
without essential modification. To be precise, the existence of a weak solution allows us to interpret the derivative of 
$e^{-i\Delta t}\psi(t)$ in the weak sense via the Duhamel formula.  
Specifically, we obtain
\[
\partial_{t}\!\,\bigl[e^{-i\Delta t}\psi(t)\bigr]
    =(-i) e^{-i\Delta t}\,\mathcal{N}(|x|,t,|\psi|)\psi(t),
\]
and by Assumption~\ref{asp: nonlocal}, this quantity is in fact bounded in $x$ for all $t>0$.  
Moreover, the relative propagation estimate introduced in Section~\ref{sec: RPE} remains valid.  
To verify this, note that for any bounded family of self-adjoint operators $B(t)$, one may interpret the identity
\begin{align*}
\frac{\partial}{\partial t}( \psi(t), B(t)\psi(t))_{L^2_x(\mathbb R^n)}
    =&
    ( (-i)\mathcal{N}\psi,\, \tilde{B}(t)e^{-i\Delta t}\psi(t))_{L^2_x(\mathbb R^n)}
    + ( e^{-i\Delta t}\psi(t),\, \tilde{B}(t)(-i)\mathcal{N}\psi)_{L^2_x(\mathbb R^n)}
    \\
    &+ \left( \psi(t), \frac{\partial}{\partial t}\tilde{B}(t)\,\psi(t)\right)_{L^2_x(\mathbb R^n)},
\end{align*}
where
\[
\tilde{B}(t) \equiv e^{-i\Delta t} B(t) e^{i\Delta t}.
\]

\end{remark}}

When $\mathcal{N}(x,t,\psi(t))$ is localized in $x$ variable, we can establish the existence of the adapted free channel wave operators in all space dimensions and provide some useful properties for the non-free part. We find that the non-free part of the solution is \emph{weakly localized} in the following sense.
\begin{definition}[The weakly localized part of the solution.]\label{def: weak} We say that part of the solution defined as $\psi_{wl}(x,t)$ is {\emph{weakly localized }} if it has non-zero mass and spreads slowly in the following sense: there exists $\beta\in (0,1)$ such that for all $t\geq 1$,
\begin{equation}
    (\psi_{wl}(x,t), |x|\psi_{wl}(x,t))_{L^2_x(\mathbb{R}^n)}\lesssim t^\beta
\end{equation}
holds true.
\end{definition}
\begin{theorem}[Space localized potentials]\label{thm} Consider $\psi(t)$ as the solution to system~\eqref{SE}. Let Assumptions~\ref{asp: global} and~\ref{asp: local} hold. Let $\delta$ be as in Assumption \ref{asp: local}. Then for all $n\geq 1$, there exist two positive constants $c_1, c_2$, which depend on $n$ and $\delta$ (see Eq.~\eqref{c1 and c2}), such that  for all $\alpha\in (0, c_1)$ and for all $\beta\in (0,\min\{c_2,\alpha\})$, the adapted free channel wave operator acting on $\psi(0)$, defined as
\eq
\Omega_{\alpha,\beta}^*\psi(0):=s\text{-}\lim\limits_{t\to \infty} e^{itH_0}\F_c(\frac{|x-2tP|}{t^\alpha}\leq 1) \F_1(t^{\beta}|P|>1)\psi(t),\label{wave-1}
\eeq
exists in $\s^2_x(\mathbb{R}^n)$ and $\Omega_{\alpha,\beta}^*\psi(0)$ is independent on the choices of $\alpha$ and $\beta$: for all $(\alpha,\beta), (\alpha',\beta')\in (0,c_1)\times (0,c_2)$ with $\beta<\alpha$ and $\beta'<\alpha'$, 
\eq\label{weak id}
\Omega_{\alpha,\beta}^*\psi(0)=\Omega_{\alpha',\beta'}^*\psi(0).
\eeq
Furthermore, if $\delta>2$, for every $\epsilon\in (0,1/2),$ we define the weakly localized component \( \psi_{wl}(t) \) by
\begin{equation}\label{def: weak: local}
    \psi_{wl}(t) := \F_c\left( \frac{|x|}{(t+1)^{1/2+\epsilon}} < 1 \right) \psi(t).
\end{equation} 
Then
\begin{enumerate}
\item the equation
\eq
\lim\limits_{t\to \infty}\|\psi(t)-e^{-itH_0}\Omega^*_{\alpha,\beta}\psi(0) -\psi_{wl}(t) \|_{\s^2_x(\mathbb{R}^n)}=0
\eeq
holds true;
\item \( \psi_{wl}(t) \) is weakly localized in the sense that
\begin{equation}\label{weakfo}
    \left( \psi_{wl}(t),\, |x|\, \psi_{wl}(t) \right)_{\mathcal{S}^2_x(\mathbb{R}^n)} \lesssim_\epsilon t^{1/2 + \epsilon}, \quad t \geq 1.
\end{equation}
\end{enumerate}
\end{theorem}
\begin{remark}
     {One expects that when 
\[
\mathcal{N}(x,t,\psi(t))
=\partial_j g_{jl}(x,t)\partial_l + I(x,t,\psi(t)),
\]
where $g$ denotes the metric tensor satisfying 
\begin{equation*}
\langle x\rangle^{\delta} g_{jl}(x,t) \in L^\infty_t W^{1,\infty}_x(\mathbb R^n\times \mathbb R)
\quad \text{for some } \delta>1,
\end{equation*}
and $I(x,t,\psi(t))$ satisfies Assumption~\ref{asp: local},
the adapted free channel wave operator acting on $\psi(0)$, defined by
\begin{equation}
\tilde \Omega_{\alpha,\beta}^*\psi(0)
:= s\text{-}\lim_{t\to \infty}
e^{itH_0}\,
\mathcal{F}_c\!\Big(\tfrac{|x-2tP|}{t^\alpha}\le 1\Big)
\mathcal{F}_1(t^\beta|P|>1)
\mathcal{F}_1(t^{-\beta}|P|\le 1)\,
\psi(t),
\end{equation}
exists in $L^2_x(\mathbb{R}^n)$ for all $n\ge 1$, provided that $\alpha,\beta>0$ are sufficiently small. Indeed, by rewriting
\[
\partial_j g_{jl}(x,t)\partial_l
= \partial_j\partial_l g_{jl}(x,t)
 - \partial_j\big[\partial_l g_{jl}(x,t)\big],
\]
and using the cutoff $\mathcal{F}_1(t^{-\beta}|P|\le 1)$ together with the observation that
\[
\partial_t\big[\mathcal{F}_1(t^{-\beta}|P|\le 1)\big]\ge 0,
\]
one can absorb the derivatives $\partial_j$ and $\partial_j\partial_l$ to obtain
\begin{equation}
\begin{aligned}
    &\Big\|
    e^{itH_0}\,
    \mathcal{F}_c\!\Big(\tfrac{|x-2tP|}{t^\alpha}\le 1\Big)
    \mathcal{F}_1(t^\beta|P|>1)
    \mathcal{F}_1(t^{-\beta}|P|\le 1)
    \big[\partial_j g_{jl}(x,t)\partial_l + I(x,t,\psi(t))\big]\psi(t)
    \Big\|_{L^2_x(\mathbb R^n)}
    \\
    &\hspace{6cm}\in L^1_t([1,\infty))
\end{aligned}
\end{equation}
for some sufficiently small $\alpha,\beta>0$.

}
\end{remark}
When $n\geq 5$ and $\mathcal{N}(x,t,\psi(t))$ is a charge-transfer interaction, we find that the non-free component consists of several moving, weakly localized parts, as defined below.
\begin{definition}[A moving weakly localized part]\label{def: weak: charge} Let $t\geq 0$ and $v\in \mathbb{R}^n$. We say that part of the solution defined as $\psi_{wl,v}(x,t)$ is {\emph{moving and weakly localized }} if it has mass and spreads slowly around $tv$, as $t\to \infty$, in the following sense: there exists $\beta\in (0,1)$ such that for all $t\geq 1$,
\begin{equation}
    (\psi_{wl,v}(x,t), |x-tv|\psi_{wl,v}(x,t))_{L^2_x(\mathbb{R}^n)}\lesssim t^\beta
\end{equation}
holds true.
    
\end{definition}

\begin{theorem}\label{thm3} Let $\psi(t)$ be the solution to system~\eqref{SE} satisfying Assumption~\ref{asp: global} with $a=1$. Let Assumption~\ref{asp: charge} be satisfied. Then when $n\geq 5$, for any $\epsilon\in (0,1/2)$, there exist $N$ moving weakly localized parts, $\psi_{wl,j}(t)\equiv\psi_{wl,\epsilon,j}(t), j=1,\cdots,N$, such that
\begin{enumerate}
\item the equation
\eq
\lim\limits_{t\to \infty}\|\psi(t)-e^{-itH_0}\Omega_{\alpha}^*\psi(0) -\sum\limits_{j=1}^N\psi_{wl,j}(t) \|_{\s^2_x(\mathbb{R}^n)}=0,\label{final}
\eeq
holds true; 
\item $\psi_{wl,j}(t), j=1,\cdots,N,$ are moving weakly localized parts around $tv_j$ satisfying
\eq
(\psi_{wl,j}(t), |x-tv_j|\psi_{wl,j}(t))_{\s^2_x(\mathbb{R}^n)}\lesssim_\epsilon t^{1/2+\epsilon},\qquad t\geq 1.
\eeq
\end{enumerate}
\end{theorem}

When the non-free part is small, we obtain scattering for system~\eqref{SE}:
\begin{proposition}\label{NLSapplication} Let $\psi(t)$ be the solution to system \eqref{SE} in $n=3$ space dimensions. Let us suppose that the following conditions are satisfied:
\begin{enumerate}
    \item Assumptions \ref{asp: global} and \ref{asp: nonlocal} are valid with $a=1$;

    \item $\mathcal{N}(x,t,\psi(t))$ is given by $\mathcal{N}(x,t,\psi(t))=I(|\psi(t)|)$ where $I: [0,\infty) \to \mathbb{R}$, is a function satisfying the following condition: for any pair $f(x), g(x)\in   \mathcal{H}_x^1(\mathbb{R}^3)$,
\begin{align}
&\|I( |f(x)|)f(x)-I(|g(x)|)g(x)\|_{ L^{6/5}_x(\mathbb{R}^{3})}\nonumber\\
\leq&  C_{I1} \| f(x)-g(x)\|_{\mathcal H^1_x}\| f(x)-g(x)\|_{L^6_x(\mathbb{R}^{3})}+C_{I2} \|g(x)\|_{L^6_x(\mathbb{R}^3)},\label{con: I}
\end{align}
where $C_{I1}=C_{I1}( \| f\|_{\mathcal{H}^1_x},\| g\|_{\mathcal{H}^1_x})>0$ and $C_{I2}=C_{I2}( \| f\|_{\mathcal{H}^1_x},\| g\|_{\mathcal{H}^1_x})>0$ are two positive constants dependent on $\| f\|_{\mathcal{H}^1_x}$ and $\| g\|_{\mathcal{H}^1_x}$.

\end{enumerate}
Then there exists {$m=m(E)> 0$ }such that 
\eq
\limsup\limits_{t\to \infty}\| \psi(t)-e^{-itH_0}\Omega_\alpha^*\psi(0)\|_{\mathcal{H}^1_x}<m
\eeq
implies 
\eq
\limsup\limits_{t\to \infty}\| \psi(t)-e^{-itH_0}\Omega_\alpha^*\psi(0)\|_{L^2_x(\mathbb{R}^3)}=0.
\eeq
\end{proposition}
\begin{remark}
The above proposition implies that if the non-free part is small, then it vanishes asymptotically. {This result is similar in nature to proofs of scattering below soliton mass threshold, but is not sharp.}
\end{remark}
{When $I(|\psi|)=\pm \lambda |\psi|^{p}$ with $\lambda>0$ for all $4/3\le p\le 4$, condition~\eqref{con: I} holds by the following argument. First observe that
\begin{align}
    |f|^{p}g - |g|^{p}f
    &= (|f|^{p}-|g|^{p})(f-g) + (|f|^{p}-|g|^{p})g + (f-g)|g|^{p}.
\end{align}
By H\"older’s inequality, using the fact that either $2(p-1)\in[2,6]$ or $6(p-1)\in[2,6]$, we have
\begin{align*}
    \bigl\| (|f|^{p}-|g|^{p})(f-g) \bigr\|_{L^{6/5}_x(\mathbb{R}^3)}
    &\le 
    \left\| \frac{|f|^{p}-|g|^{p}}{f-g} \right\|_{L^{2}_x(\mathbb{R}^3)+L^{6}_x(\mathbb{R}^3)}
    \|f-g\|_{\mathcal H^1_x}
    \|f-g\|_{L^{6}_x(\mathbb{R}^3)} \\
    &\lesssim 
    \bigl( \|f^{\,p-1}\|_{L^{2}_x+L^{6}_x}
          +\|g^{\,p-1}\|_{L^{2}_x+L^{6}_x} \bigr)
    \|f-g\|_{\mathcal H^1_x}
    \|f-g\|_{L^{6}_x} \\
    &\le C_{I1}\,\|f-g\|_{\mathcal H^1_x}\,\|f-g\|_{L^{6}_x}.
\end{align*}
Next, since $3p/2\in[2,6]$, we obtain
\begin{align*}
    \bigl\| (|f|^{p}-|g|^{p})g + (f-g)|g|^{p} \bigr\|_{L^{6/5}_x(\mathbb{R}^3)}
    &\lesssim 
    \bigl( \|f^{\,p}\|_{L^{3/2}_x(\mathbb{R}^3)}
          +\|g^{\,p}\|_{L^{3/2}_x(\mathbb{R}^3)} \bigr)
    \|g\|_{L^{6}_x(\mathbb{R}^3)} \\
    &\le C_{I2}\,\|g\|_{L^{6}_x(\mathbb{R}^3)},
\end{align*}
where $C_{I1}$ and $C_{I2}$ are positive constants depending on 
$\|f\|_{\mathcal H^1_x}$ and $\|g\|_{\mathcal H^1_x}$.
}

\section{Propagation Estimate, Relative Propagation Estimate, $tT$ potentials, estimates for interaction terms and commutator estimates}\label{sec: prop}
We let $b\in \mathbb{R}$ denote the lower bound on the time interval of interest.
\subsection{Propagation Estimate}
Given a family of self-adjoint operators $\{B(t)\}_{t\geq b}$, we define
\eq
\langle B\rangle_t:=(\psi(t), B(t)\psi(t))_{\s^2_x(\mathbb{R}^n)}=\int_{\mathbb{R}^n}  \psi(t)^*B(t)\psi(t)d^nx,\qquad t\geq b,
\eeq
where $\psi(t)$ is the solution to \eqref{SE} and $f^*$ denotes the conjugate of $f$ for any function $f$. Suppose $\langle B\rangle_t,$ for $ t\geq b,$ satisfies a boundedness condition  which is uniform over $t \in [b, \infty)$:
\eq
\sup\limits_{t\geq b} |\langle B\rangle_t|<\infty,
\eeq
and $\partial_t[\langle B\rangle_t]$ satisfies the decomposition, for all $t\geq b$,
\eqa\label{line1}
&\partial_t\langle B\rangle_t=\pm(\psi(t), C^*C\psi(t))_{\s^2_x(\mathbb{R}^n)}+g(t)\\
\label{line2}&g(t)\in L^1_t[b,\infty), \quad C^*C\geq 0.
\eeqa
We then refer to the family $\{ B(t)\}_{t\geq b}$ as a {\bf Propagation Observable} (PROB). See, for example, \cite{HWSS1999}, \cite{SS1987} and \cite{SS1988}.

From a PROB, by using {Eqs.~\eqref{line1} and ~\eqref{line2},} we derive {a} \textbf{Propagation Estimate} (PRES) for all $t_2\geq t_1\geq b$:
\begin{align}
\int_{t_1}^{t_2}\|C(t)\phi(t) \|_{\s^2_x(\mathbb{R}^n)}^2dt=&\pm ( \psi(t_2), B(t_2)\psi(t_2))_{L^2_x(\mathbb{R}^n)}\mp( \psi(t_1), B(t_1)\psi(t_1))_{L^2_x(\mathbb{R}^n)}-\int_{t_1}^{t_2}g(s) ds\nonumber\\
&\leq 2\sup\limits_{t\in [t_1,t_2]} \left|( \psi(t), B(t)\psi(t))_{\s^2_x(\mathbb{R}^n)}\right|+\|g(t)\|_{\s^1_t[b,\infty)}.\label{CC}
\end{align}

\subsection{Relative Propagation Estimate}
\label{sec: RPE}

We use a modified { PRES } as well. Given a family of self-adjoint operators $\{ \tilde B(t)\}_{t\geq b}$ and {a flow $\phi(t)$}, we define
\eq
\langle \tilde{B}: \phi(t)\rangle_t:=(\phi(t), \tilde{B}(t)\phi(t)  )_{\s^2_x(\mathbb{R}^n)}=\int_{\mathbb{R}^n}\phi(t)^*\tilde{B}(t)\phi(t)d^nx.
\eeq
Suppose $\langle \tilde B: \phi(t)\rangle_t,$ for $ t\geq b,$ satisfies a boundedness condition uniform in $t \in [b, \infty)$:
\eq
\sup\limits_{t\geq b} |\langle \tilde B: \phi(t)\rangle_t|<\infty,
\eeq
and assume that there exists a positive integer $N\in \mathbb{N}^+$ such that $\partial_t[\langle \tilde B: \phi(t)\rangle_t]$ satisfies the decomposition, for all $t\geq b$,
\eqa\label{line3}
&\partial_t[\langle \tilde B: \phi(t)\rangle_t]=\pm \sum\limits_{j=1}^N(\phi(t), C_j^*C_j\phi(t))_{\s^2_x(\mathbb{R}^n)}+g(t)\\
\label{line4}&g(t)\in L^1_t[b,\infty), \quad C_j^*C_j\geq 0.
\eeqa
We then refer to the family {$\{\tilde{B}(t)\}_{t\geq b}$} as a {\bf Relative Propagation Observable}(RPROB).

From a RPROB, by using {Eqs.~\eqref{line3} and ~\eqref{line4}}, we derive {a} {\bf Relative Propagation Estimate} (RPRES) for all $t_2\geq t_1\geq b$:
\eq\label{CC: ineq2}
{\sum\limits_{j=1}^N}\int_{t_1}^{t_2}\|C_j(t)\phi(t) \|_{\s^2_x(\mathbb{R}^n)}^2dt\leq \sup\limits_{t\geq b} \left|( \phi(t), \tilde{B}(t)\phi(t))_{\s^2_x(\mathbb{R}^n)}\right|+\|g(t)\|_{\s^1_t[b,\infty)}.
\eeq
 We will use $\phi(t)=e^{itH_0}\psi(t)$ {in the proof of Theorems~\ref{thm1} and~\ref{thm}}. 

We conclude this subsection by presenting an abstract version of the main proposition concerning the existence of the Free Channel Wave Operator. One expects that the proof of Proposition~\ref{Abs} can be derived using arguments similar to those used in Theorems~\ref{thm1} and ~\ref{thm}. A sketch of the proof is provided in Appendix ~\ref{app:3}.

\begin{proposition}\label{Abs}
Let $H_0=\omega(P)$ be the generator of the free evolution operator $U_0(t)\equiv e^{-iH_0 t}$ acting on a Hilbert space $\mathcal{H}=L_x^2(\R^n), n\geq 1.$
Let $\psi(t)$ be the solution of a Schr\"odinger type equation
$$
i\frac{\partial \psi(t)}{\partial t}=(H_0+\mathcal{N}(x,t,\psi(t)))\psi(t),\qquad t>0.
$$
Assume that for initial data $\psi(0)=\psi_0$ the solution of the above (possibly nonlinear) equation is global, uniformly bounded in $\mathcal{H}_x^1$ (that is, estimate~\eqref{con: H} is satisfied for $a=1$).
\begin{enumerate}
\item 
If the group $U_0(t)$ is bounded from $L_x^{p}(\mathbb{R}^n)$ into $L^{p'}_x(\mathbb{R}^n)$ with a bound that decays faster than $1/t^{1+\epsilon}$ for some $ \epsilon >0,$
where $1\leq p<2$ and $p'$ is the conjugate of $p$, then the following strong limit, defining the Free Channel Wave Operator acting on $\psi(0)$ exists in $L^2_x(\mathbb{R}^n)$: for all $\alpha\in (0,\frac{2p\epsilon}{(2-p)n})$,
{\begin{align}\label{Free}
\Omega^*_{free,\alpha}\psi(0)\equiv s\text{-}\lim_{ t\to \infty} U_0(-t)\F_c(\frac{|x-tv(P)|}{t^\alpha}\leq 1)\psi(t)
\end{align}}
provided the interaction term satisfies the following estimate:
\begin{equation}\label{interaction}
\sup\limits_{t\in \mathbb{R}^+}\| \mathcal{N}(x,t,\psi(t))\psi(t)\|_{L_x^p(\mathbb{R}^n)} \lesssim 1.
\end{equation}
Furthermore, $\Omega_\alpha^*\psi(0)$ is independent on the choice of $\alpha$ in the following sense: for all $\alpha,\alpha'\in (0, \frac{2p\epsilon}{(2-p)n})$, 
\eq
\Omega_{free,\alpha}^*\psi(0)=\Omega_{free,\alpha'}^*\psi(0).
\eeq
\item If there exist $ k> 0$, $p\in [1,2)$ and $\tilde p\in (2,\infty]$ such that the group $U_0(t)$ is bounded from $W_x^{k,p}(\mathbb{R}^n)$ into $L^{\tilde p}_x(\mathbb{R}^n)$ with a bound that decays faster than $1/t^{1+\epsilon}$ for some $ \epsilon >0$, then for all $\alpha\in (0, \min\{\frac{2 \epsilon \tilde p}{(\tilde p-2)n}, \frac{\epsilon \tilde p}{n}\}) $ and all $\beta\in (0, \min\{\alpha,\frac{\epsilon \tilde p-n\tilde p\alpha(1/2-1/\tilde p)}{\tilde p k} \}) $, the free channel wave operator acting on $\psi(0)$, defined by
\begin{equation}
    \Omega_{free,\alpha,\beta}^*\psi(0)=s\text{-}\lim\limits_{t\to \infty} e^{itH_0}\F_c(\frac{|x-tv(P)|}{t^\alpha}\leq 1) \F_1(|P|\leq t^\beta ) \psi(t),
\end{equation}
exists in $L^2_x(\mathbb{R}^n)$. Furthermore, $\Omega_{free,\alpha,\beta}^*\psi(0)$ is independent {\bl of} the choice of $\alpha$ and $\beta$ in the following sense: for all $\alpha,\alpha'\in (0, \min\{\frac{2 \epsilon \tilde p}{(\tilde p-2)n}, \frac{\epsilon \tilde p}{n}\}) $ and all $\beta\in (0, \frac{\epsilon \tilde p-n\alpha}{\tilde p k} ), \beta'\in (0, \frac{\epsilon \tilde p-n\alpha'}{\tilde p k} )$, 
\eq
\Omega_{free,\alpha,\beta}^*\psi(0)=\Omega_{free,\alpha',\beta'}^*\psi(0).
\eeq
\end{enumerate}
\end{proposition}

\subsection{Estimates for interaction terms}
We need the following dispersive estimates for the free flow: 
\begin{enumerate}
\item $L^p$ decay estimates, see for example Eq.~(1.1) of \cite{RSS}:
\eq
\|e^{-itH_0}f(x)\|_{\s^p_x(\mathbb{R}^n)}\lesssim_n \frac{1}{|t|^{\frac{n}{2}(\frac{1}{2}-\frac{1}{p})}}\|f(x)\|_{\s^{p'}_x(\mathbb{R}^n)},\quad f\in \s^{p'}_x(\mathbb{R}^n),\quad t\in \mathbb{R}-\{0\}\label{Lpdecay}
\eeq
where
\eq
\frac{1}{p}+\frac{1}{p'}=1, \quad 2\leq p\leq\infty.
\eeq
\item Local decay estimates, see the first proof in Appendix~\ref{app: free flows} for its proof: For $0<\alpha<1-\beta$ and $\beta\in (0,1/2)$,
\eq
\|\F_c(\frac{|x|}{t^\alpha}\leq 1)\F_1(t^\beta|P|>1)e^{\pm itH_0}\langle x\rangle^{- \sigma}\|_{L^2_x(\mathbb{R}^n)\to L^2_x(\mathbb{R}^n)}\lesssim \frac{1}{ t^{\sigma(1-\beta)}},\quad t\geq 1,\, \sigma\geq 0.\label{local: est}
\eeq
\end{enumerate}

In this section, we use the following notations
\eq
\F_1=\F_1(t^\beta |P|>1), \quad\text{  }\quad\F_c= \F_c(\frac{|x|}{t^\alpha}\leq 1)
\eeq
and 
\eq
\F_1^{(1)}=\F_1^{(1)}(k):=\frac{d}{dk}[\F_1(k)]
\eeq
when it does not lead to confusion. 
\subsubsection{Space localized $\mathcal{N}(x,t,\psi(t))$. }

\begin{proposition}\label{can1}Let $\delta$ be as in Assumption~\ref{asp: local}. Take 
{\eq\label{c1 and c2}
c_1=c_2=\frac{\delta-1}{4\delta}<\frac{1}{4}.
\eeq}
If Assumption~\ref{asp: local} is satisfied, then for all $\alpha\in (0, c_1) $ and $ \beta\in (0,c_2)$, the estimate
\begin{align}
&\|\F_c\F_1e^{itH_0}\mathcal{N}(x,t,\psi(t))\psi(t)\|_{ \s^2_x(\mathbb{R}^n)}\lesssim \frac{1}{t^{\frac{\delta+1}{2}}}\|\langle x\rangle^\delta\mathcal{N}(x,t,\psi(t))\psi(t)\|_{\s^\infty_t\s^2_x(\mathbb{R}^{n+1} )},\qquad t\geq 1\label{apsi}
\end{align}
is valid for $t\geq 1$.

\end{proposition}
\begin{proof}Since $\alpha <c_1<1-c_2<1-\beta$, by local decay estimates~\eqref{local: est}, we obtain 
\begin{align}
    \|\F_c\F_1e^{itH_0}\mathcal{N}(x,t,\psi(t))\psi(t)\|_{ \s^2_x(\mathbb{R}^n)}\lesssim& \frac{1}{t^{\delta(1-\beta)}}\| \F_c\langle x\rangle^\delta\|_{L^2_x(\mathbb{R}^n)\to L^2_x(\mathbb{R}^n)}\nonumber\\
    \lesssim & \frac{1}{t^{\delta(1-\alpha-\beta)}}\lesssim  \frac{1}{t^{\delta(1-c_1-c_2)}} \lesssim  \frac{1}{t^{\frac{\delta+1}{2}}}.
\end{align}

\end{proof}
\begin{remark} The space localization for $\mathcal{N}(x,t,\psi(t))$ is not needed in three or more space dimensions. This is because the {dispersive estimate} implies a decay rate faster than $\frac{1}{t^{1+0}}$. See Proposition \ref{can2}. Here, the discussion is mainly for the one and two space dimensional cases.
\end{remark}

\subsubsection{$L^p$ potentials.}\begin{proposition}\label{can2} If Assumptions~\ref{asp: global} and~\ref{asp: nonlocal} are valid, then for all $\alpha\in (0,1-2/n), n\geq 3$ and $t\geq 1$,
\begin{align}
\|\F_ce^{itH_0}\mathcal{N}(x,t,\psi(t))\psi(t)\|_{ \s^2_x(\mathbb{R}^n)}\lesssim_n \frac{1}{t^{1+\beta}}\|\mathcal{N}(x,t,\psi(t))\psi(t)\|_{\s^\infty_t\s^1_x(\mathbb{R}^{n+1} )}\label{can2eq}
\end{align}
where $\beta$ is given by 
\eq
\beta:=\frac{n(1-\alpha)}{2}-1>0.
\eeq
\end{proposition}
\begin{proof}
By using H\"older's inequality and $\s^\infty$ decay (estimate~\eqref{Lpdecay} with $p=\infty$), we obtain
\begin{align}
&\|\F_ce^{itH_0}\mathcal{N}(x,t,\psi(t))\psi(t)\|_{ \s^2_x(\mathbb{R}^n)}\nonumber\\
\leq & \|\F_c\|_{L^2_x(\mathbb{R}^n)}\| e^{itH_0}\|_{L^1_x(\mathbb{R}^n)\to L^\infty_x(\mathbb{R}^n)} \| \mathcal{N}(x,t,\psi(t))\psi(t)\|_{L^1_x(\mathbb{R}^n)}\nonumber\\
\lesssim_n & t^{\frac{\alpha n}{2}}\times \frac{1}{t^{n/2}}\| \mathcal{N}(x,t,\psi(t))\psi(t)\|_{L^\infty_tL^1_x(\mathbb{R}^{n+1})}\nonumber\\
\lesssim_n& \frac{1}{t^{\frac{n}{2}(1-\alpha)}}\| \mathcal{N}(x,t,\psi(t))\psi(t)\|_{L^\infty_tL^1_x(\mathbb{R}^{n+1})}.\label{eq: 3.26}
\end{align}

\end{proof}
\begin{remark} Based on the proof of Proposition \ref{can2}, $\s^\infty$ decay estimates of  the free flow are not necessary in $3$ or higher dimensions. For example, $\s^{6+\epsilon}$ decay {in $t$} will be sufficient in $3$ space dimensions provided $\alpha\ll1$. 
\end{remark}
Commutator estimates are required for identifying a positive term. Roughly speaking, consider an expression of the form $F(x)G(P)+G(P)F(x).$ In our applications both variables $x,P$ are scaled with a fractional power of $t.$
Suppose $F$ and $G$ are both positive, bounded and smooth. Then, the positive term, which is corresponding to $(\phi(t),C_jC_j^*\phi(t))_{L^2_x(\mathbb{R}^n)}$ in Eq.~\eqref{line3}, can be constructed as follows:
\eq
F(x)G(P)+G(P)F(x)=2\sqrt F G\sqrt F +[\sqrt F,[\sqrt F,G]]
\eeq
or 
\eq
F(x)G(P)+G(P)F(x)=2\sqrt G F \sqrt G +[\sqrt G,[\sqrt G,F]].
\eeq
The double commutator can be estimated using the commutator estimates provided below.
\begin{lemma}[
Commutator estimates.]\label{com} For all $\beta<\alpha,$ $l=0,1,$ and $t\geq 1$, the commutator estimate 
\eq
\|[\F_c,\F_1^{(l)}]\|_{\s^2_x(\mathbb{R}^n)\to \s^2_x(\mathbb{R}^n)}\lesssim_n \frac{1}{t^{\alpha-\beta}}
\eeq
holds true, where $\F_1^{(0)}\equiv \F_1$. 
\end{lemma}
\begin{proof} Let $\hat{\F}_1(\xi)$ denote the Fourier transform of $\F_1(x)$ in $x$ variable:
\eq
\hat{\F}_1(\xi)=\frac{1}{(2\pi)^{n/2}}\int_{\mathbb{R}^n} e^{-ix\cdot \xi} \F_1(x)d^nx.
\eeq
To compute $[  \F_c, \F_1^{(l)}]$, we find, with $\F_c=\F_c(\frac{|x|}{t^\alpha}\leq 1)$,
\begin{align}
[  \F_c, \F_1^{(l)}]=&\frac{1}{(2\pi)^{n/2}}\int\hat{\F}_1^{(l)}(\xi)e^{it^{\beta}P\cdot \xi}\times\left[e^{-it^{\beta}P\cdot \xi}\F_c e^{it^{\beta}P\cdot \xi}-\F_c \right] d^n\xi \nonumber\\
=&\frac{1}{(2\pi)^{n/2}}\int  \hat{\F}_1^{(l)}(\xi) e^{it^{\beta}P\cdot \xi}\left(\F_c(\frac{|x-t^{\beta}\xi|}{t^\alpha}\leq 1)-\F_c(\frac{|x|}{t^\alpha}\leq 1)\right) d^n\xi .
\end{align}
Using that by the mean-value Theorem,
\eq
\dfrac{\left|\F_c(\frac{|x-t^{\beta}\xi|}{t^\alpha}\leq 1)-\F_c(\frac{|x|}{t^\alpha}\leq 1) \right|}{t^{\beta-\alpha}|\xi|}\lesssim \sup\limits_{x\in \mathbb{R}^n}| \F_c'(|x|\leq 1)|\lesssim 1,
\eeq
we estimate
\begin{align}
\| [\F_c, \F_1^{(l)}]\|_{\s^2_x(\mathbb{R}^n)\to L^2_x(\mathbb{R}^n)}\lesssim_n& \frac{1}{t^{\alpha-\beta}}\int  | \hat{\F}_1^{(l)}(\xi)||\xi| d^n\xi\lesssim_n \frac{1}{t^{\alpha-\beta}}.
\end{align}

\end{proof}

\section{Proofs of Theorem~\ref{thm1} and Theorem~\ref{thm}}
In this section we prove Theorems~\ref{thm1} and ~\ref{thm}. The proof of Theorem~\ref{thm} requires the concept of forward/backward propagation waves. We adopt the following notations
\eq
\F_1^{(l)}=\F_1^{(l)}(t^\beta |P|>1) \quad\text{ or }\quad\F_1^{(l)}= \F_1^{(l)}(s^\beta |P|>1),\qquad l=0,1,
\eeq
and 
\eq
\F_c=\F_c(\frac{|x|}{t^\alpha}\leq 1)\quad\text{ or }\quad\F_c=\F_c(\frac{|x|}{s^\alpha}\leq 1)
\eeq
provided that it does not lead to any confusion.  
\subsection{Forward/backward propagation waves}
We start by discussing the concept of forward and backward propagation waves. These waves are analogous to the incoming and outgoing waves first introduced by Enss  \cite{Enss1978}. 

Let $S^{n-1}$ denote the unit sphere in $\mathbb{R}^n$. We define a class of functions on $S^{n-1}$, $\{F^{\hat{h}}(\xi)\}_{\hat{h}\in I}$, as a smooth partition of unity with an index set 
 \eq
 I=\{ \hat{h}_1,\cdots, \hat{h}_N\}\subseteq S^{n-1}\label{indexI}
 \eeq
 for some $N\in \N^+$, satisfying that there exists $c>0$ such that for every $\hat{h}_i\in I$, 
 \eq
 F^{\hat{h}}(\xi)=\begin{cases}1 & \text{ when }|\xi-\hat{h}|<c\\ 0 & \text{ when }|\xi-\hat{h}|>2c\end{cases}, \quad \xi\in S^{n-1}.\label{ceq1}
 \eeq
Given $\hat{h}\in I$, we define $\tilde F^{\hat h}: S^{n-1}\to \R,$ as another smooth cut-off function satisfying 
 \eq
 \tilde{F}^{\hat{h}}(\xi)=\begin{cases}1 & \text{ when }|\xi-\hat{h}|<4c\\ 0 & \text{ when }|\xi-\hat{h}|>8c\end{cases},\quad \xi\in S^{n-1}.\label{ceq2}
 \eeq
For $h\in \mathbb{R}^n-\{0\}$, we define $\hat h:=h/|h|$ and $\hat h=0$ when $h=0$. We also assume that $c>0$, defined in \eqref{ceq1} and \eqref{ceq2}, is properly chosen  such that for all $x,q\in \R^n$ with $x\neq 0$ and $q\neq 0$, 
{\eq
F^{\hat{h}}(\hat{x})\tilde{F}^{\hat{h}}(\hat{q})|x+q|\geq F^{\hat{h}}(\hat{x})\tilde{F}^{\hat{h}}\frac{1}{10}(|x|+|q|),\label{Feq1}
\eeq
and
\eq
F^{\hat{h}}(\hat{x})(1-\tilde{F}^{\hat{h}}(\hat{q}))|x-q|\geq F^{\hat{h}}(\hat{x})(1-\tilde{F}^{\hat{h}}(\hat{q}))\frac{1}{10^6}(|x|+|q|).\label{Feq2}
\eeq}
Now let us define {the space--velocity (or equivalently, space--frequency) smooth cutoff on the forward/backward propagation set} in terms of the phase-space $(r,v)\in \R^{n+n}$: 
\begin{definition}[{ Cutoff} on the forward/backward propagation set] {The smooth cutoffs }onto the forward and backward propagation sets, in terms of $(r,v)\in \mathbb{R}^{n+n}$, are defined as follows:
\eq
P^+(r,v):=\sum\limits_{b=1}^N F^{\hat{h}_b}(\hat{r})\tilde{F}^{\hat{h}_b}(\hat{v}),\label{Prv+}
\eeq
and 
\eq
P^-(r,v):=1-P^+(r,v),\label{Prv-}
\eeq
respectively. 
\end{definition}
\subsubsection{Estimates for Schr\"odinger operators.} With $P=-i\nabla_x$, we define the {cutoffs} $P^\pm$ as
\eq
P^\pm \equiv P^\pm (x, 2P).\label{Ppm}
\eeq
We need following estimates and their proofs can be found in Appendix~\ref{app: free flows}.
\begin{lemma}\label{Lem: Pprop} For all $\epsilon \in (0,1/2)$ and $s,t,\sigma\geq 0$, the estimates
\begin{align}
&\| \F_c(\frac{|x|}{(t+1)^{1/2+\epsilon}}\geq 1)P^\pm e^{\pm isH_0}\F_1(\sqrt{t+1}|P|\geq 1) \langle x\rangle^{-\sigma}\|_{L^2_x(\mathbb{R}^n)\to L^2_x(\mathbb{R}^n)}\nonumber\\
\lesssim_\epsilon &\frac{1}{\langle (t+1)^{1/2+\epsilon}+s/\sqrt{t+1}\rangle^\sigma}\label{free: est: 1}
\end{align}
and
\eq
\| P^\pm e^{\pm isH_0} \F_1((s+1)^{1/2-\epsilon}|P|\geq 1)\langle x\rangle^{-\sigma}\|_{L^2_x(\mathbb{R}^n)\to L^2_x(\mathbb{R}^n)}\lesssim_\epsilon \frac{1}{\langle s\rangle^{\sigma/2}}\label{free: est: 2}
\eeq
hold true.
\end{lemma}

\begin{lemma}\label{lem: Ppropfree} For all $\epsilon\in (0,1/2)$, $t,\sigma\geq 0$ and $s\in [0,t]$, the estimate 
\begin{align}
&\| \F_c(\frac{|x|}{(t+1)^{1/2+\epsilon}}\geq 1)P^+ e^{ -isH_0}\F_1(\sqrt{t+1}|P|< 1) \langle x\rangle^{-\sigma}\|_{L^2_x(\mathbb{R}^n)\to L^2_x(\mathbb{R}^n)}\nonumber\\
\lesssim_\epsilon &\frac{1}{ (t+1)^{\frac{1}{2}\sigma +\epsilon\sigma}}\label{free: est: 3}
\end{align}
holds true.
\end{lemma}

\begin{lemma}\label{lem: Ppmf} For all $f\in L^2_x(\mathbb{R}^n)$, we have
\begin{equation}
    \lim\limits_{s\to \infty}\| P^\pm e^{\pm i s H_0}f \|_{L^2_x(\mathbb{R}^n)}=0.
\end{equation}
    
\end{lemma}

\begin{lemma}\label{lem: not linear: local}
For all $f\in L^2_x(\mathbb{R}^n)$, $\alpha\in (0,1)$ and $s\geq 0$, we have
\begin{equation}
    \lim\limits_{s\to \infty}\|\chi(|x|\leq s^\alpha) P^\mp e^{\pm i s H_0}f \|_{L^2_x(\mathbb{R}^n)}=0.\label{goal: lem: eq1}
\end{equation}

\end{lemma}

\subsubsection{ Estimates for translated Schr\"odinger operators. } For $\eta\in \mathbb{R}^n$ and $t\geq 0$, we define the {cutoffs} $P_{t\eta}^\pm\equiv P_{\eta,t\eta}^\pm$ as 
\eq
P^\pm_{t\eta}\equiv P^\pm (x-t\eta, 2P-\eta).
\eeq
We need following estimates for charge-transfer problems and their proofs can be found in Appendix~\ref{app: free flows}.
\begin{lemma}\label{Lem: Pprop: charge} For all $\eta\in \mathbb{R}^n$, and $ t,\sigma\geq 0$, the estimates, 
\begin{align}
    &\| \F_c(\frac{|x-t\eta|}{(t+1)^{1/2+\epsilon}}\geq 1)P^\pm_{t\eta} e^{ i(s-t)H_0}\F_1(\sqrt{t+1}|2P-\eta|\geq 1) \langle x-s\eta\rangle^{-\sigma}\|_{L^2_x(\mathbb{R}^n)\to L^2_x(\mathbb{R}^n)}\nonumber\\
\lesssim_\epsilon &\frac{1}{\langle (t+1)^{1/2+\epsilon}+s\sqrt{t+1}\rangle^\sigma},\qquad \text{sgn}(s-t)=\pm,\label{free: charge: goal1}
\end{align}
and
\begin{align}
    \| P^-_{t\eta} e^{- itH_0}\F_1(\sqrt{t+1}|2P-\eta|\geq 1) \langle x\rangle^{-\sigma}\|_{L^2_x(\mathbb{R}^n)\to L^2_x(\mathbb{R}^n)}\lesssim_\epsilon &\frac{1}{\langle t\rangle^{\sigma/2}}\label{free: charge: goal2}
\end{align}
hold true. 
    
\end{lemma}

\begin{lemma}\label{lem: Ppropfree: charge} For all $\epsilon\in (0,1/2)$, $\eta\in \mathbb{R}^n$, $t,\sigma\geq 0$ and $s\in [0,t]$, the estimate 
\begin{align}
&\| \F_c(\frac{|x-t\eta|}{(t+1)^{1/2+\epsilon}}\geq 1)P_{t\eta}^+ e^{ -i(t-s)H_0}\F_1(\sqrt{t+1}|2P-\eta|< 1) \langle x-s\eta\rangle^{-\sigma}\|_{L^2_x(\mathbb{R}^n)\to L^2_x(\mathbb{R}^n)}\nonumber\\
\lesssim_\epsilon &\frac{1}{ (t+1)^{\frac{1}{2}\sigma +\epsilon\sigma}}\label{free: charge: eq6}
\end{align}
holds true.
\end{lemma}

\begin{lemma}\label{lem: Ppmf: charge} For all $f\in L^2_x(\mathbb{R}^n)$ and $\eta\in \mathbb{R}^n$, we have
\begin{equation}
    \lim\limits_{t\to \infty}\| P^-_{t\eta} e^{- i t H_0}f \|_{L^2_x(\mathbb{R}^n)}=0.
\end{equation}
    
\end{lemma}

\begin{lemma}\label{lem: not linear: charge}
For all $f\in L^2_x(\mathbb{R}^n)$, $\eta\in \mathbb{R}^n$, $\alpha\in (0,1)$ and $t\geq 0$, we have
\begin{equation}
  \lim\limits_{t\to \infty}  \|\chi(|x-t\eta|\leq t^\alpha) P_{t\eta}^\pm e^{- i t H_0}f \|_{L^2_x(\mathbb{R}^n)}=0.\label{charge: limit: goal1}
\end{equation}

\end{lemma}

\subsection{ Proof of Theorem~\ref{thm1}.}\label{sec: proof of Thm1}

\begin{proof}[Proof of Theorem \ref{thm1}] We define
\begin{align}
\Omega_\alpha^*(t)\psi(0) :=e^{itH_0}\F_c(\frac{|x-2tP|}{t^\alpha}\leq 1) \psi(t).\label{def: omegaalphat}
\end{align}
By Eqs.~\eqref{id: Fc} and ~\eqref{def: omegaalphat}, $\Omega_\alpha^*(t)\psi(0)$ reads, with $\F_c\equiv \F_c(\frac{|x|}{t^\alpha}\leq 1),$
\eq
\Omega_\alpha^*(t)\psi(0) =\F_ce^{itH_0} \psi(t).\label{def: omegaalphat2}
\eeq
{In what follows, we use 
\eq
\F_c=\F_c(\frac{|x|}{t^\alpha}\leq 1)\quad\text{ or }\quad\F_c= \F_c(\frac{|x|}{s^\alpha}\leq 1),
\eeq
when it does not lead to confusion. }
Using Cook's method to expand $\Omega_{\alpha}^*(t)\psi(0)$, we obtain
\begin{align}
\Omega_{\alpha}^*(t)\psi(0)=&\Omega_{\alpha}^*(1)\psi(0)+(-i)\int_1^t  \F_ce^{isH_0}\mathcal{N}(x,s,\psi(s))\psi(s)ds+\int_1^t  \partial_s[\F_c]e^{isH_0}\psi(s)ds\nonumber\\
=:&\Omega_{\alpha}^*(1)\psi(0)+\psi_{int}(t)+\psi_{p}(t).
\end{align}
By the unitarity of $e^{iH_0}$, Assumption~\ref{asp: global} and Eq.~\eqref{def: omegaalphat2}, 
{\eq
\Omega_\alpha^*(1)\psi(0)\in L^2_x(\mathbb{R}^n). \label{exist: 0: nonlocal}
\eeq}
By $\alpha\in (0, 1-2/n), n\geq 3$ and Proposition~\ref{can2}, $\psi_{int}(t)$ satisfies the estimate, with $\beta=\frac{n(1-\alpha)}{2}-1$,
\begin{align}
    \| \psi_{int}(t)\|_{L^2_x(\mathbb{R}^n)}\leq & \int_1^t \| \F_ce^{isH_0}\mathcal{N}(x,s,\psi(s))\psi(s)\|_{L^2_x(\mathbb{R}^n)}ds\nonumber\\
    \lesssim_n & \int_1^t \frac{1}{t^{1+\beta}} \| \mathcal{N}(x,s,\psi(s))\psi(s)\|_{L^\infty_tL^1_x(\mathbb{R}^{n+1})}ds\nonumber\\
    \lesssim_n & \| \mathcal{N}(x,s,\psi(s))\psi(s)\|_{L^\infty_tL^1_x(\mathbb{R}^{n+1})},
\end{align}
which implies that 
{\eq
\lim\limits_{t\to \infty}\psi_{int}(t)\text{ exists in } L^2_x(\mathbb{R}^n). \label{exist: int: nonlocal}
\eeq}
For $\psi_{p}(t)$, we use {\bf{RPRES}} by taking {$b=1$ and}
\eq
\begin{cases}
B(t):=\F_c\\
\phi(t)=e^{itH_0}\psi(t)
\end{cases},\qquad t\geq 1.
\eeq
We find that 
\begin{align}
    \partial_t\langle B: \phi(t)\rangle_t=& (\phi(t), \partial_t[\F_c] \phi(t))_{L^2_x(\mathbb{R}^n)}+(-i)(\phi(t), \F_ce^{itH_0}\mathcal{N}(x,t,\psi(t))\psi(t) )_{L^2_x(\mathbb{R}^n)}\nonumber\\
    &+i(\F_ce^{itH_0}\mathcal{N}(x,t,\psi(t))\psi(t) , \phi(t))_{L^2_x(\mathbb{R}^n)}\nonumber\\
    =& (\phi(t), C^*C\phi(t))_{L^2_x(\mathbb{R}^n)}+g(t)
\end{align}
where $C^*C$ and $g(t)$ are given by, with $\F'_c(\lambda\leq 1)\equiv \frac{d}{d\lambda}[F_c(\lambda\leq 1)]$, 
\begin{align}
C^*C:=&\partial_t[\F_c]= \F_c'(\frac{|x|}{t^\alpha}\leq 1)\times \frac{-\alpha }{t}\times \frac{|x|}{t^\alpha}\geq  0\label{Fcnneg}
\end{align}
and
\begin{align}
    g(t):=& (-i)(\phi(t), \F_ce^{itH_0}\mathcal{N}(x,t,\psi(t))\psi(t) )_{L^2_x(\mathbb{R}^n)}+i(\F_ce^{itH_0}\mathcal{N}(x,t,\psi(t))\psi(t) , \phi(t))_{L^2_x(\mathbb{R}^n)}. 
\end{align}
{We observe that $\langle B: \phi(t)\rangle_t$ is uniformly bounded over $t$. Utilizing the Cauchy-Schwarz inequality, the unitarity of $e^{itH_0}$, and Assumption~\ref{asp: global}, we have:
\begin{align}
| \langle B: \phi(t)\rangle_t|=  ( e^{itH_0}\psi(t), \F_c e^{itH_0}\psi(t))_{L^2_x(\mathbb{R}^n)} \leq  \| e^{itH_0}\psi(t)\|^2_{L^2_x(\mathbb{R}^n)}=  \| \psi(t)\|^2_{L^2_x(\mathbb{R}^n)}\lesssim  E^2.\label{supBphi}
\end{align}
Furthermore, $g(t) \in L^1_t[1,\infty)$. To be precise, by applying the Cauchy-Schwarz inequality, and using Assumptions~\ref{asp: global} and~\ref{asp: nonlocal}, along with Proposition~\ref{can2}, $g(t)$ satisfies the following estimate:}
\begin{align}
    |g(t)|\leq & 2\| \phi(t)\|_{L^2_x(\mathbb{R}^n)}\| \F_ce^{itH_0}\mathcal{N}(x,t,\psi(t))\psi(t)\|_{L^2_x(\mathbb{R}^n)}\nonumber\\
    \lesssim_n& \frac{E}{t^{1+\beta}}\| \mathcal{N}(x,t,\psi(t))\psi(t)\|_{L^\infty_tL^1_x(\mathbb{R}^{n+1})}\in  L^1_t[1,\infty).\label{ineq: g}
\end{align}
Hence, the family $\{B(t)\}_{t\in [1,\infty)}$ is a RPROB with respect to $\phi(t)=e^{itH_0}\psi(t)$ and by Eqs.~\eqref{CC: ineq2}, ~\eqref{supBphi} and ~\eqref{ineq: g}, we obtain
\begin{align}
   \int_1^\infty  |(\phi(t), \partial_t[\F_c]\phi(t))_{L^2_x(\mathbb{R}^n)}|dt=& \int_1^\infty  (\phi(t), \partial_t[\F_c]\phi(t))_{L^2_x(\mathbb{R}^n)}dt\leq  2 \sup\limits_{t\in [1,\infty)}| \langle B: \phi(t)\rangle_t|+\| g(t)\|_{L^1_t[1,\infty)}\nonumber\\
   \lesssim_n & E^2+E\| \mathcal{N}(x,t,\psi(t))\psi(t)\|_{L^\infty_tL^1_x(\mathbb{R}^{n+1})} < \infty.\label{p: nonlocal}
\end{align}
By Cauchy-Schwarz inequality and the non-negativity of $\partial_t[\F_c]$ (see~\eqref{Fcnneg}), $\psi_p(t)$ satisfies the estimate, for $T_2\geq T_1\geq 1$,
\begin{align}
  \| \psi_p(T_2)-\psi_p(T_1)\|_{L^2_x(\mathbb{R}^n)}\leq & \| \int_{T_1}^{T_2} | \partial_t[\F_c] \phi(t)|  dt \|_{L^2_x(\mathbb{R}^n)}\nonumber\\
  \leq & \| \left( \int_{T_1}^{T_2} \partial_t[\F_c]dt\right)^{1/2}  \left( \int_{T_1}^{T_2} \partial_t[\F_c]|\phi(t)|^2dt\right)^{1/2}\|_{L^2_x(\mathbb{R}^n)}.\label{psipT2T1}
\end{align}
By estimates 
\begin{align}
    \int_{T_1}^{T_2} \partial_t[\F_c]dt=& \F_c(\frac{|x|}{t^\alpha}\leq 1)\vert_{t=T_1}^{t=T_2}\leq  \F_c(\frac{|x|}{{T_2}^\alpha}\leq 1)\leq  1\label{dtFc}
\end{align}
and~\eqref{p: nonlocal}, estimate~\eqref{psipT2T1} leads to 
\begin{align}
    \| \psi_p(T_2)-\psi_p(T_1)\|_{L^2_x(\mathbb{R}^n)}\leq & \|   \left( \int_{T_1}^{T_2} \partial_t[\F_c]|\phi(t)|^2dt\right)^{1/2}\|_{L^2_x(\mathbb{R}^n)}=\left( \int_{\mathbb{R}^n}\int_{T_1}^{T_2}   \partial_t[\F_c]|\phi(t)|^2 dtd^nx\right)^{1/2} \nonumber\\
    =& \left( \int_{T_1}^{T_2}  |(\phi(t), \partial_t[\F_c]\phi(t))_{L^2_x(\mathbb{R}^n)}|dt\right)^{1/2}\to  0
\end{align}
as $T_1\to \infty$, where we have used 
\eq
\int_{\mathbb{R}^n}\int_{T_1}^{T_2}   \partial_t[\F_c]|\phi(t)|^2 dtd^nx\lesssim (T_2-T_1) \sup\limits_{t\in \mathbb{R}}\| \psi(t)\|^2_{L^2_x(\mathbb{R}^n)}<\infty
\eeq
and Fubini's Theorem to switch the order of integration. Hence, $\{ \psi_p(t)\}_{t\geq 1}$ is Cauchy in $L^2_x(\mathbb{R}^n)$ and therefore 
{\eq
\lim\limits_{t\to\infty}\psi_{p}(t)\text{ exists in } L^2_x(\mathbb{R}^n). \label{exist: p: nonlocal}
\eeq}
Eq.~\eqref{exist: p: nonlocal}, together with Eqs.~\eqref{exist: 0: nonlocal} and ~\eqref{exist: int: nonlocal}, implies that
{\eq
\Omega_\alpha^*\psi(0)\equiv \lim\limits_{t\to \infty}\Omega_\alpha^*(t)\psi(0)\text{ exists in } L^2_x(\mathbb{R}^n). \label{exist: c: nonlocal}
\eeq}
We also have that for all $\alpha,\alpha'\in (0,1-2/n)$ and {$\varphi\in L^2_x(\mathbb{R}^n)$}, by Cauchy-Schwarz inequality, Assumption~\ref{asp: global} and the unitarity of $e^{itH_0}$,
{\begin{align}
    \left|({\varphi}, \Omega_\alpha^*(t)\psi(0)-\Omega_{\alpha'}^*(t)\psi(0))_{L^2_x(\mathbb{R}^n)}\right|=& \left|( (\F_c(\frac{|x|}{t^\alpha}\leq 1)-\F_c(\frac{|x|}{t^{\alpha'}}\leq 1)){\varphi}, e^{itH_0}\psi(t))_{L^2_x(\mathbb{R}^n)}\right|\nonumber\\
    \leq & \|(\F_c(\frac{|x|}{t^\alpha}\leq 1)-\F_c(\frac{|x|}{t^{\alpha'}}\leq 1)){\varphi} \|_{L^2_x(\mathbb{R}^n)}\| e^{itH_0}\psi(t)\|_{L^2_x(\mathbb{R}^n)}\nonumber\\
    \leq & E\|(\F_c(\frac{|x|}{t^\alpha}\leq 1)-\F_c(\frac{|x|}{t^{\alpha'}}\leq 1)){\varphi} \|_{L^2_x(\mathbb{R}^n)}\nonumber\\
    \to & 0
\end{align}}
as $t\to \infty$. This implies 
\eq
w\text{-}\lim\limits_{t\to \infty} \Omega_\alpha^*(t)\psi(0)-\Omega_{\alpha'}^*(t)\psi(0)=0,\quad \text{ in }L^2_x(\mathbb{R}^n)
\eeq
and therefore, due to the existence of $\Omega_\alpha^*\psi(0)$ and $\Omega_{\alpha'}^*\psi(0)$ in $L^2_x(\mathbb{R}^n)$ in strong sense, Eq.~\eqref{nonlocal weak}. 

\end{proof}
\subsection{ Proof of Theorem~\ref{thm}.}
The proof of Theorem~\ref{thm} requires following proposition and lemma.
\begin{proposition}\label{prof: free: local}Let $\alpha,\beta$ and $n$ be as in Theorem~\ref{thm}. Let Assumptions~\ref{asp: global} and \ref{asp: local} be satisfied. Then $\Omega_{\alpha,\beta}^*\psi(0),$ defined in~\eqref{wave-1}, exists in $L^2_x(\mathbb{R}^n)$ and Eq.~\eqref{weak id} is valid for all $(\alpha,\beta), (\alpha',\beta')\in (0,c_1)\times (0,c_2)$ with $\beta<\alpha$ and $\beta'<\alpha'$.    
\end{proposition}

\begin{proof} Let
\begin{align}
\Omega_{\alpha,\beta}^*(t)\psi(0):=&e^{itH_0}\F_c(\frac{|x-2tP|}{t^\alpha}\leq 1)\F_1(t^\beta|P|>1)\psi(t).\label{omega: alpha: beta: t}
\end{align}
By Eq.~\eqref{id: Fc}, $\Omega_{\alpha,\beta}^*(t)\psi(0)$ reads, with $\F_c\equiv \F_c(\frac{|x|}{t^\alpha}\leq 1)$ and $\F_1\equiv \F_1(t^\beta|P|>1)$,
\eq
\Omega_{\alpha,\beta}^*(t)\psi(0)=\F_c\F_1e^{itH_0}\psi(t).
\eeq
{In what follows, we use 
\eq
\F_c=\F_c(\frac{|x|}{t^\alpha}\leq 1)\quad\text{ or }\quad\F_c= \F_c(\frac{|x|}{s^\alpha}\leq 1)
\eeq
and
\eq
\F_1=\F_1(t^\beta|P|> 1)\quad\text{ or }\quad\F_1=\F_1(s^\beta|P|> 1),
\eeq
when it does not lead to confusion. }
By Cook's method to expand $\Omega_{\alpha,\beta}^*(t)\psi(0)$, $\Omega_{\alpha,\beta}^*(t)\psi(0)$ can be rewritten as
\begin{align}
\Omega_{\alpha,\beta}^*(t)\psi(0)=&\Omega_{\alpha,\beta}^*(1)\psi(0)+(-i)\int_1^t  \F_c\F_1e^{isH_0}\mathcal{N}(x,s,\psi(s))\psi(s)ds\nonumber\\
&+\int_1^t  \partial_s[\F_c]\F_1e^{isH_0}\psi(s)ds+\int_1^t  \partial_s[\F_1] \F_ce^{isH_0}\psi(s)ds\nonumber\\
&+\int_1^t  [\F_c,\partial_s[\F_1]]e^{isH_0}\psi(s)ds \nonumber\\
=:&\Omega_{\alpha,\beta}^*(1)\psi(0)+\psi_{int}(t)+\psi_{p,1}(t)+\psi_{p,2}(t)+\psi_c(t).\label{Cook: local}
\end{align}
By the unitarity of $e^{iH_0}$, Assumption~\ref{asp: global} and Eq.~\eqref{omega: alpha: beta: t}, 
{\eq
\Omega_{\alpha,\beta}^*(1)\psi(0)\in L^2_x(\mathbb{R}^n). \label{exist: 0: local}
\eeq}
By Proposition \ref{can1}, we obtain that
\eq
\|\F_c\F_1e^{itH_0}\mathcal{N}(x,t,\psi(t))\psi(t)\|_{L^2_x(\mathbb{R}^n)}\in L^1_t[1,\infty),
\eeq
which implies that 
{\eq
\lim\limits_{t\to \infty}\psi_{int}(t)\text{ exists in } L^2_x(\mathbb{R}^n). \label{exist: int: local}
\eeq}
For $\psi_{p,1}(t)$ and $\psi_{p,2}(t)$, we use {\bf{RPRES}} by taking, {with $b=1$(see~\eqref{line2} for the definition of $b$),}
\eq
\begin{cases}
B_1(t):=\F_1\F_c\F_1\\
\phi(t)=e^{itH_0}\psi(t)
\end{cases}
\eeq
and
\eq
\begin{cases}
B_2(t):=\sqrt{\F_1}\F_c^2\sqrt{\F_1}\\
\phi(t)=e^{itH_0}\psi(t)
\end{cases},
\eeq
respectively. 

We begin with {\bf RPRES} {for} $\psi_{p,1}(t).$ We find that 
\begin{align}
    &\partial_t\langle B_1: \phi(t)\rangle_t\nonumber\\
    =& (\phi(t), \partial_t[\F_1\F_c\F_1] \phi(t))_{L^2_x(\mathbb{R}^n)}+(-i)(\phi(t), \F_1\F_c\F_1e^{itH_0}\mathcal{N}(x,t,\psi(t))\psi(t) )_{L^2_x(\mathbb{R}^n)}\nonumber\\
    &+i(\F_1\F_c\F_1e^{itH_0}\mathcal{N}(x,t,\psi(t))\psi(t) , \phi(t))_{L^2_x(\mathbb{R}^n)}\label{Bphi: 1}
\end{align}
By 
\begin{align}
    \partial_t[\F_1\F_c\F_1]=& \F_1\partial_t[\F_c]\F_1+\partial_t[\F_1]\F_c\F_1+\F_1\F_c\partial_t[\F_1]\nonumber\\
    =& \F_1\partial_t[\F_c]\F_1+ 2\sqrt{\F_c}\F_1\partial_t[\F_1]\sqrt{\F_c}+[\partial_t[\F_1], \sqrt{\F_c}]\sqrt{\F_c}\F_1\nonumber\\
    &+\sqrt{\F_c} \partial_t[\F_1][\sqrt{\F_c},\F_1]+\F_1\sqrt{\F_c}[\sqrt{\F_c},\partial_t[\F_1]]\nonumber\\
    &+[\F_1,\sqrt{\F_c}]\partial_t[\F_1]\sqrt{\F_c}\nonumber\\
    =: & \F_1\partial_t[\F_c]\F_1+ 2\sqrt{\F_c}\F_1\partial_t[\F_1]\sqrt{\F_c}+F_1(t),
\end{align}
where $F_1(t)$ is given by 
\begin{align}
F_1(t):=&[\partial_t[\F_1], \sqrt{\F_c}]\sqrt{\F_c}\F_1\nonumber\\
    &+\sqrt{\F_c} \partial_t[\F_1][\sqrt{\F_c},\F_1]+\F_1\sqrt{\F_c}[\sqrt{\F_c},\partial_t[\F_1]]\nonumber\\
    &+[\F_1,\sqrt{\F_c}]\partial_t[\F_1]\sqrt{\F_c},\label{F1: express}
\end{align}
Eq.~\eqref{Bphi: 1} implies 
\begin{align}
    \partial_t\langle B_1: \phi(t)\rangle_t= & \sum\limits_{j=1}^2 (\phi(t), C_j^*C_j\phi(t))_{L^2_x(\mathbb{R}^n)}+g(t)
\end{align}
where $C_j^*C_j, j=1,2,$ and $g(t)$ are given by, with $\F'_c(\lambda\leq 1)\equiv \frac{d}{d\lambda}[F_c(\lambda\leq 1)]$, 
\begin{align}
C_1^*C_1:=&\F_1\partial_t[\F_c]\F_1\nonumber\\
=&\F_1 \F_c'(\frac{|x|}{t^\alpha}\leq 1)\times \frac{-\alpha }{t}\times \frac{|x|}{t^\alpha}\F_1\nonumber\\
\geq & 0,\label{Fcnneg: 1: 1}
\end{align}
\begin{align}
C_2^*C_2:=&2\sqrt{\F_c}\partial_t[\F_1]\F_1 \sqrt{\F_c}\nonumber\\
=&\sqrt{\F_c}\F_1 \F_1'(t^\beta|P|> 1)\times \frac{\beta }{t}\times t^\beta|P| \sqrt{\F_c}\nonumber\\
\geq & 0,\label{Fcnneg: 1: 2}
\end{align}
and
\begin{align}
    g(t):=& (-i)(\phi(t), \F_1\F_c\F_1e^{itH_0}\mathcal{N}(x,t,\psi(t))\psi(t) )_{L^2_x(\mathbb{R}^n)}\nonumber\\
    &+i(\F_1\F_c\F_1e^{itH_0}\mathcal{N}(x,t,\psi(t))\psi(t) , \phi(t))_{L^2_x(\mathbb{R}^n)}+(\phi(t),F_1(t)\phi(t))_{L^2_x(\mathbb{R}^n)},
\end{align}
{respectively. }Here, we note that $\langle  B_1: \phi(t)\rangle_t$ is uniformly bounded in $t$: by Cauchy-Schwarz inequality, unitarity of $e^{itH_0}$ and Assumption~\ref{asp: global},
\begin{align}
| \langle B_1: \phi(t)\rangle_t|= & ( e^{itH_0}\psi(t), \F_1\F_c\F_1 e^{itH_0}\psi(t))_{L^2_x(\mathbb{R}^n)} \nonumber\\
\leq & \| e^{itH_0}\psi(t)\|^2_{L^2_x(\mathbb{R}^n)}\nonumber\\
= & \| \psi(t)\|^2_{L^2_x(\mathbb{R}^n)}\nonumber\\
\lesssim & E^2,\label{supBphi: 1}
\end{align}
and $g(t)\in L^1_t[1,\infty)$: by Lemma~\ref{com} and Eq.~\eqref{F1: express}, with $\alpha>\beta$,
\begin{align}
    &| (\phi(t),F_1(t)\phi(t))_{L^2_x(\mathbb{R}^n)}|\nonumber\\
    \lesssim_E & \|[\partial_t[\F_1], \sqrt{\F_c}] \|_{L^2_x(\mathbb{R}^n)\to L^2_x(\mathbb{R}^n)} +\| \partial_t[\F_1]\|_{L^2_x(\mathbb{R}^n)\to L^2_x(\mathbb{R}^n)}\| [\sqrt{\F_c}, \F_1]\|_{L^2_x(\mathbb{R}^n)\to L^2_x(\mathbb{R}^n)}\nonumber\\
    &+\| [\sqrt{\F_c}, \partial_t[\F_1]]\|_{L^2_x(\mathbb{R}^n)\to L^2_x(\mathbb{R}^n)}+\| [ \F_1,\sqrt{\F_c}]\|_{L^2_x(\mathbb{R}^n)\to L^2_x(\mathbb{R}^n)}\| \partial_t[\F_1]\|_{L^2_x(\mathbb{R}^n)\to L^2_x(\mathbb{R}^n)}\nonumber\\
    \lesssim_{E,n}& \frac{1}{t^{1+\alpha-\beta}}\nonumber\\
    \lesssim_{E,n}&\in L^1_t[1,\infty),\label{com: ineq: 1}
\end{align}
where the extra factor $\frac{1}{t}$ comes from $\partial_t[\F_1]$. Estimate~\eqref{com: ineq: 1}, together with Cauchy-Schwarz inequality, Assumptions~\ref{asp: global} and~\ref{asp: local} and Proposition~\ref{can1}, implies the estimate, with $\delta>1$ and $\alpha>\beta,$ 
\begin{align}
    |g(t)|\leq & 2\| \phi(t)\|_{L^2_x(\mathbb{R}^n)}\| \F_c\F_1e^{itH_0}\mathcal{N}(x,t,\psi(t))\psi(t)\|_{L^2_x(\mathbb{R}^n)}+|(\phi(t),F_1(t)\phi(t))_{L^2_x(\mathbb{R}^n)}|\nonumber\\
    \lesssim_{n,E}& \frac{1}{t^{\frac{\delta+1}{2}}}\|\langle x\rangle^\delta\mathcal{N}(x,t,\psi(t))\psi(t)\|_{\s^\infty_t\s^1_x(\mathbb{R}^{n+1} )}+\frac{1}{t^{1+\alpha-\beta}}\nonumber\\
    \in & L^1_t[1,\infty).\label{ineq: g: 1}
\end{align}
Hence, the family $\{B_1(t)\}_{t\in [1,\infty)}$ is a RPROB with respect to $\phi(t)=e^{itH_0}\psi(t)$ and by Eqs.~\eqref{CC: ineq2}, ~\eqref{supBphi: 1} and ~\eqref{ineq: g: 1}, we obtain
\begin{align}
   &\int_1^\infty  |(\phi(t), \F_1\partial_t[\F_c]\F_1\phi(t))_{L^2_x(\mathbb{R}^n)}|dt\nonumber\\
   =& \int_1^\infty  (\phi(t), \F_1\partial_t[\F_c]\F_1\phi(t))_{L^2_x(\mathbb{R}^n)}dt\nonumber\\
   \leq &  \sup\limits_{t\in [1,\infty)}| \langle B_1: \phi(t)\rangle_t|+\| g(t)\|_{L^1_t[1,\infty)}\nonumber\\
   \lesssim_{E,n} &{ 1+\|\langle x\rangle^\delta \mathcal{N}(x,t,\psi(t))\|_{L^\infty_tL^2_{x}(\mathbb{R}^{n+1})}}\nonumber\\
    <& \infty.\label{p1: local}
\end{align}
By Cauchy-Schwarz inequality and the non-negativity of $\partial_t[\F_c]$ (see~\eqref{Fcnneg}), $\psi_{p,1}(t)$ satisfies the estimate, for $T_2\geq T_1\geq 1$,
\begin{align}
  &\| \psi_{p,1}(T_2)-\psi_{p,1}(T_1)\|_{L^2_x(\mathbb{R}^n)}\nonumber\\
  \leq & \| \int_{T_1}^{T_2} | \partial_t[\F_c]\F_1 \phi(t)|  dt \|_{L^2_x(\mathbb{R}^n)}\nonumber\\
  \leq & \| \left( \int_{T_1}^{T_2} \partial_t[\F_c]dt\right)^{1/2}  \left( \int_{T_1}^{T_2} \partial_t[\F_c]|\F_1\phi(t)|^2dt\right)^{1/2}\|_{L^2_x(\mathbb{R}^n)}.\label{psipT2T1: 1}
\end{align}
By estimates~\eqref{dtFc} and~\eqref{p1: local}, estimate~\eqref{psipT2T1: 1} leads to 
\begin{align}
    \| \psi_{p,1}(T_2)-\psi_{p,1}(T_1)\|_{L^2_x(\mathbb{R}^n)}\leq & \|   \left( \int_{T_1}^{T_2} \partial_t[\F_c]|\F_1\phi(t)|^2dt\right)^{1/2}\|_{L^2_x(\mathbb{R}^n)}\nonumber\\
    =&\left( \int_{\mathbb{R}^n}\int_{T_1}^{T_2}   \partial_t[\F_c]|\F_1\phi(t)|^2 dtd^nx\right)^{1/2} \nonumber\\
    =& \left( \int_{T_1}^{T_2}  |(\F_1\phi(t), \partial_t[\F_c]\F_1\phi(t))_{L^2_x(\mathbb{R}^n)}|dt\right)^{1/2}\nonumber\\
    =& \left( \int_{T_1}^{T_2}  |(\phi(t), \F_1\partial_t[\F_c]\F_1\phi(t))_{L^2_x(\mathbb{R}^n)}|dt\right)^{1/2}\nonumber\\
    \to & 0
\end{align}
as $T_1\to \infty$, where we have used 
\eq
\int_{\mathbb{R}^n}\int_{T_1}^{T_2}   \partial_t[\F_c]|\F_1\phi(t)|^2 dtd^nx\lesssim (T_2-T_1) \sup\limits_{t\in \mathbb{R}}\| \psi(t)\|^2_{L^2_x(\mathbb{R}^n)}<\infty
\eeq
and Fubini's Theorem to switch the order of integration. Hence, $\{ \psi_{p,1}(t)\}_{t\geq 1}$ is Cauchy in $L^2_x(\mathbb{R}^n)$ and therefore 
{\eq
\lim\limits_{t\to \infty}\psi_{p,1}(t)\text{ exists in } L^2_x(\mathbb{R}^n). \label{exist: p1: local}
\eeq}
Next, we obtain the RPRES for $\psi_{p,2}(t)$. Similarly, by substituting $\mathcal{F}_c$ with $\mathcal{F}_c^2$ and $\mathcal{F}_1$ with $\sqrt{\mathcal{F}_1}$ in the process described above, we have the following RPRES: with (see Eq.~\eqref{Fcnneg: 1: 2})
\eq
\tilde C_2^*\tilde C_2:=\F_c \partial_t[\F_1]\F_c,
\eeq
\begin{align}
     &\int_1^\infty  |(\phi(t), \F_c\partial_t[\F_1]\F_c\phi(t))_{L^2_x(\mathbb{R}^n)}|dt< \infty,\label{p2: local}
\end{align}
which implies, for all $T_2\geq T_1\geq1$, 
\begin{align}
    \| \psi_{p,2}(T_2)-\psi_{p,2}(T_1)\|_{L^2_x(\mathbb{R}^n)}\leq & \|   \left( \int_{T_1}^{T_2} \partial_t[\F_1]|\F_c\phi(t)|^2dt\right)^{1/2}\|_{L^2_x(\mathbb{R}^n)}\nonumber\\
    =&\left( \int_{\mathbb{R}^n}\int_{T_1}^{T_2}   \partial_t[\F_1]|\F_c\phi(t)|^2 dtd^nx\right)^{1/2} \nonumber\\
    =& \left( \int_{T_1}^{T_2}  |(\F_c\phi(t), \partial_t[\F_1]\F_c\phi(t))_{L^2_x(\mathbb{R}^n)}|dt\right)^{1/2}\nonumber\\
    =& \left( \int_{T_1}^{T_2}  |(\phi(t), \F_c\partial_t[\F_1]\F_c\phi(t))_{L^2_x(\mathbb{R}^n)}|dt\right)^{1/2}\nonumber\\
    \to & 0
\end{align}
as $T_1\to \infty$. Hence, $\{ \psi_{p,2}(t)\}_{t\geq 1}$ is Cauchy in $L^2_x(\mathbb{R}^n)$ and  
{\eq
\lim\limits_{t\to \infty}\psi_{p,2}(t)\text{ exists in } L^2_x(\mathbb{R}^n). \label{exist: p2: local}
\eeq}
By Lemma~\ref{com}, with $\alpha>\beta,$
\eq
\| [\F_c,\partial_t[\F_1]]e^{itH_0}\psi(t)\|_{L^2_x(\mathbb{R}^n)}\lesssim_{n}\frac{1}{t^{1+\alpha-\beta}}\in L^1_t[1,\infty),
\eeq
which implies that
\eq
\lim\limits_{t\to \infty}\psi_{c}(t)\text{ exists in } L^2_x(\mathbb{R}^n). \label{exist: c: local}
\eeq
The conclusions \eqref{exist: 0: local},~\eqref{exist: int: local},~\eqref{exist: p1: local},~\eqref{exist: p2: local} and~\eqref{exist: c: local}, together with Eq.~\eqref{Cook: local}, imply that
{\eq
\Omega_{\alpha,\beta}^*\psi(0)\equiv\lim\limits_{t\to \infty}\Omega_{\alpha,\beta}^*\psi(0)\text{ exists in } L^2_x(\mathbb{R}^n). \label{exist: Omega: local}
\eeq}
Take ${\varphi}\in L^2_x(\mathbb{R}^n)$. We also have that for all $(\alpha,\beta), (\alpha',\beta')\in (0,c_1)\times (0,c_2)$ with $\beta<\alpha$ and $\beta'<\alpha'$, by Cauchy-Schwarz inequality, Assumption~\ref{asp: global} and the unitarity of $e^{itH_0}$,
\begin{align}
    &\left|({\varphi}, \Omega_{\alpha,\beta}^*(t)\psi(0)-\Omega_{\alpha',\beta'}^*(t)\psi(0))_{L^2_x(\mathbb{R}^n)}\right|\nonumber\\
    =& \left|( (\F_c(\frac{|x|}{t^\alpha}\leq 1)\F_1(t^\beta|P|>1)-\F_c(\frac{|x|}{t^{\alpha'}}\leq 1))\F_1(t^{\beta'}|P|>1){\varphi}, e^{itH_0}\psi(t))_{L^2_x(\mathbb{R}^n)}\right|\nonumber\\
    \leq & \|(\F_1(t^\beta|P|>1)\F_c(\frac{|x|}{t^\alpha}\leq 1)-\F_1(t^{\beta'}|P|>1)\F_c(\frac{|x|}{t^{\alpha'}}\leq 1)){\varphi} \|_{L^2_x(\mathbb{R}^n)}\| e^{itH_0}\psi(t)\|_{L^2_x(\mathbb{R}^n)}\nonumber\\
    \leq & E\|(\F_c(\frac{|x|}{t^\alpha}\leq 1)-\F_c(\frac{|x|}{t^{\alpha'}}\leq 1)){\varphi} \|_{L^2_x(\mathbb{R}^n)}+\|(\F_1(t^\beta|P|> 1)-\F_1(t^{\beta'}|P|> 1))\F_c(\frac{|x|}{t^{\alpha'}}){\varphi} \|_{L^2_x(\mathbb{R}^n)}\nonumber\\
    \to & 0\label{eq: weak limit equal}
\end{align}
as $t\to \infty$. This implies 
\eq
w\text{-}\lim\limits_{t\to \infty} \Omega_{\alpha,\beta}^*(t)\psi(0)-\Omega_{\alpha',\beta'}^*(t)\psi(0)=0,\quad \text{ in }L^2_x(\mathbb{R}^n)
\eeq
and therefore, due to the existence of $\Omega_{\alpha,\beta}^*\psi(0)$ and $\Omega_{\alpha',\beta'}^*\psi(0)$ in $L^2_x(\mathbb{R}^n)$ in strong sense, Eq.~\eqref{weak id}. 
\end{proof}
Recall that $\psi_{wl}(t)$ is given in Eq.~\eqref{def: weak: local}.
\begin{lemma}\label{lem weak} Let $\psi(t)$ be the solution to system~\eqref{SE}. Let Assumptions~\ref{asp: global} and~\ref{asp: local} be satisfied. Then
\eq\label{lim weak local}
\lim\limits_{t\to \infty} \| \psi(t)- e^{-itH_0}\Omega_{\alpha,\beta}^*\psi(0)-\psi_{wl}(t)\|_{L^2_x(\mathbb{R}^n)}=0.
\eeq

\end{lemma}
{We find that $\F_c( \frac{|x|}{(t+1)^{1/2+\epsilon}}\geq 1)\psi(t)\equiv \psi(t)-\psi_{wl}(t)$ satisfies the decomposition}
\begin{align}
\F_c( \frac{|x|}{(t+1)^{1/2+\epsilon}}\geq 1)\psi(t)= & \psi_1(t)+\psi_2(t)+\psi_3(t),\label{decom: psi largex}
\end{align}
where $\psi_j(t), j=1,2,3,$ are given by 
\eq
\psi_1(t):=\F_c( \frac{|x|}{(t+1)^{1/2+\epsilon}}\geq 1)P^+\F_1(\sqrt{t+1}|P|\geq 1)\psi(t),
\eeq
\eq
\psi_2(t):=\F_c( \frac{|x|}{(t+1)^{1/2+\epsilon}}\geq 1)P^-\psi(t)
\eeq
and
\eq
\psi_3(t):=\F_c( \frac{|x|}{(t+1)^{1/2+\epsilon}}\geq 1)P^+\F_1(\sqrt{t+1}|P|< 1)\psi(t).
\eeq
We approximate $\psi_1(t)$ by $e^{-itH_0}\Omega^*_{\alpha,\beta}\psi(0)$ and arrive at~\eqref{lim weak local} by showing that
\eq
\lim\limits_{t\to \infty}\| \psi_j(t)\|_{L^2_x(\mathbb{R}^n)}=0,\label{lim: psi2}\quad j=2,3,
\eeq
\eq
\lim\limits_{t\to \infty}\| \psi_1(t)-\F_c( \frac{|x|}{(t+1)^{1/2+\epsilon}}\geq 1)P^+\F_1(\sqrt{t+1}|P|\geq 1)e^{-itH_0} \Omega^*_{\alpha,\beta}\psi(0)\|_{L^2_x(\mathbb{R}^n)}=0\label{lim:psi3}
\eeq
and
\eq
\lim\limits_{t\to \infty}\| e^{-itH_0}\Omega^*_{\alpha,\beta}\psi(0) -\F_c( \frac{|x|}{(t+1)^{1/2+\epsilon}}\geq 1)P^+\F_1(\sqrt{t+1}|P|\geq 1)e^{-itH_0} \Omega^*_{\alpha,\beta}\psi(0)\|_{L^2_x(\mathbb{R}^n)}=0.\label{lim:psi4}
\eeq
\begin{proof}[Proof of Lemma~\ref{lem weak}.] We begin with the proof of~\eqref{lim: psi2}. By Duhamel's formula, $\psi_2(t)$ and $\psi_3(t)$ read 
\begin{align}
    \psi_2(t)=& \F_c (\frac{|x|}{(t+1)^{1/2+\epsilon}}\geq 1) P^-e^{-itH_0}\psi(0)\nonumber\\
    &+(-i)\int_0^t \F_c (\frac{|x|}{(t+1)^{1/2+\epsilon}}\geq 1) P^-e^{-i(t-s)H_0}\mathcal{N}(x,s,\psi(s)) \psi(s)ds\nonumber\\
    =: & \psi_{21}(t)+\psi_{22}(t)\label{psi21psi22}
\end{align}
and
\begin{align}
    \psi_3(t)=& \F_c (\frac{|x|}{(t+1)^{1/2+\epsilon}}\geq 1) P^+\F_1(\sqrt{t+1}|P|<1)e^{-itH_0}\psi(0)\nonumber\\
    &+(-i)\int_0^t \F_c (\frac{|x|}{(t+1)^{1/2+\epsilon}}\geq 1) P^+\F_1(\sqrt{t+1}|P|<1)e^{-i(t-s)H_0}\mathcal{N}(x,s,\psi(s)) \psi(s)ds\nonumber\\
    =: & \psi_{31}(t)+\psi_{32}(t),\label{psi31psi32}
\end{align}
respectively. By Lemma~\ref{lem: Ppmf}, $\psi_{21}(t)$ satisfies 
\begin{align}
\| \psi_{21}(t)\|_{L^2_x(\mathbb{R}^n)}\leq& C\| P^-e^{-itH_0}\psi(0)\|_{L^2_x(\mathbb{R}^n)}\to 0\label{psi21(t)}
\end{align}
as $t\to \infty$. By Lemma~\ref{Lem: Pprop} and Assumptions~\ref{asp: global} and~\ref{asp: local}, $\psi_{22}(t)$ satisfies, with $\delta>2$, 
\begin{align}
    \| \psi_{22}(t)\|_{L^2_x(\mathbb{R}^n)}\leq & \int_0^t \| \F_c (\frac{|x|}{(t+1)^{1/2+\epsilon}}\geq 1) P^-e^{-i(t-s)H_0}\F_1(\sqrt{t+1}|P|\geq 1)\langle x\rangle^{-\delta}\|_{L^2_x(\mathbb{R}^n)\to L^2_x(\mathbb{R}^n)}\nonumber\\
    &\times \| \langle x\rangle^\delta \mathcal{N}(x,s,\psi(s))\psi(s)\|_{L^2_x(\mathbb{R}^n)}ds\nonumber\\
    \lesssim_\epsilon & \int_0^t \frac{1}{\langle  (t+1)^{1/2+\epsilon}+s/\sqrt{t+1}\rangle^{\delta}} \sup\limits_{u\in \mathbb{R}} \| \langle x\rangle^\delta \mathcal{N}(x,u,\psi(u))\psi(u)\|_{L^2_x(\mathbb{R}^n)}ds \nonumber\\
    \lesssim_{\epsilon,\delta} & \frac{t}{\langle t+1\rangle^{1+2\epsilon}} \sup\limits_{u\in \mathbb{R}}\| \langle x\rangle^\delta \mathcal{N}(x,u,\psi(u))\psi(u)\|_{L^2_x(\mathbb{R}^n)}\nonumber\\
    \to & 0\label{psi22}
\end{align}
as $t\to \infty.$ {The estimates \eqref{psi21(t)} and~\eqref{psi22}}, together with Eq.~\eqref{psi21psi22}, imply 
\eq
\lim\limits_{t\to \infty} \| \psi_2(t)\|_{L^2_x(\mathbb{R}^n)}=0.\label{lim psi20}
\eeq
By 
\eq
s\text{-}\lim\limits_{t\to \infty} \F_1(\sqrt{t+1}|P|<1)=0,\quad \text{ on }L^2_x(\mathbb{R}^n),\label{F1}
\eeq
$\psi_{31}(t)$ satisfies 
\begin{align}
\limsup\limits_{t\to \infty} \| \psi_{31}(t)\|_{L^2_x(\mathbb{R}^n)} \leq & C\limsup\limits_{t\to \infty} \| \F_1(\sqrt{t+1}|P|<1)\psi(0)\|_{L^2_x(\mathbb{R}^n)}= 0.\label{psi31}
\end{align}
By Lemma~\ref{lem: Ppropfree} and Assumptions~\ref{asp: global} and~\ref{asp: local}, $\psi_{32}(t)$ satisfies, with $\delta>2$, 
\begin{align}
    \| \psi_{32}(t)\|_{L^2_x(\mathbb{R}^n)}\leq & \int_0^t C\| \F_c (\frac{|x|}{(t+1)^{1/2+\epsilon}}\geq 1) P^+e^{-i(t-s)H_0}\F_1(\sqrt{t+1}|P|<1)\langle x\rangle^{-\delta}\|_{L^2_x(\mathbb{R}^n)\to L^2_x(\mathbb{R}^n)}\nonumber\\
    &\times \| \langle x\rangle^\delta \mathcal{N}(x,s,\psi(s))\psi(s)\|_{L^2_x(\mathbb{R}^n)}ds\nonumber\\
    \lesssim_\epsilon & \int_0^t \frac{1}{ (t+1)^{\frac{1}{2}\delta+\epsilon\delta}} \sup\limits_{u\in \mathbb{R}} \| \langle x\rangle^\delta \mathcal{N}(x,u,\psi(u))\psi(u)\|_{L^2_x(\mathbb{R}^n)}ds \nonumber\\
    \lesssim_\epsilon & \frac{t}{\langle t+1\rangle^{1+2\epsilon}} \sup\limits_{u\in \mathbb{R}}\| \langle x\rangle^\delta \mathcal{N}(x,u,\psi(u))\psi(u)\|_{L^2_x(\mathbb{R}^n)} \to  0\label{psi32}
\end{align}
as $t\to \infty.$ {The estimates \eqref{psi31} and~\eqref{psi32}}, together with Eq.~\eqref{psi31psi32}, imply 
\eq
\lim\limits_{t\to \infty} \| \psi_3(t)\|_{L^2_x(\mathbb{R}^n)}=0.\label{lim psi3}
\eeq
Eqs.~\eqref{lim psi20} and ~\eqref{lim psi3} imply Eq.~\eqref{lim: psi2}.

Next, we prove~\eqref{lim:psi3}. Let 
\begin{align}
\Omega^*\psi(0):=\psi(0)+w\text{-}\lim\limits_{t\to \infty}(-i)\int_0^t e^{isH_0}\mathcal{N}(x,s,\psi(s))\psi(s)ds\qquad \text{in }L^2_x(\mathbb{R}^n).\label{eq: Omega weak}
\end{align}
{By Cook's method, $\Omega^*\psi(0)$ exists in $L^2_x(\mathbb{R}^n)$.} {Indeed, the weak limit define in~\eqref{eq: Omega weak} exists since for any $c>0, \delta>1$ and each $\varphi\in C_0^\infty(\mathbb R^n)$  with $\varphi=F_1(|P|>c)\varphi$,
\[
\begin{aligned}
 & \int_0^t  |(\varphi, e^{isH_0}\mathcal N(x,s,\psi(s))\psi(s))_{L^2_x}|ds\\
 \leq& \int_0^t\|\langle x\rangle^{-\delta}e^{-isH_0}\varphi\|_{L^2_x(\mathbb R^n)}\| \langle x\rangle^\delta\mathcal N(x,s,\psi(s))\|_{L^\infty_{x,s}(\mathbb R^{n+1})}\|\psi(s)\|_{L^2_x(\mathbb R^n)}ds\\
 \lesssim_c& \int_0^t \frac{1}{\langle s\rangle^\delta}\|\langle x\rangle^\delta \varphi\|_{L^2_x(\mathbb R^n)}\| \langle x\rangle^\delta\mathcal N(x,s,\psi(s))\|_{L^\infty_{x,s}(\mathbb R^{n+1})}\|\psi(s)\|_{L^2_x(\mathbb R^n)}ds<\infty.
\end{aligned}
\]}{In addition, following a similar process of~\eqref{eq: weak limit equal}, it is worth noting that 
\[
w\text{-}\lim\limits_{t\to \infty}\left(1-\mathcal F_c(\frac{|x|}{t^\alpha}\leq 1)\mathcal F_1(t^\beta|P|>1)\right)e^{itH_0}\psi(t)=0\qquad \text{ in }L^2_x(\mathbb R^n)
\]
holds: for all $\varphi \in L^2_x(\mathbb R^n)$,
\begin{align}
    &\left|({\varphi}, \left(1-\mathcal F_c(\frac{|x|}{t^\alpha}\leq 1)\mathcal F_1(t^\beta|P|>1)\right)e^{itH_0}\psi(t))_{L^2_x(\mathbb{R}^n)}\right|\nonumber\\
    =& \left|( \left(1-\mathcal F_c(\frac{|x|}{t^\alpha}\leq 1)\mathcal F_1(t^\beta|P|>1)\right)\varphi, e^{itH_0}\psi(t))_{L^2_x(\mathbb{R}^n)}\right|\nonumber\\
    \leq & \|\left(1-\mathcal F_c(\frac{|x|}{t^\alpha}\leq 1)\mathcal F_1(t^\beta|P|>1)\right){\varphi} \|_{L^2_x(\mathbb{R}^n)}\| e^{itH_0}\psi(t)\|_{L^2_x(\mathbb{R}^n)}\nonumber\\
    \to & 0.\nonumber
\end{align}}
By Eqs.~\eqref{wave-1} and ~\eqref{weak id}, {this, {together with the existence of $\Omega^*\psi(0)$}, implies $\Omega^*\psi(0)=\Omega_{\alpha,\beta}^*\psi(0)$ for any $\alpha\in (0,c_1)$ and $\beta\in (0,c_2)$.} This is because the channel wave operator exists in the strong {sense}, and on the complement support, where $|x|
\ge t^{\alpha},$ the weak limit exists and is equal to zero. By Duhamel's expansion and $\Omega^*\psi(0)=\Omega_{\alpha,\beta}^*\psi(0)$, we have 
\begin{align}
    & \F_c( \frac{|x|}{(t+1)^{1/2+\epsilon}}\geq 1)P^+\F_1(\sqrt{t+1}|P|>1)e^{-itH_0} \Omega^*_{\alpha,\beta}\psi(0)- \psi_1(t)\nonumber\\
    =&(-i)\F_c(\frac{|x|}{(t+1)^{1/2+\epsilon}}\geq 1) P^+\int_t^\infty \F_1(\sqrt{t+1}|P|>1) e^{i(s-t)H_0}\mathcal{N}(x,s,\psi(s))\psi(s)ds.\label{id: psi1}
\end{align}

By Lemma~\ref{Lem: Pprop} and Assumptions~\ref{asp: global} and~\ref{asp: local}, \eqref{id: psi1} implies, with $\delta>2$, 
\begin{align}
    & \|\F_c( \frac{|x|}{(t+1)^{1/2+\epsilon}}\geq 1)P^+ \F_1(\sqrt{t+1}|P|>1)e^{-itH_0} \Omega^*_{\alpha,\beta}\psi(0)- \psi_1(t)\|_{L^2_x(\mathbb{R}^n)}\nonumber\\
    \leq &\int_t^\infty\|\F_c(\frac{|x|}{(t+1)^{1/2+\epsilon}}\geq 1) P^+ \F_1(\sqrt{t+1}|P|>1)e^{i(s-t)H_0}\langle x\rangle^{-\delta}\|_{L^2_x(\mathbb{R}^n)\to L^2_x(\mathbb{R}^n)}\nonumber\\
    &\times \|\langle x\rangle^\delta\mathcal{N}(x,s,\psi(s))\psi(s)\|_{L^2_x(\mathbb{R}^n)}ds\nonumber\\
    \lesssim_\epsilon& \int_t^\infty \frac{1}{\langle (t+1)^{1/2+\epsilon}+(s-t)/\sqrt{t+1}\rangle^\delta}\|\langle x\rangle^\delta\mathcal{N}(x,s,\psi(s))\psi(s)\|_{L^2_x(\mathbb{R}^n)}ds\nonumber\\
    \lesssim_{\epsilon,\delta} & \frac{\sqrt{t+1}}{\langle t+1\rangle^{\frac{1}{2}(\delta-1)+\epsilon(\delta -1)}}{\|\langle x\rangle^\delta\mathcal{N}(x,t,\psi(t))\psi(t)\|_{L^\infty_tL^2_x(\mathbb{R}^{n+1})}}\nonumber\\
    \to & 0\label{lim psi1(t): local}
\end{align}
as $t\to \infty.$ {Estimate \eqref{lim psi1(t): local}} implies ~\eqref{lim:psi3}.

Now we prove~\eqref{lim:psi4}. Eq.~\eqref{F1}, together with Lemmas~\ref{lem: Ppmf} and~\ref{lem: not linear: local}, implies
\begin{align}
    & \|e^{-itH_0}\Omega^*_{\alpha,\beta}\psi(0) -\F_c( \frac{|x|}{(t+1)^{1/2+\epsilon}}\geq 1)P^+\F_1(\sqrt{t+1}|P|\geq 1)e^{-itH_0} \Omega^*_{\alpha,\beta}\psi(0)\|_{L^2_x(\mathbb{R}^n)}\nonumber\\
    \lesssim & \| \F_c( \frac{|x|}{(t+1)^{1/2+\epsilon}}\geq 1)P^+\F_1(\sqrt{t+1}|P|< 1)e^{-itH_0} \Omega^*_{\alpha,\beta}\psi(0) \|_{L^2_x(\mathbb{R}^n)}\nonumber\\
    & +\| \F_c( \frac{|x|}{(t+1)^{1/2+\epsilon}}< 1)P^+e^{-itH_0} \Omega^*_{\alpha,\beta}\psi(0) \|_{L^2_x(\mathbb{R}^n)}+\| P^-e^{-itH_0} \Omega^*_{\alpha,\beta}\psi(0) \|_{L^2_x(\mathbb{R}^n)}\nonumber\\
    \lesssim & \|\F_1(\sqrt{t+1}|P|< 1)e^{-itH_0} \Omega^*_{\alpha,\beta}\psi(0) \|_{L^2_x(\mathbb{R}^n)}+\| \F_c( \frac{|x|}{(t+1)^{1/2+\epsilon}}< 1)P^+e^{-itH_0} \Omega^*_{\alpha,\beta}\psi(0) \|_{L^2_x(\mathbb{R}^n)}\nonumber\\
    & +\| P^-e^{-itH_0} \Omega^*_{\alpha,\beta}\psi(0) \|_{L^2_x(\mathbb{R}^n)}\to  0\label{psi(3): local}
\end{align}
as $t\to \infty.$ {Equation} \eqref{psi(3): local} implies~\eqref{lim:psi4}. By~\eqref{lim: psi2},~\eqref{lim:psi3},~\eqref{lim:psi4} and Eq.~\eqref{decom: psi largex}, together with Eq.~\eqref{def: weak: local}, we arrive at
\begin{align}
& \| \psi(t)-  \psi_{wl}(t)-e^{-itH_0}\Omega_{\alpha,\beta}^*\psi(0) \|_{L^2_x(\mathbb{R}^n)}\nonumber\\
   = & \| \psi(t)-  \F_c( \frac{|x|}{(t+1)^{1/2+\epsilon}}< 1)\psi(t)-e^{-itH_0}\Omega_{\alpha,\beta}^*\psi(0) \|_{L^2_x(\mathbb{R}^n)}\nonumber\\
    \leq & \|  \F_c( \frac{|x|}{(t+1)^{1/2+\epsilon}}\geq 1)\psi(t)-e^{-itH_0}\Omega_{\alpha,\beta}^*\psi(0) \|_{L^2_x(\mathbb{R}^n)}\nonumber\\
    \leq &\| \psi_1(t)-\F_c( \frac{|x|}{(t+1)^{1/2+\epsilon}}\geq 1)P^+\F_1(\sqrt{t+1}|P|\geq 1)e^{-itH_0} \Omega^*_{\alpha,\beta}\psi(0)\|_{L^2_x(\mathbb{R}^n)}\nonumber\\
    &+\| e^{-itH_0}\Omega^*_{\alpha,\beta}\psi(0) - \F_c( \frac{|x|}{(t+1)^{1/2+\epsilon}}\geq 1)P^+\F_1(\sqrt{t+1}|P|\geq 1)e^{-itH_0} \Omega^*_{\alpha,\beta}\psi(0)\|_{L^2_x(\mathbb{R}^n)}\nonumber\\
    &+\| \psi_2(t)\|_{L^2_x(\mathbb{R}^n)}+\| \psi_3(t)\|_{L^2_x(\mathbb{R}^n)} \to  0,
\end{align}
as $t\to \infty.$

\end{proof}
\begin{remark}Here, it is important to emphasize that the integral
\[
\phi(t) := \int_t^\infty e^{isH_0}\, V\, \mathcal{N}(x,s,\psi(s))\, \psi(s)\, ds
\]
is well-defined in the \emph{weak} sense, namely,
\[
\langle f, \phi(t) \rangle = \int_t^\infty \langle f, e^{isH_0} V \mathcal{N}(x,s,\psi(s)) \psi(s) \rangle\, ds
\qquad \text{for all } f \in C_0^\infty(\mathbb R^n).
\]
Furthermore, the estimate \eqref{lim psi1(t): local} implies that, after applying the localization operator
\[
\mathcal{L}_t := \mathcal{F}_c\!\left( \frac{|x|}{(t+1)^{1/2+\epsilon}} \right) P^+ \mathcal{F}_1\!\left( \sqrt{t+1}\, |P| > 1 \right),
\]
the truncated expression
\[
\mathcal{L}_t\, \phi(t)
= 
\mathcal{F}_c\!\left( \frac{|x|}{(t+1)^{1/2+\epsilon}} \right)
P^+\,
\mathcal{F}_1\!\left( \sqrt{t+1}\, |P| > 1 \right)
\phi(t)
\]
admits a \emph{strong} interpretation. In particular, the integral defining $\mathcal{L}_t\, \phi(t)$ converges in the strong sense. In addition, the limit of a weakly convergent sequence is always smaller than or equal to the limit inferior of the norms of the sequence.\end{remark}

\begin{proof}[Proof of Theorem \ref{thm}] It follows from Proposition~\ref{prof: free: local}, Lemma~\ref{lem weak} and Eq.~\eqref{def: weak: local}.
    
\end{proof}

\section{Charge Transfer Potentials}
In this section, we prove Theorem~\ref{thm3}. The proof of Theorem~\ref{thm3} requires the use of the  following proposition and lemmas. Let 
\eq\label{def: psij(t)}
\psi_j(t):=(-i)\int_0^t  e^{i(-t+s)H_0} V_j(x-sv_j,s)\psi(s) ds,\quad j=1,\cdots,N.
\eeq
Then $\psi_j(t),j=1,\cdots,N$, satisfy the differential equations
\eq
i\partial_t\psi_j(t)=H_0\psi_j(t)+V_j(x-tv_j,t)\psi(t),\qquad t>0.\label{df: psij}
\eeq
\begin{proposition}\label{sec: charge: psij} {If all the assumptions of Theorem~\ref{thm3} are satisfied, then for every $j=1,\cdots,N$, $\psi_j(t)$ is uniformly bounded in $L^2_x(\mathbb{R}^n)$ for all $t\in [0,\infty)$:
\eq\label{lim: psijuniform}
\sup\limits_{t\geq 0}\|\psi_j(t) \|_{L^2_x(\mathbb{R}^n)}\lesssim_E 1.
\eeq}
\end{proposition}
The proof of Proposition~\ref{sec: charge: psij} requires Lemmas~\ref{sec: charge: psij1} and~\ref{sec: charge: psij2} listed below.
\begin{lemma}\label{sec: charge: psij1} Let $v_j\in \mathbb{R}^n, j=1,\cdots,N$, be $N$ {non-congruent }vectors: for some $\epsilon>0$,
{$|v_j-v_{j'}|\geq \epsilon, \forall  j\neq j'.$}
For all $j,j'\in \{1,\cdots,N\}$ with $j\neq j'$, and all {$\sigma > n/2+2$}, the estimate 
\eq
\|\langle x\rangle^{-\sigma}\langle P\rangle^{-2} e^{i((s_1-s_2)H_0-s_1 v_j\cdot P+s_2 v_{j'}\cdot P)}\langle x\rangle^{-\sigma}\|_{L^2_x(\mathbb{R}^n)}\lesssim_{ n/2+2,|v_j|,|v_j-v_{j'}|} \frac{1}{\langle s_1-s_2\rangle^{n/2-2}}\frac{1}{\langle s_1\rangle}\frac{1}{\langle s_2\rangle}\label{charge: free: estimate}
\eeq
holds.
\end{lemma}
\begin{proof} {Let 
\begin{equation}
    \mathcal{A}:=\langle x\rangle^{-\sigma}\langle P\rangle^{-2}
    e^{\,i((s_1-s_2)H_0-s_1 v_j\cdot P+s_2 v_{j'}\cdot P)}\langle x\rangle^{-\sigma}.
\end{equation}
We estimate the $L^2$ operator norm of $\mathcal A$ in three cases.

\medskip
\noindent\textbf{Case 1:} $\max\{|s_1|,|s_2|\}\leq 10$.  
Then $|s_1-s_2|\leq 20$. Since 
$e^{i((s_1-s_2)H_0-s_1 v_j\cdot P+s_2 v_{j'}\cdot P)}$ is unitary,
\begin{equation}
   \|\mathcal A\|_{L^2_x(\mathbb R^n)\to L^2_x(\mathbb R^n)}
   \leq  \|\langle x\rangle^{-\sigma}\langle P\rangle^{-2}\|_{L^2_x(\mathbb R^n)\to L^2_x(\mathbb R^n)}
   \|\langle x\rangle^{-\sigma}\|_{L^2_x(\mathbb R^n)\to L^2_x(\mathbb R^n)}
   \lesssim 1.
\end{equation}

\medskip
\noindent\textbf{Case 2:} $|s_1-s_2|\geq \frac{1}{10}\max\{|s_1|,|s_2|\}>1$.  
In this regime we have 
$|s_2|\leq |s_1-s_2|+|s_1|\leq 11|s_1-s_2|$.
Using the $L^\infty$ decay estimate of the free Schr\"odinger propagator, we obtain
\begin{equation}
    \begin{aligned}
        \|\mathcal A\|_{L^2_x(\mathbb R^n)\to L^2_x(\mathbb R^n)}
        &\leq 
        \|\langle x\rangle^{-\sigma}\langle P\rangle^{-2}\|_{L^\infty_x(\mathbb R^n)\to L^2_x(\mathbb R^n)}
        \|e^{i(s_1-s_2)H_0}\|_{L^1_x(\mathbb R^n)\to L^\infty_x(\mathbb R^n)}\\
        &\quad\times
        \|e^{i(-s_1v_j\cdot P+s_2v_{j'}\cdot P)}\|_{L^1_x(\mathbb R^n)\to L^1_x(\mathbb R^n)}
        \|\langle x\rangle^{-\sigma}\|_{L^2_x(\mathbb R^n)\to L^1_x(\mathbb R^n)}\\
        &\lesssim \frac{1}{|s_1-s_2|^{n/2}}\\
        &\lesssim \frac{1}{\langle s_1-s_2\rangle^{n/2-2}}
        \frac{1}{\langle s_1\rangle}
        \frac{1}{\langle s_2\rangle}.
    \end{aligned}
\end{equation}

\medskip
\noindent\textbf{Case 3:} $\max\{|s_1|,|s_2|\}>10|s_1-s_2|$ and $\max\{|s_1|,|s_2|\}>10$.  
Then 
\[
\min\{|s_1|,|s_2|\}\geq \max\{|s_1|,|s_2|\}-|s_1-s_2|>\tfrac{9}{10}\max\{|s_1|,|s_2|\}.
\]
In Fourier space, the phase function can be written as
\begin{equation}
    f(q)=(s_1-s_2)|q|^2-(s_1-s_2)v_j\cdot q+s_2(v_{j'}-v_j)\cdot q.
\end{equation}
Let $e_1:=(v_{j'}-v_j)/|v_{j'}-v_j|$ and $q_1:=e_1\cdot q$.  
Using
\begin{equation}
   e^{is_2|v_{j'}-v_j|q_1}
   =\frac{1}{is_2|v_{j'}-v_j|}\partial_{q_1}\big[e^{is_2|v_{j'}-v_j|q_1}\big],
\end{equation}
we integrate by parts twice in $q_1$ to gain $(s_2|v_{j'}-v_j|)^{-2}$ decay, and then, for $|s_1-s_2|\geq 1$, use $L^\infty$ decay of the free flow to obtain
\begin{equation}
    \begin{aligned}
        \|\mathcal A\|_{L^2_x(\mathbb R^n)\to L^2_x(\mathbb R^n)}
        &\leq 
        \chi(|s_1-s_2|\geq 1)\|\langle x\rangle^{-\sigma}\langle P\rangle^{-2}\langle x\rangle^2\|_{L^\infty_x(\mathbb R^n)\to L^2_x(\mathbb R^n)}\|\langle x\rangle^{2-\sigma}\|_{L^2_x(\mathbb R^n)\to L^1_x(\mathbb R^n)}\\
        &\times \frac{(1+|v_j|+|v_j|^2)}{\langle s_1-s_2\rangle^{n/2}(s_2|v_j-v_{j'}|)^2}+ \frac{\chi(|s_1-s_2|<1)(1+|v_j|+|v_j|^2)}{\langle s_1-s_2\rangle^{n/2}(s_2|v_j-v_{j'}|)^2}\\
        &\lesssim_{ n/2+2,|v_j|,|v_j-v_{j'}|} \frac{1}{\langle s_1-s_2\rangle^{n/2-2}}
        \frac{1}{\langle s_1\rangle}
        \frac{1}{\langle s_2\rangle}.
    \end{aligned}
\end{equation}

\medskip
Combining the three cases yields the desired bound
\begin{equation}
 \|\mathcal A\|_{L^2_x(\mathbb R^n)\to L^2_x(\mathbb R^n)}
 \lesssim_{ n/2+2,|v_j|,|v_j-v_{j'}|}
 \frac{1}{\langle s_1-s_2\rangle^{n/2-2}}
 \frac{1}{\langle s_1\rangle}
 \frac{1}{\langle s_2\rangle}.
\end{equation}
This completes the proof.}\end{proof}
\begin{remark}
    To prove estimate~\eqref{charge: free: estimate}, one may alternatively establish a pointwise kernel bound for the inner operator
\begin{equation}
    K(x,y)
    :=\langle x\rangle^{-2}\langle y\rangle^{-2}
    \int_{\mathbb R^n}
    e^{\,i(x-y)\cdot \xi}
    e^{\,i\big((s_1-s_2)|\xi|^2+(s_2v_{j'}-s_1v_j)\cdot \xi\big)}
    \langle \xi\rangle^{-2}\,d\xi,
\end{equation}
by applying the standard stationary/nonstationary phase analysis
(see~\cite[Section 1.1]{Sogge2017}).

\end{remark}

\begin{lemma}\label{sec: charge: psij2} Under the assumptions of Theorem~\ref{thm3}, for all $j,l\in \{1,\cdots,N\}$ with $j\neq l$,
\eq
\left| (\psi_j(t), \psi_{l}(t))_{L^2_x(\mathbb{R}^n)}\right|\lesssim_{E,|v_j|,|v_j-v_l|} 1\quad t\geq 0.
\eeq
\end{lemma}
\begin{proof} Take $j,l\in \{1,\cdots,N\}$ with $j<l$. By Cauchy-Schwarz inequality, the unitarity of $e^{i(t-s_k)H_0}, k=1,2,$ and Assumptions~\ref{asp: global} and~\ref{asp: charge}, we have
\begin{align}
&\int_0^t \int_0^t\int_{\mathbb{R}^n}|  e^{-i(t-s_1)H_0}V_{j}(x-s_1v_j,s_1) \psi(s_1)||e^{-i(t-s_2)H_0}V_{l}(x-s_2v_l,s_2) \psi(s_2)|dx ds_1ds_2\nonumber\\
\leq &  \int_0^t \int_0^t\|e^{-i(t-s_1)H_0}V_{j}(x-s_1v_j,s_1) \psi(s_1)\|_{L^2_x(\mathbb{R}^n)}\|e^{-i(t-s_2)H_0}V_{l}(x-s_2v_l,s_2) \psi(s_2)\|_{L^2_x(\mathbb{R}^n)}ds_1ds_2\nonumber\\
= &  \int_0^t \int_0^t\|V_{j}(x-s_1v_j,s_1) \psi(s_1)\|_{L^2_x(\mathbb{R}^n)}\|V_{l}(x-s_2v_l,s_2) \psi(s_2)\|_{L^2_x(\mathbb{R}^n)}ds_1ds_2\nonumber\\
\leq & t^2\|V_j(x,s)\|_{L^\infty_{x,s}(\mathbb{R}^n)}\|V_l(x,s)\|_{L^\infty_{x,s}(\mathbb{R}^n)} E^2<\infty.\label{Lem5.3: eq1}
\end{align}
By Eq.~\eqref{Lem5.3: eq1} and Fubini's Theorem, $(\psi_j(t),\psi_l(t))_{L^2_x(\mathbb{R}^n)}$ reads
\begin{align}
    &(\psi_j(t),\psi_l(t))_{L^2_x(\mathbb{R}^n)}\nonumber\\
    =& ( \int_0^t e^{-i(t-s_1)H_0}V_{j}(x-s_1v_j,s_1) \psi(s_1)ds_1, \int_0^t e^{-i(t-s_2)H_0}V_{l}(x-s_2v_l,s_2) \psi(s_2)ds_2)_{L^2_x(\mathbb{R}^n)}\nonumber\\
    =& \int_0^t \int_0^t  ( V_{j}(x-s_1v_j,s_1) \psi(s_1),e^{i((s_2-s_1)H_0)}V_{l}(x-s_2v_l,s_2) \psi(s_2))_{L^2_x(\mathbb{R}^n)} ds_1ds_2\nonumber\\
    =& \int_0^t \int_0^t  ( V_{j}(x,s_1)e^{is_1 P\cdot v_j} \psi(s_1),e^{i((s_2-s_1)H_0+s_1v_j\cdot P-s_2v_l\cdot P)}V_{l}(x,s_2) e^{is_2v_l\cdot P}\psi(s_2))_{L^2_x(\mathbb{R}^n)} ds_1ds_2,
\end{align}
where we have used the unitarity of $e^{-i(t-s_1)H_0}$ and the equations
\eq
V_j(x-s_1v_j, s_1)=e^{-is_1 v_j\cdot P}V_j(x, s_1)e^{is_1 v_j\cdot P}
\eeq
and
\eq
V_l(x-s_2v_l, s_2)=e^{-is_2 v_l\cdot P}V_l(x, s_2)e^{is_2 v_l\cdot P}.
\eeq
By Lemma~\ref{Lem: B2}, {with $\delta$ as in Assumption~\ref{asp: charge}} we have
\begin{align}
    \| \langle x\rangle^\delta \langle P\rangle V_j(x,s_1)e^{is_1 P\cdot v_j}\psi(s_1)\|_{L^2_x(\mathbb{R}^n)}\lesssim_\delta& \|\langle x\rangle^\delta  V_j(x,s_1)e^{is_1 P\cdot v_j}\psi(s_1) \|_{\mathcal H^1_x}\nonumber\\
    \lesssim_E& \sup\limits_{t\in \mathbb{R}} \| \langle x\rangle^\delta  V_j(x,t)\|_{W^{1,\infty}_x(\mathbb{R}^n)}\label{Beq1}
\end{align}
and
\begin{align}
    \| \langle x\rangle^\delta \langle P\rangle V_l(x,s_2)e^{is_2 P\cdot v_l}\psi(s_2)\|_{L^2_x(\mathbb{R}^n)}\lesssim_\delta& \|\langle x\rangle^\delta  V_l(x,s_2)e^{is_2 P\cdot v_l}\psi(s_2) \|_{\mathcal H^1_x}\nonumber\\
    \lesssim_E& \sup\limits_{t\in \mathbb{R}} \| \langle x\rangle^\delta  V_l(x,t)\|_{W^{1,\infty}_x(\mathbb{R}^n)}.\label{Beq2}
\end{align}
By Lemma~\ref{sec: charge: psij1}, Cauchy-Schwarz inequality and the estimates~\eqref{Beq1} and~\eqref{Beq2}, together with assumptions of Theorem~\ref{thm3}, we have, for $n\geq 5$, 
\begin{align}
    | (\psi_j(t),\psi_l(t))_{L^2_x(\mathbb{R}^n)}|\leq & \int_0^t \int_0^t \| \langle x\rangle^\delta \langle P\rangle V_j(x,s_1)e^{is_1 P\cdot v_j}\psi(s_1)\|_{L^2_x(\mathbb{R}^n)}\nonumber\\
    &\qquad\qquad\qquad\times\| \langle x\rangle^{-\delta} \langle P\rangle^{-2} e^{i((s_2-s_1)H_0+s_1v_j\cdot P-s_2v_l\cdot P)}\langle x\rangle^{-\delta}\|_{L^2_x(\mathbb{R}^n)\to L^2_x(\mathbb{R}^n)}\nonumber\\
    &\qquad\qquad\qquad\qquad\times\| \langle x\rangle^\delta \langle P\rangle V_l(x,s_2)e^{is_2 P\cdot v_l}\psi(s_2)\|_{L^2_x(\mathbb{R}^n)}ds_1ds_2\nonumber\\
    \lesssim_{E,|v_j|,|v_l-v_j|}& \sup\limits_{t\in \mathbb{R}} \| \langle x\rangle^\delta  V_l(x,t)\|_{W^{1,\infty}_x(\mathbb{R}^n)}\sup\limits_{t\in \mathbb{R}} \| \langle x\rangle^\delta  V_j(x,t)\|_{W^{1,\infty}_x(\mathbb{R}^n)}\nonumber\\
    &\times\int_0^t \int_0^t \frac{1}{\langle s_1-s_2\rangle^{n/2-2}}\frac{1}{\langle s_1\rangle}\frac{1}{\langle s_2\rangle}ds_1ds_2\nonumber\\
    \lesssim_{E,|v_j|,|v_l-v_j|}&\sup\limits_{t\in \mathbb{R}} \| \langle x\rangle^\delta  V_l(x,t)\|_{W^{1,\infty}_x(\mathbb{R}^n)}\sup\limits_{t\in \mathbb{R}} \| \langle x\rangle^\delta  V_j(x,t)\|_{W^{1,\infty}_x(\mathbb{R}^n)}.
\end{align}

\end{proof}
\begin{proof}[Proof of Proposition~\ref{sec: charge: psij}] By Duhamel's formula, $\psi(t)-e^{-itH_0}\psi(0)$ reads 
\eq
\psi(t)-e^{-itH_0}\psi(0)=\sum\limits_{j=1}^N \psi_j(t).\label{eq: psi(t): charge}
\eeq
Expanding $(\psi(t)-e^{-itH_0}\psi(0),\psi(t)-e^{-itH_0}\psi(0))_{L^2_x(\mathbb{R}^n)}$ by Eq.~\eqref{eq: psi(t): charge}, we obtain 
\begin{align}
    &(\psi(t)-e^{-itH_0}\psi(0),\psi(t)-e^{-itH_0}\psi(0))_{L^2_x(\mathbb{R}^n)}\nonumber\\
    =&\sum\limits_{j=1}^N (\psi_j(t),\psi_j(t))_{L^2_x(\mathbb{R}^n)}+\sum\limits_{1\leq j<l\leq N} 2\text{Re}(\psi_j(t),\psi_l(t))_{L^2_x(\mathbb{R}^n)},
\end{align}
which is equivalent to 
\begin{align}
    &\sum\limits_{j=1}^N (\psi_j(t),\psi_j(t))_{L^2_x(\mathbb{R}^n)}\nonumber\\
    =&(\psi(t)-e^{-itH_0}\psi(0),\psi(t)-e^{-itH_0}\psi(0))_{L^2_x(\mathbb{R}^n)}-\sum\limits_{1\leq j<l\leq N} 2\text{Re}(\psi_j(t),\psi_l(t))_{L^2_x(\mathbb{R}^n)}.\label{expand: charge}
\end{align}
By Cauchy-Schwarz inequality, the unitarity of $e^{-itH_0}$ and Assumption~\ref{asp: global}, Eq.~\eqref{expand: charge}, together with Lemma~\ref{sec: charge: psij2}, implies, 
\begin{align}
     \sum\limits_{j=1}^N (\psi_j(t),\psi_j(t))_{L^2_x(\mathbb{R}^n)}\leq & 2\| \psi(t)\|^2_{L^2_x(\mathbb{R}^n)}+2\| \psi(0)\|^2_{L^2_x(\mathbb{R}^n)}\nonumber\\
     &+2\sum\limits_{1\leq j<l\leq N} \left| (\psi_j(t),\psi_l(t))_{L^2_x(\mathbb{R}^n)}\right|\lesssim_E  1+N^2\lesssim_E  1,\qquad t\geq 1.
\end{align}
This together with the positivity of $(\psi_j(t),\psi_j(t))_{L^2_x(\mathbb{R}^n)}, j=1,\cdots,N$, yields estimate~\eqref{lim: psijuniform}.

\end{proof}

\begin{proposition}\label{prop: charge} Let assumptions in Theorem~\ref{thm3} be satisfied. Then for every $j=1,\cdots,N$, $\epsilon\in (0,1/2)$ and $\alpha\in (0,1-2/n), n\geq 5$,
\begin{enumerate}
    \item \eq\label{def: psijfree}
\psi_{j,+}(x):=s\text{-}\lim\limits_{t\to \infty} e^{itH_0}\F_c(\frac{|x-2tP|}{t^\alpha}\leq 1)\psi_j(t)
\eeq
exists in $L^2_x(\mathbb{R}^n)$ and for all $\alpha,\alpha'\in (0,1-2/n)$,
\eq
s\text{-}\lim\limits_{t\to \infty} e^{itH_0}\F_c(\frac{|x-2tP|}{t^\alpha}\leq 1)\psi_j(t)=s\text{-}\lim\limits_{t\to \infty} e^{itH_0}\F_c(\frac{|x-2tP|}{t^{\alpha'}}\leq 1)\psi_j(t);\label{weak: charge}
\eeq

\item there exist $N$ moving weakly localized parts, $\psi_{wl,j}(t)\equiv \psi_{wl,j,\epsilon}(t)$ such that the equation
\eq
\lim\limits_{t\to \infty}\| \psi_j(t)-e^{-itH_0}\psi_{j,+}(x)-\psi_{wl,j}(t)\|_{L^2_x(\mathbb{R}^n)}=0\label{AC: de: charge}
\eeq
holds true, and $\psi_{wl,j}(t), j=1,\cdots,N,$ are moving weakly localized parts around $tv_j$ satisfying
\eq
(e^{itP\cdot v_j}\psi_{wl,j}(t), |x|e^{itP\cdot v_j}\psi_{wl,j}(t))_{\s^2_x(\mathbb{R}^n)}\lesssim_\epsilon t^{1/2+\epsilon},\qquad t\geq 1.
\eeq
\end{enumerate}

\end{proposition}
\begin{proof} {\textbf{Proof of the existence of $\psi_{j,+}(x)$:}} We use Cook's method and the process similar to the proof of Theorem~\ref{thm1} to show the existence of $\psi_{j,+}(x)$ in $L^2_x(\mathbb{R}^n)$.

Let
\begin{align}
\psi_{j,+}(t) :=e^{itH_0}\F_c(\frac{|x-2tP|}{t^\alpha}\leq 1) \psi_j(t).
\end{align}
By Eq.~\eqref{id: Fc}, $\psi_{j,+}(t)$ reads
\eq
\psi_{j,+}(t) =F_c(\frac{|x|}{t^\alpha}\leq 1)e^{itH_0} \psi_{j}(t).
\eeq
Using Cook's method to expand $\psi_{j,+}(t)$ and, by Eqs.~\eqref{def: psij(t)} and 
\begin{align}
    \partial_t[e^{itH_0}\psi_j(t)]=& e^{itH_0}(iH_0-iH_0)\psi_j(t) +(-i)e^{itH_0}V_j(x-tv_j,t)\psi(t),
\end{align}
we obtain
\begin{align}
\psi_{j,+}(t)=&\psi_{j,+}(1)+(-i)\int_1^t  F_c(\frac{|x|}{s^\alpha}\leq 1)e^{isH_0}V_j(x-sv_j,s)\psi(s)ds\nonumber\\
&+\int_1^t  \partial_s[F_c(\frac{|x|}{s^\alpha}\leq 1)]e^{isH_0}\psi_j(s)ds\nonumber\\
=:&\psi_{j,+}(1)+\psi_{j,int}(t)+\psi_{j,p}(t).
\end{align}
By the unitarity of $e^{iH_0}$ and Lemma~\ref{sec: charge: psij}, $\psi_{j,+}(1)\in L^2_x(\mathbb{R}^n)$. By letting $\alpha\in (0, 1-2/n), n\geq 5$, and using Cauchy-Schwarz inequality and Assumption~\ref{asp: charge}, {similarly to~\eqref{eq: 3.26}} we find that $\psi_{j,int}(t)$ satisfies the estimate ($\beta=\frac{n(1-\alpha)}{2}-1$),
\begin{align}
    \| \psi_{j,int}(t)\|_{L^2_x(\mathbb{R}^n)}\leq & \int_1^t \| F_c(\frac{|x|}{s^\alpha}\leq 1)e^{isH_0}V_j(x-sv_j,s)\psi(s)\|_{L^2_x(\mathbb{R}^n)}ds\nonumber\\
    \lesssim_n & \int_1^t{\frac{1}{s^{1+\beta}} }\| V_j(x-sv_j,s)\psi(s)\|_{L^\infty_sL^1_x(\mathbb{R}^{n+1})}ds\nonumber\\
    \lesssim_n & \|\langle x\rangle^{n+1}V_j(x,s) \|_{L^\infty_{x,s}(\mathbb{R}^{n+1})}\|\psi(s)\|_{L^\infty_sL^2_x(\mathbb{R}^{n+1})},\label{eq: psi j int}
\end{align}
which leads to the existence of $\psi_{j,int}(\infty)$ in $\s^2_x(\mathbb{R}^n)$. For $\psi_{j,p}(t)$, we use {\bf{RPRES}} by taking $b=1$ and
\eq
\begin{cases}
B(t):=\F_c(\frac{|x|}{t^\alpha}\leq 1)\\
\phi(t)=e^{itH_0}\psi_j(t)
\end{cases}.
\eeq
We find that 
\begin{align}
    &\partial_t\langle B: \phi(t)\rangle_t\nonumber\\
    =& (\phi(t), \partial_t[\F_c(\frac{|x|}{t^\alpha}\leq 1)] \phi(t))_{L^2_x(\mathbb{R}^n)}+(-i)(\phi(t), \F_c(\frac{|x|}{t^\alpha}\leq 1)e^{itH_0}V_j(x-tv_j,t)\psi(t) )_{L^2_x(\mathbb{R}^n)}\nonumber\\
    &+i(\F_c(\frac{|x|}{t^\alpha}\leq 1)e^{itH_0}V_j(x-tv_j,t)\psi(t) , \phi(t))_{L^2_x(\mathbb{R}^n)}\nonumber\\
    =& (\phi(t), C^*C\phi(t))_{L^2_x(\mathbb{R}^n)}+g(t)
\end{align}
where $C^*C$ and $g(t)$ are given by, with $\F'_c(\lambda\leq 1)\equiv \frac{d}{d\lambda}[F_c(\lambda\leq 1)]$, 
\begin{align}
C^*C:=&\partial_t[\F_c(\frac{|x|}{t^\alpha}\leq 1)]\nonumber\\
=& \F_c'(\frac{|x|}{t^\alpha}\leq 1)\times \frac{-\alpha }{t}\times \frac{|x|}{t^\alpha}
\geq  0\label{Fcnneg: charge}
\end{align}
and
\begin{align}
    g(t):=& (-i)(\phi(t), \F_c(\frac{|x|}{t^\alpha}\leq 1)e^{itH_0}V_j(x-tv_j,t)\psi(t) )_{L^2_x(\mathbb{R}^n)}\nonumber\\
    &+i(\F_c(\frac{|x|}{t^\alpha}\leq 1)e^{itH_0}V_j(x-tv_j,t)\psi(t) , \phi(t))_{L^2_x(\mathbb{R}^n)},
\end{align}
respectively. We note that $\langle B: \phi(t)\rangle_t$ is uniformly bounded in $t$. This follows from the Cauchy-Schwarz inequality, the unitarity of $e^{itH_0}$ and Lemma~\ref{sec: charge: psij}:
\begin{align}
| \langle B: \phi(t)\rangle_t|= & ( e^{itH_0}\psi_j(t), F_c(\frac{|x|}{t^\alpha}\leq 1) e^{itH_0}\psi_j(t))_{L^2_x(\mathbb{R}^n)} \nonumber\\
\leq & \| e^{itH_0}\psi_j(t)\|^2_{L^2_x(\mathbb{R}^n)}
=  \| \psi_j(t)\|^2_{L^2_x(\mathbb{R}^n)}
\lesssim_E  1.\label{supBphi: charge}
\end{align}
Additionally, {similarly to~\eqref{eq: 3.26} and~\eqref{eq: psi j int},} $g(t)$ satisfies the following estimate due to Assumptions~\ref{asp: global} and~\ref{asp: charge}, and the Cauchy-Schwarz inequality:
\begin{align}
    |g(t)|\leq & 2\| \phi(t)\|_{L^2_x(\mathbb{R}^n)}\| \F_c(\frac{|x|}{t^\alpha}\leq 1)e^{itH_0}V_j(x-tv_j,t)\psi(t)\|_{L^2_x(\mathbb{R}^n)}\nonumber\\
    \lesssim_{n,E}& \frac{1}{t^{1+\beta}}\| V_j(x-tv_j,t)\psi(t)\|_{L^\infty_tL^1_x(\mathbb{R}^{n+1})}\lesssim_{n,E} \frac{1}{t^{1+\beta}}\| \langle x\rangle^n V_j(x,t)\|_{L^\infty_{x,t}(\mathbb{R}^{n+1})}\|\psi(t)\|_{L^\infty_tL^2_x(\mathbb{R}^{n+1})}.\label{ineq: g: charge}
\end{align}
This implies $g\in L^1_t[1,\infty)$. Hence, the family $\{B(t)\}_{t\in [1,\infty)}$ is a RPROB with respect to $\phi(t)=e^{itH_0}\psi_j(t)$ and by Eq.~\eqref{CC: ineq2}, Eqs.~\eqref{supBphi: charge} and~\eqref{ineq: g: charge}, we obtain
\begin{align}
   &\int_1^\infty  |(\phi(t), \partial_t[\F_c(\frac{|x|}{t^\alpha}\leq 1)]\phi(t))_{L^2_x(\mathbb{R}^n)}|dt\nonumber\\
   =& \int_1^\infty  (\phi(t), \partial_t[\F_c(\frac{|x|}{t^\alpha}\leq 1)]\phi(t))_{L^2_x(\mathbb{R}^n)}dt
   \leq   2\sup\limits_{t\in [1,\infty)}| \langle B: \phi(t)\rangle_t|+\| g(t)\|_{L^1_t[1,\infty)}\nonumber\\
   \lesssim_{n,E} & 1+\| \langle x\rangle^n V_j(x,t)\|_{L^\infty_{x,t}(\mathbb{R}^{n+1})}
    <\infty.\label{p: nonlocal: charge}
\end{align}
By Cauchy-Schwarz inequality and the non-negativity of $\partial_t[\F_c(\frac{|x|}{t^\alpha}\leq 1)]$ (see~\eqref{Fcnneg: charge}), $\psi_{j,p}(t)$ satisfies the estimate, for $T_2\geq T_1$,
\begin{align}
  &\| \psi_{j,p}(T_2)-\psi_{j,p}(T_1)\|_{L^2_x(\mathbb{R}^n)}\nonumber\\
  \leq & \| \int_{T_1}^{T_2} | \partial_t[\F_c(\frac{|x|}{t^\alpha}\leq 1)] \phi(t)|  dt \|_{L^2_x(\mathbb{R}^n)}\nonumber\\
  \leq & \| \left( \int_{T_1}^{T_2} \partial_t[\F_c(\frac{|x|}{t^\alpha}\leq 1)]dt\right)^{1/2}  \left( \int_{T_1}^{T_2} \partial_t[\F_c(\frac{|x|}{t^\alpha}\leq 1)]|\phi(t)|^2dt\right)^{1/2}\|_{L^2_x(\mathbb{R}^n)}\label{psipT2T1: charge}
\end{align}
By estimates~\eqref{dtFc} and~\eqref{p: nonlocal: charge}, estimate~\eqref{psipT2T1: charge} leads to 
\begin{align}
    \| \psi_{j,p}(T_2)-\psi_{j,p}(T_1)\|_{L^2_x(\mathbb{R}^n)}\leq & \|   \left( \int_{T_1}^{T_2} \partial_t[\F_c(\frac{|x|}{t^\alpha}\leq 1)]|\phi(t)|^2dt\right)^{1/2}\|_{L^2_x(\mathbb{R}^n)}\nonumber\\
    =& \left( \int_{T_1}^\infty  |(\phi(t), \partial_t[\F_c(\frac{|x|}{t^\alpha}\leq 1)]\phi(t))_{L^2_x(\mathbb{R}^n)}|dt\right)^{1/2} \to  0,
\end{align}
as $T_1\to \infty$. Hence, $\{ \psi_{j,p}(t)\}_{t\geq 1}$ is Cauchy in $L^2_x(\mathbb{R}^n)$ and therefore $\psi_{j,p}(\infty)$ exists in $L^2_x(\mathbb{R}^n)$. This together with the existence of $\psi_{j,int}(\infty)$ in $L^2_x(\mathbb{R}^n)$ and that $\psi_{j,+}(1)\in L^2_x(\mathbb{R}^n)$, implies the existence of $\psi_{j,+}(x)$ in $L^2_x(\mathbb{R}^n)$. {Here, we use the notation $\psi_{j,+}(x)$ similar to $\Omega_{\alpha}^*\psi(0)$ (defined in Eq.~\eqref{wave-4}) since $\psi_{j,+}(x)$ is independent on $\alpha\in (0,1-2/n)$. See the proof of Theorem~\ref{thm1} in Section~\ref{sec: proof of Thm1}.}

\textbf{Existence of $\psi_{wl,j}(t)$:} Given $\epsilon\in (0,1/2)$, take 
\eq
\psi_{wl,j}(t)\equiv \psi_{wl,j,\epsilon}(t)=\F_\alpha(\frac{|x-tv_j|}{(t+1)^{1/2+\epsilon}}<1)\psi_j(t), \quad j=1,\cdots,N.\label{weak: def: charge}
\eeq
Break $\psi_j(t)-\psi_{wl,j}(t)$ into three pieces:
\eq
\psi_j(t)-\psi_{wl,j}(t)=\psi_{j1}(t)+\psi_{j2}(t)+\psi_{j3}(t),\label{decomp: charge: weak}
\eeq
where $\psi_{jk}(t), k=1,2,3,$ are given by 
\eq
\psi_{j1}(t)=\F_c(\frac{|x-tv_j|}{(t+1)^{1/2+\epsilon}}\geq 1)P^+_{tv_j}\F_1(\sqrt{t+1}|2P-v_j|\geq 1)\psi_j(t),
\eeq
\eq
\psi_{j2}(t)=\F_c(\frac{|x-tv_j|}{(t+1)^{1/2+\epsilon}}\geq 1)P^-_{tv_j}\psi_j(t)
\eeq
and
\eq
\psi_{j3}(t)=\F_c(\frac{|x-tv_j|}{(t+1)^{1/2+\epsilon}}\geq 1)P^+_{tv_j}\F_1(\sqrt{t+1}|2P-v_j|< 1)\psi_j(t),
\eeq
respectively. We approximate $\psi_{j1}(t)$ by $e^{-itH_0}\psi_{j,+}(x)$ and arrive at~\eqref{AC: de: charge} by showing that
\eq
\lim\limits_{t\to \infty}\| \psi_{jl}(t)\|_{L^2_x(\mathbb{R}^n)}=0,\label{lim: psi2: charge}\quad l=2,3,
\eeq
\eq
\lim\limits_{t\to \infty}\| \psi_{j1}(t)-\F_c( \frac{|x-tv_j|}{(t+1)^{1/2+\epsilon}}\geq 1)P^+\F_1(\sqrt{t+1}|2P-v_j|\geq 1)e^{-itH_0} \psi_{j,+}(x)\|_{L^2_x(\mathbb{R}^n)}=0\label{lim:psi3: charge}
\eeq
and
\begin{align}
\lim\limits_{t\to \infty}&\| \F_c( \frac{|x-tv_j|}{(t+1)^{1/2+\epsilon}}\geq 1)P^+\F_1(\sqrt{t+1}|2P-v_j|\geq 1)e^{-itH_0} \psi_{j,+}(x) \nonumber\\
&\qquad\qquad\qquad\qquad\qquad-e^{-itH_0} \psi_{j,+}(x)\|_{L^2_x(\mathbb{R}^n)}=0.\label{lim:psi4: charge}
\end{align}
By Lemma~\ref{Lem: Pprop: charge} and Assumptions~\ref{asp: global} and~\ref{asp: charge}, $\psi_{j2}(t)$ satisfies, with $\delta>2$, 
\begin{align}
    \| \psi_{j2}(t)\|_{L^2_x(\mathbb{R}^n)}\lesssim & \int_0^t \| \F_c (\frac{|x-tv_j|}{(t+1)^{1/2+\epsilon}}\geq 1) P^-_{tv_j}e^{-i(t-s)H_0}\F_1(\sqrt{t+1}|2P-v_j|\geq 1)\nonumber\\
    &\times\langle x-sv_j\rangle^{-\delta}\|_{L^2_x(\mathbb{R}^n)\to L^2_x(\mathbb{R}^n)} \| \langle x-sv_j\rangle^\delta V_j(x-sv_j,s)\|_{L^\infty_{x}(\mathbb{R}^n)} \|\psi(s)\|_{L^2_x(\mathbb{R}^n)}ds\nonumber\\
    \lesssim_{E,\epsilon} & \int_0^t \frac{1}{\langle  (t+1)^{1/2+\epsilon}+s/\sqrt{t+1}\rangle^{\delta}} \sup\limits_{u\in \mathbb{R}} \| \langle x\rangle^\delta V_j(x,u)\|_{L^\infty_x(\mathbb{R}^n)}ds \nonumber\\
    \lesssim_{E,\epsilon, \delta} & \frac{t}{\langle t+1\rangle^{1+2\epsilon}} \sup\limits_{u\in \mathbb{R}} \| \langle x\rangle^\delta V_j(x,u)\|_{L^\infty_x(\mathbb{R}^n)}\to  0,\label{psi22: charge}
\end{align}
as $t\to \infty.$ Eq.~\eqref{psi22: charge} implies 
\eq
\lim\limits_{t\to \infty} \| \psi_2(t)\|_{L^2_x(\mathbb{R}^n)}=0.\label{lim psi20: charge}
\eeq
By Lemma~\ref{lem: Ppropfree: charge} and Assumptions~\ref{asp: global} and~\ref{asp: charge}, $\psi_{j3}(t)$ satisfies, with $\delta>2$, 
\begin{align}
   & \| \psi_{j3}(t)\|_{L^2_x(\mathbb{R}^n)}\leq  \int_0^t \| \F_c (\frac{|x-tv_j|}{(t+1)^{1/2+\epsilon}}\geq 1) P_{tv_j}^+e^{-i(t-s)H_0}\F_1(\sqrt{t+1}|2P-v_j|<1)\nonumber\\
    &\times\langle x-sv_j\rangle^{-\delta}\|_{L^2_x(\mathbb{R}^n)\to L^2_x(\mathbb{R}^n)}\| \langle x-sv_j\rangle^\delta V_j(x-sv_j,s)\|_{L^\infty_x(\mathbb{R}^n)}\|\psi(s)\|_{L^2_x(\mathbb{R}^n)}ds\nonumber\\
    \lesssim_{E,\epsilon} & \int_0^t \frac{1}{ (t+1)^{\frac{1}{2}\delta+\epsilon\delta}} \sup\limits_{u\in \mathbb{R}} \| \langle x\rangle^\delta V_j(x,u)\|_{L^\infty_x(\mathbb{R}^n)}ds  \lesssim_{E,\epsilon,\delta}  \frac{t}{\langle t+1\rangle^{1+2\epsilon}} \sup\limits_{u\in \mathbb{R}}\| \langle x\rangle^\delta V_j(x,u)\|_{L^\infty_x(\mathbb{R}^n)}
    \to  0\label{psi32: charge}
\end{align}
as $t\to \infty.$ Eq.~\eqref{psi32: charge} implies 
\eq
\lim\limits_{t\to \infty} \| \psi_3(t)\|_{L^2_x(\mathbb{R}^n)}=0.\label{lim psi3: charge}
\eeq
Eqs.~\eqref{lim psi20: charge} and ~\eqref{lim psi3: charge} imply Eq.~\eqref{lim: psi2: charge}.

Next, we prove~\eqref{lim:psi3: charge}. Let 
\begin{align}
\tilde\psi_{j,+}(x):=w\text{-}\lim\limits_{t\to \infty}(-i)\int_0^t e^{isH_0}V_j(x-sv_j,s)\psi(s)ds.
\end{align}
By the existence of $\psi_{j,+}(x)$ and Eq.~\eqref{weak: charge}, $\tilde\psi_{j,+}(x)$ exists in $L^2_x(\mathbb{R}^n)$ and $\tilde\psi_{j,+}(x)=\psi_{j,+}(x)$. By Duhamel's expansion and $\tilde \psi_{j,+}=\psi_{j,+}$, we have 
\begin{align}
    & \F_c( \frac{|x-tv_j|}{(t+1)^{1/2+\epsilon}}\geq 1)P_{tv_j}^+\F_1(\sqrt{t+1}|2P-v_j|<1)e^{-itH_0} \psi_{j,+}- \psi_{j1}(t)\nonumber\\
    =&(-i)\F_c(\frac{|x-tv_j|}{(t+1)^{1/2+\epsilon}}\geq 1) P^+\int_t^\infty \F_1(\sqrt{t+1}|2P-v_j|<1) e^{i(s-t)H_0}V_j(x-sv_j,s)\psi(s)ds.\label{id: psi1: charge}
\end{align}
By Lemma~\ref{Lem: Pprop: charge} and Assumptions~\ref{asp: global} and~\ref{asp: charge}, \eqref{id: psi1: charge} implies, with $\delta>2$, 
\begin{align}
    & \|\F_c( \frac{|x-tv_j|}{(t+1)^{1/2+\epsilon}}\geq 1)P^+ \F_1(\sqrt{t+1}|2P-v_j|<1)e^{-itH_0} \psi_{j,+}- \psi_1(t)\|_{L^2_x(\mathbb{R}^n)}\nonumber\\
    \leq &\int_t^\infty\|\F_c(\frac{|x-tv_j|}{(t+1)^{1/2+\epsilon}}\geq 1) P_{tv_j}^+ \F_1(\sqrt{t+1}|2P-v_j|<1)e^{i(s-t)H_0}\langle x-sv_j\rangle^{-\delta}\|_{L^2_x(\mathbb{R}^n)\to L^2_x(\mathbb{R}^n)}\nonumber\\
    &\times \|\langle x-sv_j\rangle^\delta V_j(x-sv_j,s)\|_{L^\infty_x(\mathbb{R}^n)}\|\psi(s)\|_{L^2_x(\mathbb{R}^n)}ds\nonumber\\
    \lesssim_{E,\epsilon}& \int_t^\infty \frac{1}{\langle (t+1)^{1/2+\epsilon}+(s-t)/\sqrt{t+1}\rangle^\delta}\sup\limits_{u\in \mathbb{R}}\|\langle x\rangle^\delta V_j(x,u)\|_{L^\infty_x(\mathbb{R}^n)}ds\nonumber\\
    \lesssim_{E,\epsilon} & \frac{\sqrt{t+1}}{\langle t+1\rangle^{\frac{1}{2}(\delta-1)+\epsilon(\delta -1)}}\sup\limits_{u\in \mathbb{R}}\|\langle x\rangle^\delta V_j(x,u)\|_{L^\infty_x(\mathbb{R}^n)}
    \to 0\label{lim psi1(t): local: charge}
\end{align}
as $t\to \infty.$ {Estimate \eqref{lim psi1(t): local: charge}} implies ~\eqref{lim:psi3: charge}. 

Now we prove~\eqref{lim:psi4: charge}. Equation 
\eq
s\text{-}\lim\limits_{t\to \infty} \F_1(\sqrt{t+1}|2P-\eta|< 1)=0\quad \text{ on }L^2_x(\mathbb{R}^n),
\eeq
together with Lemmas~\ref{lem: Ppmf: charge} and~\ref{lem: not linear: charge}, implies
\begin{align}
    & \|e^{-itH_0}\psi_{j,+}(x) -\F_c( \frac{|x-tv_j|}{(t+1)^{1/2+\epsilon}}\geq 1)P_{tv_j}^+\F_1(\sqrt{t+1}|2P-v_j|\geq 1)e^{-itH_0} \psi_{j,+}(x)\|_{L^2_x(\mathbb{R}^n)}\nonumber\\
    \leq & \| \F_c( \frac{|x-tv_j|}{(t+1)^{1/2+\epsilon}}\geq 1)P_{tv_j}^+\F_1(\sqrt{t+1}|2P-\eta|< 1)e^{-itH_0} \psi_{j,+}(x) \|_{L^2_x(\mathbb{R}^n)}\nonumber\\
    & +\| \F_c( \frac{|x-tv_j|}{(t+1)^{1/2+\epsilon}}< 1)P_{tv_j}^+e^{-itH_0} \psi_{j,+}(x)\|_{L^2_x(\mathbb{R}^n)}+\| P_{tv_j}^-e^{-itH_0} \psi_{j,+}(x) \|_{L^2_x(\mathbb{R}^n)}\nonumber\\
    \leq & \|\F_1(\sqrt{t+1}|2P-\eta|< 1) \psi_{j,+}(x) \|_{L^2_x(\mathbb{R}^n)}+\| \F_c( \frac{|x-tv_j|}{(t+1)^{1/2+\epsilon}}< 1)P_{tv_j}^+e^{-itH_0} \psi_{j,+}(x) \|_{L^2_x(\mathbb{R}^n)}\nonumber\\
    & +\| P_{tv_j}^-e^{-itH_0} \psi_{j,+} \|_{L^2_x(\mathbb{R}^n)} \to  0\label{psi(3): local: charge}
\end{align}
as $t\to \infty.$ {Estimate \eqref{psi(3): local: charge}} implies~\eqref{lim:psi4: charge}. By~\eqref{lim: psi2: charge},~\eqref{lim:psi3: charge},~\eqref{lim:psi4: charge} and Eq.~\eqref{decomp: charge: weak}, together with Eq.~\eqref{weak: def: charge}, we arrive at
\begin{align}
& \| \psi_j(t)-  \psi_{wl,j}(t)-e^{-itH_0}\psi_{j,+} \|_{L^2_x(\mathbb{R}^n)}\nonumber\\
   = & \| \psi_j(t)-  \F_c( \frac{|x-tv_j|}{(t+1)^{1/2+\epsilon}}< 1)\psi_j(t)-e^{-itH_0}\psi_{j,+}(x) \|_{L^2_x(\mathbb{R}^n)}\nonumber\\
    \leq & \|  \F_c( \frac{|x-tv_j|}{(t+1)^{1/2+\epsilon}}\geq 1)\psi_j(t)-e^{-itH_0}\psi_{j,+}(x) \|_{L^2_x(\mathbb{R}^n)}\nonumber\\
    \leq &\| \psi_{j1}(t)-\F_c( \frac{|x-tv_j|}{(t+1)^{1/2+\epsilon}}\geq 1)P^+\F_1(\sqrt{t+1}|2P-v_j|\geq 1)e^{-itH_0} \psi_{j,+}\|_{L^2_x(\mathbb{R}^n)}\nonumber\\
    &+\| e^{-itH_0}\psi_{j,+} - \F_c( \frac{|x-tv_j|}{(t+1)^{1/2+\epsilon}}\geq 1)P^+\F_1(\sqrt{t+1}|2P-v_j|\geq 1)e^{-itH_0} \psi_{j,+}\|_{L^2_x(\mathbb{R}^n)}\nonumber\\
    &+\| \psi_{j2}(t)\|_{L^2_x(\mathbb{R}^n)}+\| \psi_{j3}(t)\|_{L^2_x(\mathbb{R}^n)}
    \to  0
\end{align}
as $t\to \infty.$

\end{proof}

\begin{proof}[Proof of Theorem~\ref{thm3}] By Duhamel's formula, $\psi(t)$ reads
\begin{equation}
    \psi(t)=e^{-itH_0}\psi(0)+\sum\limits_{j=1}^N \psi_j(t),
\end{equation}
which, together with Theorem~\ref{thm1}, Eq.~\eqref{def: psij(t)} and Eq.~\eqref{def: psijfree} of Proposition~\ref{prop: charge}, implies 
\begin{align}\label{charge: free: id}
\Omega_{\alpha}^*\psi(0)=&s\text{-}\lim\limits_{t\to \infty} e^{itH_0}\F_c(\frac{|x-2tP|}{t^\alpha}\leq 1) \psi(t)\nonumber\\
=& \lim\limits_{t\to \infty} e^{itH_0}\F_c(\frac{|x-2tP|}{t^\alpha}\leq 1) e^{-itH_0} \psi(0)+\sum\limits_{j=1}^N \lim\limits_{t\to \infty}e^{itH_0}\F_c(\frac{|x-2tP|}{t^\alpha}\leq 1) \psi_j(t)\nonumber\\
=&\lim\limits_{t\to \infty} \F_c(\frac{|x|}{t^\alpha}\leq 1)  \psi(0)+\sum\limits_{j=1}^N \lim\limits_{t\to \infty}e^{itH_0}\F_c(\frac{|x-2tP|}{t^\alpha}\leq 1) \psi_j(t) \nonumber\\
=& \psi(0)+\sum\limits_{j=1}^N \psi_{j,+}(x).
\end{align}
By Proposition~\ref{prop: charge} and Eq.~\eqref{charge: free: id}, we have
\begin{align}
    & \limsup\limits_{t\to \infty}\| \psi(t)-e^{-itH_0}\Omega_{\alpha}^*\psi(0)-\sum\limits_{j=1}^N \psi_{wl,j}(t)\|_{L^2_x(\mathbb{R}^n)}\nonumber\\
    =&  \limsup\limits_{t\to \infty}\| e^{-itH_0}\psi(0)+\sum\limits_{j=1}^N\psi_j(t)-e^{-itH_0}\Omega_{\alpha}^*\psi(0)-\sum\limits_{j=1}^N \psi_{wl,j}(t)\|_{L^2_x(\mathbb{R}^n)}\nonumber\\
    \leq & \sum\limits_{j=1}^N\limsup\limits_{t\to \infty}\| \psi_j(t)-e^{-itH_0}\psi_{j,+}(x)-\psi_{wl,j}(t)\|_{L^2_x(\mathbb{R}^n)}
    =       0.
\end{align}

\end{proof}

\section{Proof of Proposition~\ref{NLSapplication}}\label{app}
\begin{proof}[Proof for Proposition \ref{NLSapplication}] We note that, due to Assumption~\ref{asp: global} with $a=1$, 
\begin{equation}
    \|\F_c(\frac{|x|}{t^\alpha}\leq 1) e^{itH_0}\psi(t)\|_{\mathcal{H}^1_x}\lesssim_E 1,\qquad \forall t\geq 1.
\end{equation}
This, together with the existence of $\Omega_\alpha^*\psi(0)$ in $L^2_x(\mathbb{R}^n)$, implies 
\begin{equation}
    \Omega_\alpha^*\psi(0)\in \mathcal{H}^1_x. 
\end{equation}
Next, we prove that $\psi(t)$ satisfies the endpoint Strichartz estimate 
\eq
\| \psi(t)\|_{L^2_tL^{6}_x(\mathbb{R}^{3+1})}\lesssim_E 1\label{Strichartz global}
\eeq
provided that
\eq
\limsup\limits_{t\to \infty}\| \psi(t)-e^{-itH_0}\Omega_\alpha^*\psi(0)\|_{\mathcal{H}^1_x} \text{ is sufficiently small.}
\eeq
We note that by Assumption~\ref{asp: global},
\begin{equation}
    \| \psi(t)\|_{L^\infty_tL^6_x(\mathbb{R}^{n+1})}\lesssim_E 1,
\end{equation}
which implies the endpoint Strichartz estimate locally in $t$:
\begin{equation}
    \| \psi(t)\|_{L^2_tL^6_x(\mathbb{R}^3\times [0,T])}\lesssim_{E,T}1,\qquad T>0.\label{local Stri}
\end{equation}
{Next, we} use estimate~\eqref{local Stri} to prove estimate~\eqref{Strichartz global}. The endpoint Strichartz estimate of $\psi(t)$ follows from a standard contraction argument. By Duhamel's formula, $\psi(t)-e^{-itH_0}\Omega_\alpha^*\psi(0)$ reads 
\begin{align}
    \psi(t)-e^{-itH_0}\Omega_\alpha^*\psi(0)= & e^{-itH_0}\psi(0)+(-i)\int_0^t e^{-i(t-s)H_0} I(\psi(s))\psi(s) ds -e^{-itH_0}\Omega_\alpha^*\psi(0)\nonumber\\
    =& e^{-itH_0}(\psi(0)-\Omega_\alpha^*\psi(0))+(-i)\int_0^t e^{-i(t-s)H_0} I(\psi(s))\psi(s) ds.\label{Duhamel}
\end{align}
By writing 
\begin{align}
    &(-i)\int_0^t e^{-i(t-s)H_0} I(|\psi(s)|)\psi(s) ds\nonumber\\
    =&(-i)\int_0^t e^{-i(t-s)H_0} (I(|\psi(s)|)\psi(s)-I(|e^{-isH_0}\Omega_{\alpha}^*\psi(0)|)e^{-isH_0}\Omega_{\alpha}^*\psi(0)) ds\nonumber\\
    &+(-i)\int_0^t e^{-i(t-s)H_0} I(|e^{-isH_0}\Omega_{\alpha}^*\psi(0)|)e^{-isH_0}\Omega_{\alpha}^*\psi(0) ds \nonumber\\
    =:& \psi_{r}(t)+\psi_{f}(t),
\end{align}
according to Eq.~\eqref{Duhamel}, we arrive at 
\begin{align}
    \psi(t)-e^{-itH_0}\Omega_\alpha^*\psi(0)=& e^{-itH_0}(\psi(0)-\Omega_\alpha^*\psi(0))+\psi_{r}(t)+\psi_{f}(t).
\end{align}
By Strichartz estimate of free flows, we obtain 
\begin{align}
\| e^{-itH_0}(\psi(0)-\Omega_\alpha^*\psi(0))\|_{L^2_tL^6_x(\mathbb{R}^{3+1})}\lesssim &\|\psi(0)-\Omega_\alpha^*\psi(0))\|_{L^2_x(\mathbb{R}^3)}\nonumber\\
\lesssim & \|\psi(0)\|_{L^2_x(\mathbb{R}^3)}. \label{Apr.16.1}
\end{align}
By the dual homogeneous Strichartz estimates 
\eq
\|\int_{t>s} e^{isH_0}F(s)ds \|_{L^2_tL^6_x(\mathbb{R}^{3+1})}\lesssim \| F(s)\|_{L^2_sL^{6/5}_x(\mathbb{R}^{3+1})},
\eeq
H\"older's inequality and Eq.~\eqref{con: I} ({with $g=0$}), $\psi_{f}(t)$ satisfies 
\begin{align}
 \| \psi_f(t)\|_{L^2_tL^6_x(\mathbb{R}^{3+1})}=& \| \int_0^t e^{-i(t-s)H_0} I(e^{-isH_0}\Omega_\alpha^*\psi(0))e^{-isH_0}\Omega_\alpha^*\psi(0))ds\|_{L^2_tL^6_x(\mathbb{R}^{3+1}) }\nonumber\\
 \lesssim & \| I(e^{-isH_0}\Omega_\alpha^*\psi(0))e^{-isH_0}\Omega_\alpha^*\psi(0) \|_{L^2_sL^{6/5}_x(\mathbb{R}^{3+1})}\nonumber\\
 \lesssim_{E} & \|e^{-isH_0}\Omega_\alpha^*\psi(0) \|_{L^2_sL^{6}_x(\mathbb{R}^{3+1})}\nonumber\\
 \lesssim_E & \| \psi(0) \|_{L^2_x(\mathbb{R}^3)}. \label{Apr.16.2}
\end{align}
By inequality~\eqref{con: I}, estimate~\eqref{local Stri} and the inhomogeneous Strichartz estimate 
\eq
\|\int_{s\leq t} e^{-i(t-s)H_0} F(s)ds\|_{L^2_tL^6_x(\mathbb{R}^{3+1})}\lesssim \| F(s)\|_{L^2_tL^{6/5}_x(\mathbb{R}^{3+1})},
\eeq
$\psi_r(t)$ satisfies, for $T\geq T_1\geq 0$, 
\begin{align}
   & \| \psi_r(t)\|_{L^2_tL^6_x(\mathbb{R}^{3}\times[0,T])}\nonumber\\
   \leq& \| \int_0^t e^{-i(t-s)H_0} (I(|\psi(s)|)\psi(s)-I(|e^{-isH_0}\Omega_{\alpha}^*\psi(0)|)e^{-isH_0}\Omega_{\alpha}^*\psi(0)) ds\|_{L^2_tL^6_x(\mathbb{R}^{3}\times[0,T])}\nonumber\\
    \lesssim & \|I(|\psi(s)|)\psi(s)-I(|e^{-isH_0}\Omega_{\alpha}^*\psi(0)|)e^{-isH_0}\Omega_{\alpha}^*\psi(0) \|_{L^2_sL^{6/5}_x(\mathbb{R}^{3}\times [0,T])}\nonumber\\
    \lesssim & C_{I1}\| \psi(s)-e^{-isH_0}\Omega_\alpha^*\psi(0)\|_{L^\infty_s\mathcal H^1_x(\mathbb{R}^{3}\times [T_1,T])} \| \psi(s)-e^{-isH_0}\Omega_\alpha^* \psi(0)\|_{L^2_sL^6_x(\mathbb{R}^{3}\times[0,T])}\nonumber\\
    &+ C_{I1}\| \psi(s)-e^{-isH_0}\Omega_\alpha^*\psi(0)\|_{L^\infty_s\mathcal H^1_x(\mathbb{R}^{3}\times [0,T_1])} \| \psi(s)-e^{-isH_0}\Omega_\alpha^* \psi(0)\|_{L^2_sL^6_x(\mathbb{R}^{3}\times [0,T_1])}\nonumber\\
    & + C_{I2} \| e^{-isH_0}\Omega_\alpha^* \psi(0)\|_{L^2_sL^6_x(\mathbb{R}^{3+1})}\nonumber\\
    \leq & C_1\| \psi(s)-e^{-isH_0}\Omega_\alpha^*\psi(0)\|_{L^\infty_s\mathcal H^1_x(\mathbb{R}^{3}\times [T_1,T])}\| \psi(s)-e^{-isH_0}\Omega_\alpha^* \psi(0)\|_{L^2_sL^6_x(\mathbb{R}^{3}\times[0,T])}+C_2\|\psi(0)\|_{L^2_x(\mathbb{R}^3)}\label{Apr.16.3}
\end{align}
where $C_1=C_1(E)>0$ and $C_2=C_2(E,T_1)>0$ denote two positive constants. Estimates~\eqref{Apr.16.1},~\eqref{Apr.16.2} and~\eqref{Apr.16.3} imply, for all $T\geq T_1$,
\begin{align}
    &\| \psi(t)-e^{-itH_0}\Omega_\alpha^*\psi(0)\|_{L^2_tL^6_x(\mathbb{R}^3\times [0,T])}\nonumber\\
    \leq& C_1\| \psi(s)-e^{-isH_0}\Omega_\alpha^*\psi(0)\|_{L^\infty_s\mathcal H^1_x(\mathbb{R}^{3}\times [T_1,T])}\| \psi(s)-e^{-isH_0}\Omega_\alpha^* \psi(0)\|_{L^2_sL^6_x(\mathbb{R}^{3}\times[0,T])}\nonumber\\
    &+C_2\|\psi(0)\|_{L^2_x(\mathbb{R}^3)},
\end{align}
where $C_1=C(E)>0$ and $C_2=C(E,T_1)>0$ are two positive constants. By taking 
\eq
m{\equiv m(E):=\frac{1}{C_1(E)}}
\eeq
and by taking $T_1$ large enough such that $t\geq T_1$ implies 
\begin{align}
    &\| \psi(t)-e^{-itH_0}\Omega_\alpha^*\psi(0)\|_{\mathcal H^1_x(\mathbb{R}^3)}<m,
\end{align}
we have, with
\eq
\tilde C:= C_1 \| \psi(t)-e^{-itH_0}\Omega_\alpha^*\psi(0)\|_{L^\infty_t\mathcal H^1_x(\mathbb{R}^3\times[T_1,\infty))}<1,
\eeq
for all $T\geq T_1$,
\begin{align}
    \| \psi(t)-e^{-itH_0}\Omega_\alpha^*\psi(0)\|_{L^2_tL^6_x(\mathbb{R}^3\times [0,T])}\leq C_2 \|\psi(0)\|_{L^2_x(\mathbb{R}^3)}+\tilde C \| \psi(t)-e^{-itH_0}\Omega_\alpha^*\psi(0)\|_{L^2_tL^6_x(\mathbb{R}^3\times [0,T])},
\end{align}
which leads to 
\begin{align}
     \| \psi(t)-e^{-itH_0}\Omega_\alpha^*\psi(0)\|_{L^2_tL^6_x(\mathbb{R}^3\times [0,T])}\lesssim_{E} \|\psi(0)\|_{L^2_x(\mathbb{R}^3)},\quad \forall t\geq T_1.
\end{align}
By taking $T\to\infty$, we arrive at 
\begin{align}
     \| \psi(t)-e^{-itH_0}\Omega_\alpha^*\psi(0)\|_{L^2_tL^6_x(\mathbb{R}^{3+1})}\lesssim_{E} \|\psi(0)\|_{L^2_x(\mathbb{R}^3)}
\end{align}
and therefore, 
\begin{align}
     \| \psi(t)\|_{L^2_tL^6_x(\mathbb{R}^{3+1})}\lesssim_{E} \|\psi(0)\|_{L^2_x(\mathbb{R}^3)}.
\end{align}
By using Duhamel's formula, the dual homogeneous Strichartz estimate
\eq
\|\int_{\mathbb{R}} e^{isH_0}F(s)ds \|_{L^2_x(\mathbb{R}^3)}\lesssim \| F(s)\|_{L^2_sL^{6/5}_x(\mathbb{R}^{3+1})},
\eeq
H\"older's inequality, Eq.~\eqref{con: I} (with $f=\psi$ and $g=0$) and Strichartz estimate of $\psi(t)$, we obtain
\begin{align}
    \|\psi(t)-e^{-itH_0}\Omega_\alpha^*\psi(0)\|_{\mathcal{H}^1_x}= & \| \int_t^\infty e^{-i(t-s)H_0} I(|\psi(s)|)\psi(s)ds\|_{L^2_x(\mathbb{R}^3)}\nonumber\\
    \lesssim& \|  \chi(|s|\geq t)I(|\psi(s)|)\psi(s) \|_{L^2_sL^{6/5}_x(\mathbb{R}^3)}\nonumber\\
    \lesssim_E & \|\chi(s\geq t) \psi(s)\|_{L^2_sL^6_x(\mathbb{R}^{3+1})}\nonumber\\
    \to & 0
\end{align}
as $t\to \infty$. We finish the proof.

\end{proof}

\section*{Acknowledgements}
The authors thank Scipio Cuccagna, Ting Zhang and Jiayan Wu for carefully reading the manuscript and for providing helpful and valuable comments. They also thank the anonymous referees for their many useful suggestions. X.W. acknowledges the support from Australian Laureate Fellowships, grant FL220100072. A.S. acknowledges the support from National Science Foundation via the grant DMS-2205931. Part of this work was completed by the second author while at Rutgers University, Texas A\&M University and the Fields Institute.

\appendix
\section{Phase-space operators}\label{app: phase-space}
\begin{proof}[Proof of Eq.~\eqref{id: Fc}] Let $P=(P_1,\cdots,P_n)$ and $x=(x_1,\cdots,x_n)$. To compute $i[H_0,x_j]$, $j=1,\cdots,n$, we find 
\begin{align}
i[H_0,x_j]=& i (H_0x_j)+2i (P_jx_j)P_j\nonumber\\
=& 2P_j,
\end{align}
which implies 
\eq
i[H_0,x]=2P
\eeq
and therefore, 
\begin{align}
\partial_t[e^{-itH_0} xe^{itH_0}]=&e^{-itH_0}(-i) [H_0,x]e^{itH_0}\nonumber\\
=& -2P.\label{free: H_0}
\end{align}
Eq.~\eqref{free: H_0} implies 
\begin{align}
e^{-itH_0} xe^{itH_0}=&e^{-isH_0} xe^{isH_0}\vert_{s=0}+\int_0^t  \partial_s[ e^{-isH_0} xe^{isH_0}]ds\nonumber\\
=& x-2tP,
\end{align}
which leads to
\eq
e^{-itH_0} e^{ix\cdot \xi}e^{itH_0}=e^{-i(x-2tP)\cdot\xi},\qquad \xi\in \mathbb{R}^n.\label{free: apen: eq2}
\eeq
Therefore, by Fourier transform and Eq.~\eqref{free: apen: eq2}, $e^{-itH_0}\F_c(\frac{|x|}{t^\alpha}> 1)e^{itH_0}$ reads
\begin{align}
    e^{-itH_0}\F_c(\frac{|x|}{t^\alpha}> 1)e^{itH_0}= & \frac{1}{(2\pi)^{n/2}}\int_{\mathbb{R}^n} \hat{\F}_c(\xi)e^{-itH_0} e^{i\xi\cdot \frac{x}{t^\alpha}} e^{itH_0}d^n \xi\nonumber\\
    =& \frac{1}{(2\pi)^{n/2}}\int_{\mathbb{R}^n} \hat{\F}_c(\xi) e^{i\xi\cdot \frac{x-2tP}{t^\alpha}} d^n \xi\nonumber\\
    =& \F_c(\frac{|x-2tP|}{t^\alpha}> 1),
\end{align}
which implies Eq.~\eqref{id: Fc}. This completes the proof.

\end{proof}
\section{Estimations of Free Flows and Differential Operators}\label{app: free flows}
\begin{proof}[Proof of estimate~\eqref{local: est}] Let 
\eq
\mathcal O(t):= \F_c(\frac{|x|}{t^\alpha}\leq 1)\F_1(t^\beta|P|>1)e^{\pm itH_0}\langle x\rangle^{- \sigma}.
\eeq
Break $\mathcal O(t)$ into two pieces:
\begin{align}
    \mathcal O(t)=\mathcal O_1(t)+\mathcal O_2(t),
\end{align}
where 
\begin{align}
    \mathcal O_1(t):=\F_c(\frac{|x|}{t^\alpha}\leq 1)\F_1(t^\beta|P|>1)e^{\pm itH_0}\langle x\rangle^{- \sigma}\chi(|x|\geq \frac{1}{10}t^{1-\beta})
\end{align}
and
\begin{align}
    \mathcal O_2(t):=\F_c(\frac{|x|}{t^\alpha}\leq 1)\F_1(t^\beta|P|>1)e^{\pm itH_0}\langle x\rangle^{- \sigma}\chi(|x|< \frac{1}{10}t^{1-\beta}).
\end{align}
By using the weight $\langle x\rangle^{-\sigma }$ and the unitarity of $e^{\pm i tH_0}$, we obtain, for $t\geq 1$, 
\begin{align}
   & \| \mathcal O_1(t)\|_{L^2_x(\mathbb{R}^n)\to L^2_x(\mathbb{R}^n)}\nonumber\\
   \leq & \| \F_c(\frac{|x|}{t^\alpha}\leq 1)\F_1(t^\beta|P|>1)e^{\pm itH_0}\|_{L^2_x(\mathbb{R}^n)\to L^2_x(\mathbb{R}^n)}\| \langle x\rangle^{-\sigma} \chi(|x|\geq \frac{1}{10}t^{1-\beta})\|_{L^2_x(\mathbb{R}^n)\to L^2_x(\mathbb{R}^n)}\nonumber\\
   \leq & \| \langle x\rangle^{-\sigma} \chi(|x|\geq \frac{1}{10}t^{1-\beta})\|_{L^2_x(\mathbb{R}^n)\to L^2_x(\mathbb{R}^n)}\lesssim  \frac{1}{t^{\sigma (1-\beta)}}.\label{Oeq1}
\end{align}
For $\mathcal O_2(t)$, take $f\in L^2_x(\mathbb{R}^n)$. By Fourier transform, $\mathcal O_2(t)f$ reads, with $\F_c\equiv \F_c(\frac{|x|}{t^\alpha}\leq 1)$ and $\F_1\equiv \F_1(t^\beta|q|>1)$, 
\begin{equation}
    \mathcal O_2f=\frac{1}{(2\pi)^n} \int \F_c e^{\pm itq^2} \F_1(t^\beta |q|>1) e^{-iy\cdot q}\langle y\rangle^{-\sigma} f(y)dy dq,\label{def: O2f}
\end{equation}
where we have used 
\begin{equation}
   [\langle x\rangle^{-\sigma} f](\hat{q})=\frac{1}{(2\pi)^{n/2}}\int e^{-iq\cdot y} \langle y\rangle^{-\sigma}f(y) dy.
\end{equation}
When $q$ is in the support of $\F_1(t^\beta |q|>1)$, $x$ is in the support of $\F_c(\frac{|x|}{t^\alpha}\leq 1)$ and $|y|\leq \frac{1}{10}t^{1-\beta}$, $x\pm 2tq-y$ satisfies, for $t\gg 1$ and $\alpha \in (0,1-\beta)$ 
\begin{align}
|x\pm 2tq-y|\geq & |2tq|-|x|-|y|\geq\frac{9}{10} t^{1-\beta}-2t^\alpha \geq  \frac{1}{2} t^{1-\beta}\label{t: lem1}
\end{align}
and
\begin{equation}
    |x\pm 2tq-y|\geq|2tq|-|x|-|y|\geq |tq|.\label{tq: lem1}
\end{equation}
We define an orthogonal basis $\{e_1,\cdots,e_n\}$ in $\mathbb{R}^n$ with $e_1$ satisfying  
\eq
|x_1\pm 2tq_1-y_1|\geq C|x\pm 2tq-y|,
\eeq
for some constant $C=C(n)>0$, where $x_1:=x\cdot e_1$, $q_1:=q\cdot e_1$ and $y_1:=y\cdot e_1$. Let $\tilde \F_1(k\sim 1)$ denote a smooth cut-off function satisfying, with $\F_1'(k>1)\equiv \frac{d}{dk}[\F_1(k>1)]$ and $[a,b]:=\text{supp}(\F_1'(k>1))$, 
\begin{equation}
    \tilde \F_1(k\sim 1)=\begin{cases}
        1 & \text{ when }k \in [a,b]\\
        0 & \text{ when }k \in (-\infty, a/2)\cup (2b,\infty)
    \end{cases}.\label{def: tilde F1}
\end{equation}
By using
\eq
e^{i(x_1q_1\pm itq_1^2-y_1q_1)}=\frac{1}{i(x_1\pm2tq_1-y_1)}\partial_{q_1}[e^{i(x_1q_1\pm itq_1^2-y_1q_1)}]
\eeq
and estimates, by estimates~\eqref{t: lem1} and~\eqref{tq: lem1}, 
\begin{align}
|\partial_{q_1}[\frac{1}{(x_1\pm 2tq_1-y_1)}]|=& |\frac{2t}{(x_1\pm 2tq_1-y_1)^2}|\lesssim \frac{1}{t^{1-\beta}|q|} (\lesssim \frac{1}{t^{1-2\beta}}),
\end{align}
\begin{align}
|\partial_{q_1}[\F_1 (t^\beta|q|>1)]|=|\frac{t^\beta q_1}{|q|} t^{l\beta} \F_1'(t^\beta|q|> 1)|\leq 2t^\beta |\F_1'(t^\beta|q|> 1)|
\end{align}
and similarly by Eq.~\eqref{def: tilde F1}, 
\begin{align}
|\partial_{q_1}^l[\F_1 (t^\beta|q|>1)]|\lesssim_l t^{l\beta}\tilde \F_1(t^\beta|q|\sim 1),\qquad l=1,2,\cdots,
\end{align}
we integrate by parts the right-hand side of Eq.~\eqref{def: O2f} for many times and each integration brings up $\frac{1}{t^{1-\beta}|q|}$ (when $|q|\geq 1$) or $\frac{\chi(|q|\leq 1)}{t^{1-2\beta}}$ (when the derivative hits $\F_1$, the support of $\F_1'$ implies $|q|\leq 1$) up to a constant, and therefore, we obtain with $\beta\in (0,1/2)$,
\begin{equation}
    \|\mathcal O_2(t)f\|_{L^2_x(\mathbb{R}^n)}\lesssim_N \frac{1}{t^N}\|f\|_{L^2_x(\mathbb{R}^n)}.\label{Oeq2}
\end{equation}
Estimates~\eqref{Oeq1} and ~\eqref{Oeq2} imply 
\begin{equation}
    \|\mathcal O(t)\|_{L^2_x(\mathbb{R}^n)\to L^2_x(\mathbb{R}^n)}\lesssim \frac{1}{t^{\sigma (1-\beta)}}.
\end{equation}
    
\end{proof}
\begin{proof}[Proof of Lemma~\ref{Lem: Pprop}] In this proof, we use notations
\begin{equation}
    \F_c=\F_c(\frac{|x|}{(t+1)^{1/2+\epsilon}}\geq 1) 
\end{equation}
and
\begin{equation}
    \F_1=\F_1(\sqrt{t+1}|P|\geq 1)\text{ or } \F_1=\F_1(\sqrt{t+1}|q|\geq 1)
\end{equation}
We start with proving estimate~\eqref{free: est: 1}. Let
\eq
A^\pm(t,s):=  P^\pm e^{\pm isH_0}\F_1\langle x\rangle^{-\sigma}.
\eeq
Break $A^\pm (t,s)$ into two pieces:
\begin{align}
    A^\pm(t,s)=&A^\pm_1(t,s)+A^\pm_2(t,s),
\end{align}
where $A^\pm_j(t,s),j=1,2,$ are given by
\begin{align}
A^\pm_1(t,s):=& P^\pm e^{\pm isH_0}\F_1 \langle x\rangle^{-\sigma}\chi(|x|> \frac{1}{10^{10}}((t+1)^{1/2+\epsilon}+s/\sqrt{t+1}))
\end{align}
and
\begin{align}
A^\pm_2(t,s):=& P^\pm e^{\pm isH_0}\F_1 \langle x\rangle^{-\sigma}\chi(|x|\leq \frac{1}{10^{10}}((t+1)^{1/2+\epsilon}+s/\sqrt{t+1})).
\end{align}
$\|A^\pm_1(t,s)\|_{L^2_x(\mathbb{R}^n)\to L^2_x(\mathbb{R}^n)}$ satisfies 
\begin{align}
\|A^\pm_1(t,s)\|_{L^2_x(\mathbb{R}^n)\to L^2_x(\mathbb{R}^n)}\leq &\| \langle x\rangle^{-\sigma} \chi(|x|\leq \frac{1}{10^{10}}((t+1)^{1/2+\epsilon}+s/\sqrt{t+1}))\|_{L^2_x(\mathbb{R}^n)\to L^2_x(\mathbb{R}^n)}\nonumber\\
\lesssim & \frac{1}{\langle (t+1)^{1/2+\epsilon}+s/\sqrt{t+1}\rangle^\sigma}.\label{Apm1}
\end{align}
For $A^\pm_2(t,s)$, choose $f\in L^2_x(\mathbb{R}^n)$ and by Fourier transform, $A^\pm_2(t,s)f$ reads 
\begin{align}
    A^\pm_2(t,s)f= & \frac{1}{(2\pi)^n} \int  e^{ix\cdot q} P^\pm(x,2q)e^{\pm i sq^2}\F_1  e^{-iy\cdot q} \nonumber\\
    & \times \langle y\rangle^{-\sigma}\chi(|y|\leq \frac{1}{10^{10}}((t+1)^{1/2+\epsilon}+s/\sqrt{t+1}))f(y)d^nyd^nq\label{Apm2ts},
\end{align}
where we have used that $P^\pm=P^\pm(x,2P)$ (see Eq.~\eqref{Ppm}). We note that the phase function in the right-hand side of Eq.~\eqref{Apm2ts} is given by:
\eq
f(q)=(x-y)\cdot q \pm s |q|^2.
\eeq
When 
\eq
|y|\leq \frac{1}{10^{10}}((t+1)^{1/2+\epsilon}+s/\sqrt{t+1}),
\eeq
\eq
|x|\geq \frac{(t+1)^{1/2+\epsilon}}{2}
\eeq
and
\eq
|q|\geq \frac{1}{2\sqrt{t+1}},
\eeq
on the support of the cutoffs $P^\pm(x,2P)$, by Eqs.~\eqref{Feq1}-~\eqref{Ppm}, we have
\begin{align}
|\nabla_q[f(q)]|=&|x-y\pm 2sq|\geq  \frac{1}{2}|x\pm 2sq|-|y|\geq  \frac{1}{10^6}(|x|+2s|q|)-|y|\nonumber\\
\geq & \frac{1}{10^6}(\frac{(t+1)^{1/2+\epsilon}}{2}+s/(2\sqrt{1+t}) )-|y|\nonumber\\
\geq & \frac{1}{10^7}(|x|+s|q| ).\label{Ap: B:est1}
\end{align}
We choose an orthogonal basis $\{e_1,\cdots,e_n\}$ with $e_1:=\frac{x-y\pm 2sq}{|x-y\pm 2sq|}$. Let $x_1:=x\cdot e_1, y_1:=x\cdot e_1$ and $q_1:=q\cdot e_1$. 
We also have the estimate
\begin{align}
    \| \F_cA^\pm_2(t,s)f\|_{L^2_x(\mathbb{R}^n)}\leq & \| \langle x\rangle^{-n}\|_{L^2_x(\mathbb{R}^n)}\| \langle x\rangle^n \F_c A^\pm_2(t,s)f\|_{ L^\infty_x(\mathbb{R}^n)}\lesssim \| \langle x\rangle^n A^\pm_2(t,s)f\|_{ L^\infty_x(\mathbb{R}^n)} .\label{2toinfty}
\end{align}
By estimate~\eqref{Ap: B:est1}, we have 
\begin{align}
|\partial_{q_1}[\frac{1}{(x_1-y_1\pm 2s q_1)}]|=& |\frac{\mp 2s}{ (x_1-y_1\pm 2sq_1)^2}|\nonumber\\
\lesssim  & \frac{s}{(|x|+s|q|)^2}\lesssim  \frac{1}{|q|(|x|+s|q|)}.
\end{align}
This, together with estimates
\begin{align}
    | \partial_{q_1}[\F_1(\sqrt{t+1}|q|\geq 1)]|= & \frac{\sqrt{t+1}|q_1|}{|q|} |\F_1'(\sqrt{t+1}|q|\geq 1)|\nonumber\\
    \lesssim  &\frac{1}{|q|}, 
\end{align}
\begin{align}
    | \partial^j_{q_1}[\F_1(\sqrt{t+1}|q|\geq 1)]|\lesssim_j  &\frac{1}{|q|^j} 
\end{align}
and (recall that $P^\pm(r,v)$ is defined in terms of $\hat{r}$ and $\hat{v}$. See Eqs.~\eqref{Prv+} and~\eqref{Prv-}.)
\begin{align}
    |\partial_{q_1}^j[P^\pm(x,2q)]|\lesssim_j \frac{1}{|q|^j},\quad j=1,2,\cdots,
\end{align}
implies,  by taking integration by parts in $q_1$ variable for $N$ times,
\begin{align}
|\langle x\rangle^n A_2^\pm(t,s)f|\lesssim & \int \frac{\langle x\rangle^n \chi(|q|\geq \frac{1}{2\sqrt{t+1}})}{ |q|^N\langle |x|+s|q|\rangle^N} \nonumber\\
&\times\langle y\rangle^{-\sigma}\chi(|y|\leq \frac{1}{10^{10}}((t+1)^{1/2+\epsilon}+s/\sqrt{t+1}))|f(y)|d^nqd^ny .\label{Apm2}
\end{align}
Taking the integral over $q$ in the right-hand side of estimate~\eqref{Apm2} and using estimates (with $|q|\geq 1/(2\sqrt{t+1})$)
\eq
\frac{1}{\langle |x|+s|q|\rangle}\lesssim \frac{1}{\langle |x|+s/(2\sqrt{t+1})\rangle}
\eeq
and
\eq
\frac{\langle x\rangle^n}{\langle |x|+s/(2\sqrt{t+1})\rangle^N}\lesssim  \frac{1}{\langle |x|+s/(2\sqrt{t+1})\rangle^{N-n}},
\eeq
we obtain 
\begin{align}
    \|\langle x\rangle^n \F_cA_2^\pm(t,s)f\|_{L^\infty_x(\mathbb{R}^n)}\lesssim &\int \frac{\langle x\rangle^n\F_c \chi(|q|\geq \frac{1}{2\sqrt{t+1}})}{ |q|^N\langle |x|+s/(2\sqrt{t+1})\rangle^N} \nonumber\\
&\times\langle y\rangle^{-\sigma}\chi(|y|\leq \frac{1}{10^{10}}((t+1)^{1/2+\epsilon}+s/\sqrt{t+1}))|f(y)|d^nqd^ny\nonumber\\
\lesssim & \int \frac{\F_c 
 (t+1)^{\frac{N-n}{2}}}{ \langle |x|+s/(2\sqrt{t+1})\rangle^{N-n}} \nonumber\\
&\times\langle y\rangle^{-\sigma}\chi(|y|\leq \frac{1}{10^{10}}((t+1)^{1/2+\epsilon}+s/\sqrt{t+1}))|f(y)|d^ny.\label{ApB: est2}
\end{align}
By using Cauchy-Schwarz inequality in estimate~\eqref{ApB: est2}, we arrive at 
\begin{align}
    & \|\langle x\rangle^n \F_cA_2^\pm(t,s)f\|_{L^\infty_x(\mathbb{R}^n)}\nonumber\\
    \lesssim & \frac{\F_c 
 (t+1)^{\frac{N-n}{2}}}{ \langle |x|+s/(2\sqrt{t+1})\rangle^{N-n}} \|f(y)\|_{L^2_y(\mathbb{R}^n)}\| \langle y\rangle^{-\sigma}\chi(|y|\leq \frac{1}{10^{10}}((t+1)^{1/2+\epsilon}+s/\sqrt{t+1}))\|_{L^2_y(\mathbb{R}^n)}\nonumber\\
 \lesssim & \frac{ 
 (t+1)^{\frac{N-n}{2}}}{ \langle (t+1)^{1/2+\epsilon}+s/(2\sqrt{t+1})\rangle^{N-3n/2}} \|f(y)\|_{L^2_y(\mathbb{R}^n)}\nonumber\\
 \lesssim& \frac{ 
 1}{ \langle (t+1)^{1/2+\epsilon}+s/(2\sqrt{t+1})\rangle^{\frac{2\epsilon N}{1+2\epsilon}-\frac{n}{2}-\frac{2\epsilon}{1+2\epsilon}n}} \|f(y)\|_{L^2_y(\mathbb{R}^n)}\nonumber\\
 \lesssim & \frac{ 
 1}{ \langle (t+1)^{1/2+\epsilon}+s/(2\sqrt{t+1})\rangle^{\sigma}} \|f(y)\|_{L^2_y(\mathbb{R}^n)}\label{estimate: A2pm}
\end{align}
with $N=[\frac{1+2\epsilon}{2\epsilon}\sigma+n+\frac{(1+2\epsilon)n}{4\epsilon}]+1$. Estimates~\eqref{Apm1},~\eqref{2toinfty} and ~\eqref{estimate: A2pm} imply estimate~\eqref{free: est: 1}. Similarly, we have estimate~\eqref{free: est: 2}. 
\end{proof}


\begin{proof}[Proof of Lemma~\ref{lem: Ppropfree}] Let $x$ and $y$ denote the position variables on the left-hand side and the right-hand side, respectively. The velocity is given by $\nabla_P[H_0]=2P$. Let $q$ denote the variable in the Fourier space. When $|x|\geq (t+1)^{1/2+\epsilon}/2$, $s|q|\leq \frac{2t}{\sqrt{t+1}}$ and $|y|\leq (t+1)^{1/2+\epsilon}/4$, 
\eq
|x-y-2tq|\geq |x|-|y|-2t|q|\geq (t+1)^{1/2+\epsilon}/20.
\eeq
Therefore, by using a similar argument of Lemma~\ref{Lem: Pprop}, we get estimate~\eqref{free: est: 3}. This completes the proof.
    
\end{proof}

\begin{proof}[Proof of Lemma~\ref{lem: Ppmf}] We fix $s\geq 0$. $P^\pm e^{\pm isH_0}f$ satisfies, for all $M\geq 1$ and $\epsilon\in (0,1/2)$, 
\begin{align}
    \| P^\pm e^{\pm isH_0}f\|_{L^2_x(\mathbb{R}^n)}\leq & \| P^\pm e^{\pm isH_0}F_1(\sqrt{s+1}|P|\geq 1)\chi(|x|\leq M) f \|_{L^2_x(\mathbb{R}^n)}\nonumber\\
    &+\| P^\pm e^{\pm isH_0}F_1(\sqrt{s+1}|P|\geq 1)\chi(|x|> M) f \|_{L^2_x(\mathbb{R}^n)}\nonumber\\
    &+\| P^\pm e^{\pm isH_0}(1-F_1(\sqrt{s+1}|P|\geq 1) )f \|_{L^2_x(\mathbb{R}^n)}\nonumber\\
    \leq & \|  P^\pm e^{\pm isH_0}F_1(\sqrt{s+1}|P|\geq 1)\langle x\rangle^{-1}\|_{L^2_x(\mathbb{R}^n)\to L^2_x(\mathbb{R}^n)}\nonumber\\
    &\times \|\langle x\rangle\chi(|x|\leq M) f \|_{L^2_x(\mathbb{R}^n)}+\| \chi(|x|>M)f\|_{L^2_x(\mathbb{R}^n)}\nonumber\\
    &+\| (1-F_1(\sqrt{s+1}|P|\geq 1) )f \|_{L^2_x(\mathbb{R}^n)}.
\end{align}
By taking $M=(1+s)^{1/100}$ and by using Lemma~\ref{Lem: Pprop}, we obtain that 
\begin{align}
    \| P^\pm e^{\pm isH_0}f\|_{L^2_x(\mathbb{R}^n)}\lesssim_\epsilon & \frac{1}{\langle s\rangle^{1/2}} (1+s)^{1/100} \| f\|_{L^2_x(\mathbb{R}^n)} +\| \chi(|x|>(s+1)^{1/100})f\|_{L^2_x(\mathbb{R}^n)}\nonumber\\
    &+\| (1-F_1(\sqrt{s+1}|P|\geq 1) )f \|_{L^2_x(\mathbb{R}^n)}\nonumber\\
    \to & 0\label{psi21}
\end{align}
as $s\to \infty.$  
    
\end{proof}

\begin{proof}[Proof of Lemma~\ref{lem: not linear: local}] Let $x$ and $y$ denote the position variables on the left-hand side and the right-hand side, respectively. The velocity is given by $\nabla_P[H_0]=2P$. We break $\chi(|x|\leq s^\alpha)P^\mp e^{\pm i sH_0}f$ into three pieces:
\begin{align}
\chi(|x|\leq s^\alpha)P^\mp e^{\pm i sH_0}f=&\chi(|x|\leq s^\alpha)P^\mp e^{\pm i sH_0}\F_1(s^{(1-\alpha)/100}|P|>1 )\chi(|x|\leq s^\alpha)f\nonumber\\
&+\chi(|x|\leq s^\alpha)P^\mp e^{\pm i sH_0}\F_1(s^{(1-\alpha)/100}|P|\leq1 )\chi(|x|\leq s^\alpha)f\nonumber\\
&+\chi(|x|\leq s^\alpha)P^\mp e^{\pm i sH_0}\chi(|x|> s^\alpha)f\nonumber\\
=: & f_1(s)+f_2(s)+f_3(s).
\end{align}
We have 
\begin{align}
    \limsup\limits_{s\to \infty}\|f_2(s)\|_{L^2_x(\mathbb{R}^n)}\leq&\limsup\limits_{s\to \infty} \| \F_1(s^{(1-\alpha)/100}|P|\leq 1) \chi(|x|\leq s^\alpha)f\|_{L^2_x(\mathbb{R}^n)}\nonumber\\
    \leq& \limsup\limits_{s\to \infty}\|\F_1(s^{(1-\alpha)/100}|P|\leq1 )f\|_{L^2_x(\mathbb{R}^n)}+\limsup\limits_{s\to \infty}\| \chi(|x|> s^\alpha)f\|_{L^2_x(\mathbb{R}^n)}\nonumber\\
    =& 0\label{free: f2}
\end{align}
and
\begin{align}
    \limsup\limits_{s\to \infty}\|f_3(s)\|_{L^2_x(\mathbb{R}^n)}\leq& \limsup\limits_{s\to \infty}\| \chi(|x|> s^\alpha)f\|_{L^2_x(\mathbb{R}^n)}\nonumber\\
    = & 0.\label{free: f3}
\end{align}
To estimate $f_1(s)$, we follow a similar argument of Lemma~\ref{Lem: Pprop}. by Fourier transform, $A^\pm_2(t,s)f$ reads 
\begin{align}
    f_1(s)= & \frac{1}{(2\pi)^n} \int  \chi(|x|\leq s^\alpha)e^{ix\cdot q} P^\pm(x,2q)e^{\pm i sq^2}\F_1(s^{(1-\alpha)/100}|q|>1)  e^{-iy\cdot q} \nonumber\\
    & \times\chi(|y|\leq s^\alpha)f(y)d^nyd^nq\label{fts},
\end{align}
where we have used that $P^\pm=P^\pm(x,2P)$ (see Eq.~\eqref{Ppm}). We note that the phase function in the right-hand side of Eq.~\eqref{fts} is given by:
\eq
f(q)=(x-y)\cdot q \pm s |q|^2.
\eeq
When $|y|\leq s^\alpha$, $|x|\geq s^\alpha$ and $|q|\geq \frac{1}{2}s^{(\alpha-1)/100}$, by Eqs.~\eqref{Feq1}-~\eqref{Ppm} and estimate
\begin{align}
2s|q|\geq s^{(99+\alpha)/100}\geq s^\alpha,\quad s\geq 1,
\end{align}
we have
\begin{align}
|\nabla_q[f(q)]|=&|x-y\pm 2sq|\geq  \frac{1}{2}|x\pm 2sq|-|y|\geq  \frac{1}{10^6}(|x|+2s|q|)-|y|\nonumber\\
\geq & \frac{1}{10^7}(|x|+s|q| ).\label{Ap: B:est10}
\end{align}
We choose an orthogonal basis $\{e_1,\cdots,e_n\}$ with $e_1:=\frac{x-y\pm 2sq}{|x-y\pm 2sq|}$. Let $x_1:=x\cdot e_1, y_1:=x\cdot e_1$ and $q_1:=q\cdot e_1$. 
We also have the estimate
\begin{align}
    \| f_1(s)\|_{L^2_x(\mathbb{R}^n)}\leq & \| \langle x\rangle^{-n}\|_{L^2_x(\mathbb{R}^n)}\| \langle x\rangle^n f_1(s)\|_{ L^\infty_x(\mathbb{R}^n)}\lesssim \| \langle x\rangle^n f_1(s)\|_{ L^\infty_x(\mathbb{R}^n)} .\label{B2toinfty}
\end{align}
By estimate~\eqref{Ap: B:est1}, we have 
\begin{align}
|\partial_{q_1}[\frac{1}{(x_1-y_1\pm 2s q_1)}]|=& |\frac{\mp 2s}{ (x_1-y_1\pm 2sq_1)^2}|\nonumber\\
\lesssim  & \frac{s}{(|x|+s|q|)^2}\lesssim  \frac{1}{|q|(|x|+s|q|)}.
\end{align}
This, together with estimates
\begin{align}
    | \partial_{q_1}[\F_1(s^{(1-\alpha)/100}|q|\geq 1)]|= & \frac{s^{(1-\alpha)/100}|q_1|}{|q|} |\F_1'(s^{(1-\alpha)/100}|q|\geq 1)|\lesssim  \frac{1}{|q|}, 
\end{align}
\begin{align}
    | \partial^j_{q_1}[\F_1(s^{(1-\alpha)/100}|q|\geq 1)]|= & \frac{s^{(1-\alpha)/100}|q_1|}{|q|} |\F_1'(s^{(1-\alpha)/100}|q|\geq 1)|\lesssim_j  \frac{1}{|q|^j}
\end{align}
and (recall that $P^\pm(r,v)$ is defined in terms of $\hat{r}$ and $\hat{v}$. See Eqs.~\eqref{Prv+} and~\eqref{Prv-}.)
\begin{align}
    |\partial_{q_1}^j[P^\pm(x,2q)]|\lesssim_j \frac{1}{|q|^j},\quad j=1,2,\cdots,
\end{align}
implies,  by taking integration by parts in $q_1$ variable for $N$ times,
\begin{align}
|\langle x\rangle^n f_1(s)|\lesssim & \int \frac{\langle x\rangle^n \chi(|q|\geq \frac{1}{2}s^{(\alpha-1)/100})}{ |q|^N\langle |x|+s|q|\rangle^N} \chi(|y|\leq s^\alpha)|f(y)|d^nqd^ny .\label{BApm2}
\end{align}
Taking the integral over $q$ in the right-hand side of estimate~\eqref{BApm2} and using estimates (with $|q|\geq \frac{1}{2}s^{(\alpha-1)/100}$)
\eq
\frac{1}{\langle |x|+s|q|\rangle}\lesssim \frac{1}{\langle |x|+s^{(99+\alpha)/100}\rangle}
\eeq
and
\eq
\frac{\langle x\rangle^n}{\langle |x|+s^{(99+\alpha)/100}\rangle^{N}}\lesssim  \frac{1}{\langle s^{(99+\alpha)/100}\rangle^{N-n}},
\eeq
we obtain 
\begin{align}
    \|\langle x\rangle^n f_1(s)\|_{L^\infty_x(\mathbb{R}^n)}\lesssim &\int \frac{\langle x\rangle^n \chi(|q|\geq \frac{1}{2}s^{(\alpha-1)/100})}{ |q|^N\langle |x|+s|q|\rangle^N} \chi(|y|\leq s^\alpha)|f(y)|d^nqd^ny\nonumber\\
\lesssim & \int \frac{ 
 s^{ (\alpha-1)(N-n)/100}}{ \langle s^{(99+\alpha)/100}\rangle^{N-n}} \chi(|y|\leq s^\alpha)|f(y)|d^ny\nonumber\\
 \lesssim & \frac{1}{\langle s\rangle^{N-n}}\int  \chi(|y|\leq s^\alpha)|f(y)|d^ny.\label{ApB: est20}
\end{align}
By using Cauchy-Schwarz inequality in estimate~\eqref{ApB: est2}, we arrive at, as $s\to \infty$, 
\begin{align}
     \|\langle x\rangle^n f_1(s)\|_{L^\infty_x(\mathbb{R}^n)}\lesssim & \frac{1}{\langle s\rangle^{N-n}} \| f(y)\|_{L^2_y(\mathbb{R}^n)}\| \chi(|y|\leq s^\alpha)\|_{L^2_y(\mathbb{R}^n)}\nonumber\\
     \lesssim &\frac{1}{\langle s\rangle^{N-n-\alpha n/2}} \| f(y)\|_{L^2_y(\mathbb{R}^n)} \nonumber\\
     \to & 0
\end{align}
with $N=n+\alpha n/2+1$. This, together with estimates~\eqref{free: f2} and~\eqref{free: f3}, implies ~\eqref{goal: lem: eq1}. 

\end{proof}


\begin{proof}[Proof of Lemma~\ref{Lem: Pprop: charge}] By equations
\begin{align}
&\F_c(\frac{|x-t\eta|}{(t+1)^{1/2+\epsilon}}\geq 1)P^\pm_{t\eta} e^{ i(s-t)H_0}\F_1(\sqrt{t+1}|2P-\eta|\geq 1) \langle x-s\eta\rangle^{-\sigma}\nonumber\\
=& e^{-it\eta\cdot P} \F_c(\frac{|x|}{(t+1)^{1/2+\epsilon}}\geq 1)P^\pm(x,2P-\eta)e^{it\eta \cdot P} e^{ i(s-t)H_0}\F_1(\sqrt{t+1}|2P-\eta|\geq 1)\nonumber\\
&\times e^{-is\eta\cdot P} \langle x\rangle^{-\sigma}e^{is\eta\cdot P}\nonumber\\
=& e^{-it\eta\cdot P} \F_c(\frac{|x|}{(t+1)^{1/2+\epsilon}}\geq 1)P^\pm(x,2P-\eta) e^{i (s-t)(H_0 -\eta\cdot P)}\F_1(\sqrt{t+1}|2P-\eta|\geq 1)\langle x\rangle^{-\sigma}e^{is\eta\cdot P}\nonumber
\end{align}
and
\begin{align}
&P^-_{t\eta} e^{ -itH_0}\F_1(\sqrt{t+1}|2P-\eta|\geq 1) \langle x\rangle^{-\sigma}\nonumber\\
=& e^{-it\eta\cdot P} P^-(x,2P-\eta)e^{it\eta \cdot P} e^{ -itH_0}\F_1(\sqrt{t+1}|2P-\eta|\geq 1)  \langle x\rangle^{-\sigma}\nonumber\\
=& e^{-it\eta\cdot P} P^-(x,2P-\eta) e^{-it(H_0 -\eta\cdot P)}\F_1(\sqrt{t+1}|2P-\eta|\geq 1)\langle x\rangle^{-\sigma}\nonumber,
\end{align}
we obtain
\begin{align}
    &\| \F_c(\frac{|x-t\eta|}{(t+1)^{1/2+\epsilon}}\geq 1)P^\pm_{t\eta} e^{ i(s-t)H_0}\F_1(\sqrt{t+1}|2P-\eta|\geq 1) \langle x-s\eta\rangle^{-\sigma}\|_{L^2_x(\mathbb{R}^n)\to L^2_x(\mathbb{R}^n)}\nonumber\\
    =& \| \F_c(\frac{|x|}{(t+1)^{1/2+\epsilon}}\geq 1)P^\pm(x,2P-\eta) e^{i(s-t)(H_0 -\eta\cdot P)}\F_1(\sqrt{t+1}|2P-\eta|\geq 1)\langle x\rangle^{-\sigma}\|_{L^2_x(\mathbb{R}^n)\to L^2_x(\mathbb{R}^n)}\label{free: charge: eq1}
\end{align}
and
\begin{align}
&\|P^-_{t\eta} e^{ -itH_0}\F_1(\sqrt{t+1}|2P-\eta|\geq 1) \langle x\rangle^{-\sigma}\|_{L^2_x(\mathbb{R}^n)\to L^2_x(\mathbb{R}^n)}\nonumber\\
=& \| P^-(x,2P-\eta) e^{-it(H_0 -\eta\cdot P)}\F_1(\sqrt{t+1}|2P-\eta|\geq 1)\langle x\rangle^{-\sigma}\|_{L^2_x(\mathbb{R}^n)\to L^2_x(\mathbb{R}^n)}.\label{free: charge: eq2}
\end{align}
By a similar argument of Lemma~\ref{Lem: Pprop}, we obtain, with $u\geq0$, 
\begin{align}
&\| \F_c(\frac{|x|}{(t+1)^{1/2+\epsilon}}\geq 1)P^\pm (x,2P-\eta)e^{\pm iu(H_0-\eta\cdot P)}\F_1(\sqrt{t+1}|2P-\eta|\geq 1) \langle x\rangle^{-\sigma}\|_{L^2_x(\mathbb{R}^n)\to L^2_x(\mathbb{R}^n)}\nonumber\\
\lesssim_\epsilon &\frac{1}{\langle (t+1)^{1/2+\epsilon}+u/\sqrt{t+1}\rangle^\sigma}\label{free: est: 1: charge}
\end{align}
and
\eq
\| P^-(x,2P-\eta) e^{- iu(H_0-\eta\cdot P)} \F_1((u+1)^{1/2-\epsilon}|2P-\eta|\geq 1)\langle x\rangle^{-\sigma}\|_{L^2_x(\mathbb{R}^n)\to L^2_x(\mathbb{R}^n)}\lesssim_\epsilon \frac{1}{\langle u\rangle^{\sigma/2}}\label{free: est: 2: charge}.
\eeq
Estimates~\eqref{free: est: 1: charge} and~\eqref{free: est: 2: charge}, together with Eqs.~\eqref{free: charge: eq1} and ~\eqref{free: charge: eq2}, imply estimates~\eqref{free: charge: goal1} and~\eqref{free: charge: goal2}. This completes the proof.
    
\end{proof}


\begin{proof}[Proof of Lemma~\ref{lem: Ppropfree: charge}] By equation
\begin{align}
    & \F_c(\frac{|x-t\eta|}{(t+1)^{1/2+\epsilon}}\geq 1)P_{t\eta}^+ e^{ -i(t-s)H_0}\F_1(\sqrt{t+1}|2P-\eta|< 1) \langle x-s\eta\rangle^{-\sigma}\nonumber\\
    =& e^{-it\eta\cdot P} \F_c(\frac{|x|}{(t+1)^{1/2+\epsilon}}\geq 1)P^+(x,2P-\eta) e^{-i(t-s)(H_0-\eta\cdot P)}\F_1(\sqrt{t+1}|2P-\eta|< 1)\langle x\rangle^{-\sigma}e^{is\eta\cdot P},
\end{align}
we obtain 
\begin{align}
&\| \F_c(\frac{|x-t\eta|}{(t+1)^{1/2+\epsilon}}\geq 1)P_{t\eta}^+ e^{ -i(t-s)H_0}\F_1(\sqrt{t+1}|2P-\eta|< 1) \langle x-s\eta\rangle^{-\sigma}\|_{L^2_x(\mathbb{R}^n)\to L^2_x(\mathbb{R}^n)}\nonumber\\
=&\| \F_c(\frac{|x|}{(t+1)^{1/2+\epsilon}}\geq 1)P^+(x,2P-\eta) e^{-i(t-s)(H_0-\eta\cdot P)}\F_1(\sqrt{t+1}|2P-\eta|< 1)\langle x\rangle^{-\sigma}\|_{L^2_x(\mathbb{R}^n)\to L^2_x(\mathbb{R}^n)}.\label{free: charge: eq5}
\end{align}
By a similar argument of Lemma~\ref{lem: Ppropfree}, we obtain, with $s\in [0,t]$, 
\begin{align}
&\| \F_c(\frac{|x|}{(t+1)^{1/2+\epsilon}}\geq 1)P^+(x,2P-\eta) e^{-i(t-s)(H_0-\eta\cdot P)}\F_1(\sqrt{t+1}|2P-\eta|< 1)\langle x\rangle^{-\sigma}\|_{L^2_x(\mathbb{R}^n)\to L^2_x(\mathbb{R}^n)}\nonumber\\
\lesssim_\epsilon& \frac{1}{(t+1)^{\frac{1}{2}\sigma+\epsilon\sigma}}.\label{free: charge: eq4}
\end{align}
Estimate~\eqref{free: charge: eq4} and Eq.~\eqref{free: charge: eq5} imply estimate~\eqref{free: charge: eq6}. This completes the proof.
    
\end{proof}

\begin{proof}[Proof of Lemma~\ref{lem: Ppmf: charge}] $P^-_{t\eta} e^{- itH_0}f$ satisfies, for all $M\geq 1$ and $\epsilon\in (0,1/2)$, with $t\geq 0$,
\begin{align}
    \| P^-_{t\eta} e^{- itH_0}f\|_{L^2_x(\mathbb{R}^n)}\leq & \| P^-_{t\eta} e^{- itH_0}F_1(\sqrt{t+1}|2P-\eta|\geq 1)\chi(|x|\leq M) \psi(0) \|_{L^2_x(\mathbb{R}^n)}\nonumber\\
    &+\| P^-_{t\eta} e^{- itH_0}F_1(\sqrt{t+1}|2P-\eta|\geq 1)\chi(|x|> M) \psi(0) \|_{L^2_x(\mathbb{R}^n)}\nonumber\\
    &+\| P^-_{t\eta} e^{- itH_0}(1-F_1(\sqrt{t+1}|2P-\eta|\geq 1) )\psi(0) \|_{L^2_x(\mathbb{R}^n)}\nonumber\\
    \leq & \|  P^-_{t\eta} e^{- itH_0}F_1(\sqrt{t+1}|2P-\eta|\geq 1)\langle x\rangle^{-1}\|_{L^2_x(\mathbb{R}^n)\to L^2_x(\mathbb{R}^n)}\nonumber\\
    &\times \|\langle x\rangle\chi(|x|\leq M) \psi(0) \|_{L^2_x(\mathbb{R}^n)}+\| \chi(|x|>M)\psi(0)\|_{L^2_x(\mathbb{R}^n)}\nonumber\\
    &+\| (1-F_1(\sqrt{t+1}|2P-\eta|\geq 1) )\psi(0) \|_{L^2_x(\mathbb{R}^n)}.
\end{align}
By taking $M=(1+t)^{1/100}$ and by using Lemma~\ref{Lem: Pprop: charge}, we obtain that 
\begin{align}
    \| P^-_{t\eta} e^{-itH_0}f\|_{L^2_x(\mathbb{R}^n)}\lesssim_\epsilon & \frac{1}{\langle t\rangle^{1/2}} (1+t)^{1/100} \| \psi(0)\|_{L^2_x(\mathbb{R}^n)} +\| \chi(|x|>(t+1)^{1/100})\psi(0)\|_{L^2_x(\mathbb{R}^n)}\nonumber\\
    &+\| (1-F_1(\sqrt{t+1}|2P-\eta|\geq 1) )\psi(0) \|_{L^2_x(\mathbb{R}^n)}\nonumber\\
    \to & 0\label{psi21: charge}
\end{align}
as $t\to \infty.$ This completes the proof. 
    
\end{proof}


\begin{proof}[Proof of Lemma~\ref{lem: not linear: charge}]By equation
\begin{align}
    \chi(|x-t\eta|\leq t^\alpha) P_{t\eta}^\pm e^{- i t H_0}f=& e^{-it\eta\cdot P} \chi(|x|\leq t^\alpha) P^-(x,2P-\eta) e^{it\eta\cdot P}e^{-itH_0}f \nonumber\\
    =& e^{-it\eta\cdot P} \chi(|x|\leq t^\alpha) P^-(x,2P-\eta)e^{-it(H_0-\eta\cdot P)}f,
\end{align}
we obtain
\begin{equation}
   \|\chi(|x-t\eta|\leq t^\alpha) P_{t\eta}^- e^{- i t H_0}f \|_{L^2_x(\mathbb{R}^n)}=\|\chi(|x|\leq t^\alpha) P^-(x,2P-\eta)e^{-it(H_0-\eta\cdot P)}f \|_{L^2_x(\mathbb{R}^n)}.
\end{equation}
This, together with Lemma~\ref{Lem: Pprop: charge}, implies
\begin{align}
    &\|\chi(|x|\leq t^\alpha) P^-(x,2P-\eta)e^{-it(H_0-\eta\cdot P)}f \|_{L^2_x(\mathbb{R}^n)}\nonumber\\
    \leq & \| P^-(x,2P-\eta)e^{-it(H_0-\eta\cdot P)}\F_1(\sqrt{t+1}|2P-\eta|\geq 1) \chi(|x|<t^{1/10})f \|_{L^2_x(\mathbb{R}^n)}\nonumber\\
    &+\| P^-(x,2P-\eta)e^{-it(H_0-\eta\cdot P)}\F_1(\sqrt{t+1}|2P-\eta|\geq 1) \chi(|x|\geq t^{1/10})f \|_{L^2_x(\mathbb{R}^n)}\nonumber\\
    & +\| \F_1(\sqrt{t+1}|2P-\eta|< 1) f \|_{L^2_x(\mathbb{R}^n)}\nonumber\\
    \lesssim& \frac{1}{\langle t\rangle}\|\langle x\rangle^2  \chi(|x|<t^{1/10})f \|_{L^2_x(\mathbb{R}^n)}+ \|\chi(|x|\geq t^{1/10})f \|_{L^2_x(\mathbb{R}^n)}+\| \F_1(\sqrt{t+1}|2P-\eta|< 1) f \|_{L^2_x(\mathbb{R}^n)}\nonumber\\
    \lesssim& \frac{1}{\langle t\rangle^{4/5}}\|f \|_{L^2_x(\mathbb{R}^n)}+ \|\chi(|x|\geq t^{1/10})f \|_{L^2_x(\mathbb{R}^n)}+\| \F_1(\sqrt{t+1}|2P-\eta|< 1) f \|_{L^2_x(\mathbb{R}^n)}\nonumber\\
    \to & 0
\end{align}
as $t\to \infty$. This leads to 
\begin{equation}
  \limsup\limits_{t\to \infty}  \|\chi(|x-t\eta|\leq t^\alpha) P_{t\eta}^- e^{- i t H_0}f \|_{L^2_x(\mathbb{R}^n)}=0\label{charge: limit: goal2},
\end{equation}
which, together with 
\begin{equation}
  \lim\limits_{t\to \infty}  \|\chi(|x-t\eta|\leq t^\alpha) e^{- i t H_0}f \|_{L^2_x(\mathbb{R}^n)}=0,
\end{equation}
implies 
\begin{equation}
  \limsup\limits_{t\to \infty}  \|\chi(|x-t\eta|\leq t^\alpha) P_{t\eta}^+ e^{- i t H_0}f \|_{L^2_x(\mathbb{R}^n)}=0\label{charge: limit: goal3}.
\end{equation}
Both Eqs. ~\eqref{charge: limit: goal2} and ~\eqref{charge: limit: goal3} imply Eq.~\eqref{charge: limit: goal1}. We finish the proof.

\end{proof}
\begin{lemma}\label{Lem: B1} Let $\{e_1,\cdots,e_n\}$ denote an orthogonal basis in $\mathbb{R}^n$. For all $\sigma>0$, we have 
\begin{equation}
    \| [\langle x\rangle^\sigma, \frac{P_j}{\langle P\rangle}]\langle x\rangle^{-\sigma}\|_{L^2_x(\mathbb{R}^n)\to L^2_x(\mathbb{R}^n)}\lesssim 1,\qquad j=1,\cdots,n.\label{B1: goal eq}
\end{equation}
\end{lemma}
 \begin{proof} {It follows directly from basic pseudodifferential calculus.}

\end{proof}
    

\begin{lemma}\label{Lem: B2} For all $\delta>0$, $\langle x\rangle^\delta f\in \mathcal H^1_x$, we have 
\begin{equation}
    \| \langle x\rangle^\delta \langle P\rangle f \|_{L^2_x(\mathbb{R}^n)}\lesssim \| \langle x\rangle^\delta f\|_{\mathcal H^1_x}.\label{Lem:B2:eq}
\end{equation}
    
\end{lemma}
\begin{proof} It follows straightforwardly from pseudodifferential calculus.
\end{proof}

\section{The sketch of Proof of Proposition~\ref{Abs}}\label{app:3}

The proof of Proposition~\ref{Abs} is the same as the proof of the original existence proof of the Free Channel Wave Operator. The use of Cook's argument reduces as before the problem to proving the integrability in time of the contribution of the Interaction term.
We use the dispersive estimate for $U_0(t)$ similar to the one one we used for the free flow: for some $p\in (2,\infty]$,
$$
\|U_0(-t)\mathcal{N}(x,t,\psi(t))\psi(t)\|_{L^p_x(\mathbb{R}^n)}\leq \frac{1}{t^{1+\epsilon}},\qquad t\geq 1.
$$
Then the integrability of the interaction term follows from, with $1/q+1/p=1/2$,
\begin{align}
&\|\F_c(|x|/t^{\alpha}\leq 1)U_0(-t)\mathcal{N}(x,t,\psi(t))\psi(t)\|_{L^2_x(\mathbb{R}^n)}\nonumber \\
\leq &\|\F_c(|x|/t^{\alpha}\leq 1)\|_{L^q_x(\mathbb{R}^n)} \|U_0(-t)\mathcal{N}(x,t,\psi(t))\psi(t)\|_{L^p_x(\mathbb{R}^n)}\nonumber\\
\lesssim&  \frac{1}{t^{1+\epsilon-n\alpha/q}} \in L^1_t[1,\infty)
\end{align}
provided that
\begin{align}
   \alpha< \frac{\epsilon q}{n}=\frac{\epsilon}{n(1/2-1/p)}=\frac{2\epsilon p}{n(p-2)}.
\end{align}
This, together with estimate~\eqref{Fcnneg} and propagation estimate, implies the existence of free channel wave operator in part $(1)$ of Proposition~\ref{can2}.

 When the dimension is low, and the thresholds of $U_0(t)$ results in too low dispersive decay rate, or when the decay requires a smooth initial data (as is the case for the wave equation) we add to the definition of the Free Channel Wave operator cut-off functions that vanish in a shrinking in time neighborhood of the thresholds, and also cut off a neighborhood of infinite frequency. Then, the above estimates hold as well.
 This frequency cutoff does NOT change the wave operator. This is because on the complement, the wave operator is zero (by taking the weak limit in t that defines the operator.) This holds provided the solution is uniformly bounded in $\mathcal{H}^1_x.$
 In part $(2)$ of Proposition~\ref{can2}, we cut off a neighborhood of infinite frequency. The interaction term satisfies the estimate, with $\frac{1}{\tilde q}+\frac{1}{\tilde p}=\frac{1}{2}$, 
 \begin{align}
&\|\F_c(|x|/t^{\alpha}\leq 1)\F_1(|P|\leq t^\beta)U_0(-t)\mathcal{N}(x,t,\psi(t))\psi(t)\|_{L^2_x(\mathbb{R}^n)}\nonumber \\
\leq &\|\F_c(|x|/t^{\alpha}\leq 1)\|_{L^{\tilde q}_x(\mathbb{R}^n)} \|\F_1(|P|\leq t^\beta)U_0(-t)\mathcal{N}(x,t,\psi(t))\psi(t)\|_{L^{\tilde p}_x(\mathbb{R}^n)}\nonumber\\
\lesssim&  \frac{1}{t^{1+\epsilon-n\alpha/{\tilde q}-\beta k}} \in L^1_t[1,\infty)
\end{align}
provided that
\begin{align}
    \epsilon-n\alpha/{\tilde q}-\beta k>0.\label{check: last}
\end{align}
Inequality~\eqref{check: last} is satisfied since when $\alpha\in (0, \min\{\frac{2 \epsilon \tilde p}{(\tilde p-2)n}, \frac{\epsilon \tilde p}{n}\}) $ and 
 $\beta\in (0, \min\{\alpha,\frac{\epsilon \tilde p-n\tilde p\alpha(1/2-1/\tilde p)}{\tilde p k} \}) $,
\begin{align}
\epsilon \tilde p -n\alpha>0,
\end{align}
\begin{align}
\epsilon \tilde p-n\tilde p\alpha(1/2-1/\tilde p)=\epsilon \tilde p -n \alpha (\tilde p-2)/2=\epsilon \tilde p -n\tilde p\alpha/\tilde q>0
\end{align}
and
\begin{align}
    \epsilon-n\alpha/{\tilde q}-\beta k>& \epsilon-n\alpha/{\tilde q}- (\epsilon \tilde p -n\alpha)/\tilde p \nonumber\\
=& 0.
\end{align}
 Additionally, $\beta<\alpha$ and both space and frequency cut-off functions satisfy the non-negativity property: Estimates~\eqref{Fcnneg} and 
 \begin{equation}
     \partial_{t}[\F_1(|q|\leq t^\beta)]=-\beta t^{-1-\beta}|q|\F_1'(|q|\leq t^\beta)\geq 0
 \end{equation}
 are valid. Then by using propagation estimates twice, we obtain the existence of the free channel wave operator.

\bibliographystyle{uncrt}
\def\bibfont{\HUGE}

\end{document}